\documentclass[english,compress,square,super]{elsarticle}

\usepackage{lineno}
\modulolinenumbers[1]
\usepackage{geometry}
\usepackage[colorlinks=true]{hyperref}

\journal{arXiv.}

\usepackage{amsmath}
\usepackage{xcolor}
\newcommand\tb[1]{\boldsymbol{#1}}
\newcommand\td{\mathrm{d}}
\newcommand\pd{\partial}
\usepackage{mathrsfs}

\newcommand\ddfrac[2]{{\displaystyle\frac{\displaystyle #1}{\displaystyle #2}}}
\usepackage{stmaryrd}
\usepackage{enumitem}
\usepackage{caption}
\usepackage{subcaption}
\usepackage{cleveref}
\usepackage{algorithm}  
\usepackage{algpseudocode}  
\usepackage{amssymb}
\usepackage{upgreek}
\usepackage{pgfplots}
\pgfplotsset{compat=1.16}
\usepackage{svg}
\usepackage{multirow}
\usepackage{booktabs}
\usepackage{mathtools}
\usepackage{adjustbox}
\usepackage[para,online,flushleft]{threeparttable}
\renewcommand{\arraystretch}{1.2}
\makeatletter
\usepackage{xspace}
\algnewcommand{\algorithmicgoto}{\textbf{go to}}%
\algnewcommand{\Goto}{\algorithmicgoto\xspace}%
\DeclareRobustCommand\bigop[1]{%
  \mathop{\vphantom{\sum}\mathpalette\bigop@{#1}}\slimits@
}
\newcommand{\bigop@}[2]{%
  \vcenter{%
    \sbox\z@{$#1\sum$}%
    \hbox{\resizebox{\ifx#1\displaystyle.9\fi\dimexpr\ht\z@+\dp\z@}{!}{$\m@th#2$}}%
  }%
}

\usepackage{amsthm}
\newtheorem{remark}{Remark}
\usepackage{accents}

\begin{document}

\begin{frontmatter}

\title{Mixed material point method formulation, stabilization, and validation for a unified analysis of free-surface and seepage flow}

\author[1]{Bodhinanda {Chandra}}
\ead{bchandra@berkeley.edu}
\author[2]{Ryota {Hashimoto}\corref{cor1}}
\ead{hashimoto.ryota.6e@kyoto-u.ac.jp}
\cortext[cor1]{Corresponding authors}
\author[4]{Ken {Kamrin}}
\ead{kkamrin@mit.edu}
\author[1]{Kenichi {Soga}\corref{cor1}}
\ead{soga@berkeley.edu}

\address[1]{Department of Civil and Environmental Engineering, University of California, Berkeley, CA, 94720, USA}
\address[2]{Department of Civil and Earth Resources Engineering, Kyoto University, Kyoto, Japan}
\address[4]{Department of Mechanical Engineering, Massachusetts Institute of Technology, Cambridge, MA, 02139, USA}

\begin{abstract}
This paper presents a novel stabilized mixed material point method (MPM) designed for the unified modeling of free-surface and seepage flow. The unified formulation integrates the Navier-Stokes equation with the Darcy-Brinkman-Forchheimer equation, effectively capturing flows in both non-porous and porous domains. In contrast to the conventional Eulerian computational fluid dynamics (CFD) solver, which solves the velocity and pressure fields as unknown variables, the proposed method employs a monolithic displacement-pressure formulation adopted from the mixed-form updated-Lagrangian finite element method (FEM). To satisfy the discrete \textit{inf-sup} stability condition, a stabilization strategy based on the variational multiscale method (VMS) is derived and integrated into the proposed formulation. Another distinctive feature is the implementation of blurred interfaces, which facilitate a seamless and stable transition of flows between free and porous domains, as well as across two distinct porous media. The efficacy of the proposed formulation is verified and validated through several benchmark cases in 1D, 2D, and 3D scenarios. Conducted numerical examples demonstrate enhanced accuracy and stability compared to analytical, experimental, and other numerical solutions.
\end{abstract}


\end{frontmatter}

\section{Introduction}

Coupled porous-medium and incompressible free-flow systems are prevalent in many natural and industrial settings. Examples include interactions between surface and subsurface water, geothermal and hydrocarbon flow in drilled wells, industrial air filters, as well as the filtration of blood through arterial walls. From a macro-scale perspective, the theoretical framework governing the motion of free fluid and porous-media flow has been distinctively developed into two sets of equations. The incompressible free flow is typically described by the incompressible Navier-Stokes equation. Meanwhile, the Darcy \cite{darcy1856fontaines}, Brinkman \cite{brinkman1949calculation}, and Forchheimer \cite{forchheimer1901wasserbewegung} equations have played a central role in investigating seepage flow in saturated porous media \cite{ehlers2020darcy}. The study of the interaction between these two flow regimes was pioneered by \citet{beavers1967boundary}, who described flow conditions across and along the porous interface. Since then, extensive investigations have been conducted through numerous experimental observations \cite{neale1974practical, vafai1987analysis, ochoa1995momentum_i, ochoa1995momentum_ii, ochoa1998, alazmi2001analysis}. While the two flows inherently form unified systems, they are often separated for logical reasons, such as time scales, and technical considerations, including analytical and computational solvability \cite{furman2008modeling}. A substantial body of literature exists, covering the theoretical derivation and numerical implementation of various interface coupling conditions, e.g.~\cite{saffman1971boundary, mikelic2000interface, discacciati2002mathematical, discacciati2009navier, jager2009modeling, lacis2017framework, eggenweiler2020unsuitability}.

Numerical analysis of coupled free-fluid and porous media flow is typically conducted using Eulerian mesh-based CFD methods, such as FEM \cite{badia2009unified} and the finite volume method (FVM) \cite{angot2017asymptotic}. In these methods, necessary boundary treatments are imposed at the porous interface to connect the two governing equations. Two classes of boundary treatment exist in the literature: \textit{sharp} and \textit{blurred interface}\cite{goyeau2003momentum, valdes2013velocity}, each with their respective advantages and drawbacks. The sharp interface approach generally requires parametrization of the interface, specifying its geometry and exact location. This can be achieved through conforming (\textit{body-fitted}) or nonconforming discretization, employing explicit or implicit geometric models. However, for complex geometries in 3D, the surface representation and discretization are challenging to automate and error-prone due to issues like \textit{small-cut instabilities}, making them computationally demanding. On the other hand, the blurred interface method presents itself as an appealing option by abandoning the idea of precisely defining the interface. Instead, it implicitly approximates the domain using a phase-field function that gradually transitions from a nonzero value inside the domain to zero in the exterior \cite{stoter2017diffuse}. Because of this feature, the blurred interface is often more computationally convenient and stable (due to its improved smoothness), though with a price of decreasing accuracy as it introduces an extra length scale into the formulation \cite{rycroft2020reference}.  

In many engineering fields, such as geotechnical and coastal engineering, the handling of evolving free-surface flow poses another computational challenge. Over the past decades, the volume-of-fluid (VOF) method \cite{hirt1981volume} has been developed within the Eulerian framework to address the tracking of free surfaces. This framework has undergone further developments, leading to the introduction of a unified equation that incorporates porosity and drag forces into the balance equations, where the Brinkman correction is also considered \cite{liu1999numerical, hieu2006verification, del2012three, jensen2014investigations, higuera2014three, losada2016modeling}. A similar strategy has also been adopted to FEM using the edge-based level-set method by \citet{larese2015finite}. At the same time, an alternative approach for modeling free surface flows using the \textit{meshless} methods has also been gaining attention \cite{lucy1977numerical, gingold1977smoothed, koshizuka1996moving}. Akbari \cite{akbari2013moving, akbari2014modified} pioneered the incorporation of the unified framework into the smoothed particle hydrodynamics (SPH) to study the interaction of waves and porous structures. Many similar types of research using the meshless method have followed thereafter \cite{aly2015three, peng2017multiphase, khayyer2018development, sun2019numerical}.

As methodologies for resolving the aforementioned coupled fluid dynamics problems continue to evolve, recent advances in hybrid Eulerian-Lagrangian methods have made noteworthy contributions to the field. Among various approaches, the material point method (MPM) \cite{Sulsky1994, Sulsky1995b}, derived from the particle-in-cell (PIC) method \cite{harlow1964particle}, has gained popularity in simulating large-deformation geophysical flows. The MPM harnesses both the strengths of mesh-based and meshless descriptions by discretizing the material domain into Lagrangian moving particles while simultaneously utilizing a constantly reset background grid to solve the balance equations and apply boundary conditions. This attribute enhances MPM's versatility to simulate various challenging problems in engineering and computer graphics \cite{jiang2016material, zhang2016material, soga2016trends, wolper2021glacier, liang2022shear}.

Research efforts in modeling flow within porous media using MPM have gained traction since the late 2000s due to the increasing interest in multi-phase coupled problems \cite{zhang2008material, mackenzie2010modeling}. Following that, numerous hydro-mechanical formulations have been proposed, which can be categorized into two approaches based on the use of material points to represent solid and fluid phases \cite{soga2016trends, solowski2021material}: the \textit{one-layer} \cite{zhang2009material, zabala2011progressive, zhao2020stabilized, kularathna2021semi} and \textit{double-layer} \cite{bandara2015coupling, baumgarten2019general, yamaguchi2020solid} formulations. The one-layer formulation focuses on simulating the deformation of porous solids given the presence of pore fluid. It does not properly track the motion of the fluid phase, hence, fluid mass conservation is only ensured locally through the continuity equation. In contrast, the double-layer formulation fully tracks the Lagrangian motion of fluid material points, thus ensuring global mass conservation. However, the presence of two material point layers increases computational costs and necessitates careful modeling of interfaces. 

The double-layer MPM formulation facilitates seamless coupling of free-surface and seepage flow, in addition to computing the deformation of the solid phase. The standard formulation proposed by \citet{abe2014material} and \citet{bandara2015coupling} is rooted in the theory of porous media \cite{coussy2004poromechanics}, assuming Terzaghi's effective stress \cite{terzaghi1951theoretical} and Darcy's drag model. \citet{martinelli2016soil} then incorporated the Forchheimer term into their formulation and proposed modeling the blurred porous interface using a linear function. The formulation is further extended by \citet{baumgarten2019general} to include the packing-fraction-dependent Brinkman correction term \cite{einstein1906calculation} and the nonlinear drag force model suggested by \citet{beetstra2007drag}. While the majority of prior formulations explicitly estimate fluid pressure through an equation of state, following the weakly-compressible assumption, e.g.~\cite{molinos2023derivation}, \citet{yamaguchi2020solid} adopted the fractional-step approach \cite{kularathna2017implicit, zhang2017incompressible}, which solves pressure implicitly via a pressure Poisson equation \cite{chorin1968numerical}. This significantly mitigates pressure oscillations due to \textit{volumetric locking}, although its stability appears to be highly influenced by the application of pressure boundary conditions \cite{chandra2023stabilized} and the selected time step size, which can be particularly problematic for low-permeability porous media \cite{kularathna2021semi, mieremet2016numerical, morikawa2022soil}.

While MPM has demonstrated effectiveness in modeling various solid materials subjected to substantial deformations, MPM is noted for encountering notable integration errors as deformations persist, particularly evident in scenarios involving fluid flow \cite{baumgarten2021coupled, baumgarten2023analysis}. These integration errors contribute to an accelerated accumulation of overall errors during simulations, consequently limiting MPM's applicability in addressing long-period fluid-related problems \cite{baumgarten2021coupled}. To avoid such issues, recent trends in coupled MPM simulations often involve the use of Eulerian mesh-based methods, such as FVM \cite{baumgarten2021coupled} and FEM \cite{pan2021mpm, tran2022mpmice}, equipped with a free surface tracking algorithm. In our recent work \cite{chandra2023stabilized}, we presented several strategies to enhance the stability and accuracy of free-surface fluid simulation in MPM using the stabilized mixed displacement-pressure formulation.

The primary objective of this study is to extend the stabilized mixed MPM formulation previously proposed for incompressible fluid dynamics \cite{chandra2023stabilized} to incorporate the concurrent flow of free-surface fluid and fluid within porous media. In this work, we initially develop a unified formulation of coupled Navier-Stokes and Darcy-Brinkman-Forchheimer equations and derive consistent matrix systems within the previously converged updated Lagrangian framework. These discrete equations are then implemented within the MPM formulation and solved iteratively using the Newton-Krylov solver. To enhance computational efficiency and reduce costs, we opt for equal-order elements and propose a stabilization strategy based on the variational multiscale method (VMS) \cite{hughes1995multiscale, hughes1998variational, codina2001stabilized, codina2002stabilized, codina2007time} to satisfy the discrete inf–sup condition. 

Furthermore, to characterize the non-porous and porous domains, we construct porosity and permeability fields using the blurred interface approach. As we enter the porous domain, the resisting force (drag) is gradually applied through empirical closures based on the Darcy-Forchheimer drag model, following the particle-based \textit{smooth volume-fraction} approach suggested by \citet{mackenzie2010modeling}. A gradual porosity-dependent transition is also considered for the fluid deviatoric stress, particularly for cases with changing effective viscosity. For a comprehensive characterization of the spatiotemporal errors of the fluid solver, we assume the porous material to be rigid and isotropic. Here, the normal and tangential fluid tractions at the porous interface are assumed to be in equilibrium \cite{vafai1987analysis}, and the effect of capillary force at the moving front inside the porous media is neglected. In this study, we also conducted in-depth verification and validation tests to assess the efficacy of the proposed formulation and stabilization. Several steady and unsteady analytical solutions for 1D gravity-driven flow under various conditions are derived and validation benchmarks in 2D and 3D are performed. The obtained numerical results are then compared with experimental data and results from other numerical methods.

The current manuscript is organized as follows. First, the governing equations of the unified Navier--Stokes/Darcy--Brinkman--Forchheimer flow are presented in \Cref{sec:stab_mixed_MPM}, along with the spatiotemporal discretization and stabilization strategy in mixed MPM. \Cref{sec:num_ex} presents several numerical examples involving complex flows in porous and non-porous regions to assess the accuracy, stability, and efficiency of the proposed formulation. Finally, \Cref{sec:conclusions} presents the conclusions and suggestions for future work.

\section{Stabilized mixed MPM formulation for unified free-surface--seepage flow}
\label{sec:stab_mixed_MPM}

In this section, we introduce the stabilized mixed MPM formulation for a unified analysis of free-surface--seepage flow. Subsequently, we establish the notation to be used throughout. Three types of sub- or superscripts are defined: $\square_p$ for particle variables, $\square_I$ or $\square_J$ for nodal variables, and $\square^t$ for time indices where $t=n$ and $t=n+1$ denote current and future variables, respectively. Additionally, the notation for operators is defined as follows: $\dot{\square}$ and $\ddot{\square}$ represent first and second-order material time derivatives, $\square\otimes\square$ denotes the dyadic operator, and $\square\cdot\square$ and $\square:\square$ denote single and double contractions of tensor indices.

\subsection{Geometrical settings and discretization assumptions}

\begin{figure}[h!]
    \centering
    \includegraphics[scale=0.48]{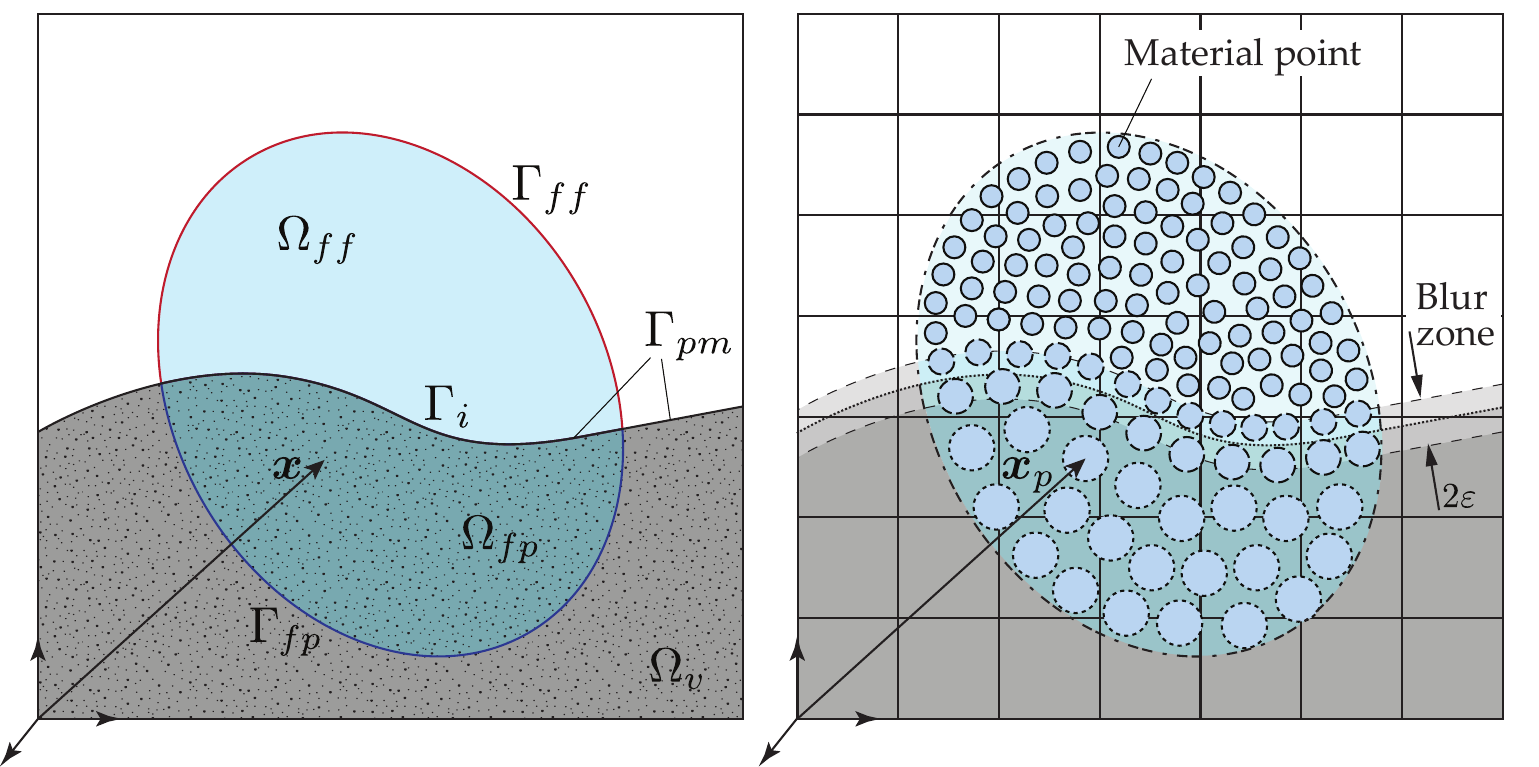}
    \caption{(left) Geometrical settings of the problem domain and (right) a schematic illustration of their discretization in MPM.}
    \label{fig:2_porous_flow_MPM}
\end{figure}

Let us examine a fluid body $\mathcal{B}_f$, which occupies a volume $\Omega_f$ within the three-dimensional Euclidean space $\mathcal{R}^{3}$. The volume $\Omega_f$ consists of free-fluid flow region $\Omega_{ff}\subset\mathcal{R}^{3}$ and porous flow region $\Omega_{fp}\subset\mathcal{R}^{3}$, where $\Omega_f = \Omega_{ff}\cup \Omega_{fp}$ and $\Omega_{ff} \cap \Omega_{fp} =  \varnothing$ (see Fig.~\ref{fig:2_porous_flow_MPM} (left)). These two volumes are assumed to be closed with smooth boundaries, $\pd \Omega_{ff}$ and $\pd \Omega_{fp}$ corresponding to $ \Omega_{ff}$ and $ \Omega_{fp}$. Here, three surfaces can be defined accordingly: the free-fluid-porous interface $\Gamma_i:= \pd \Omega_{ff} \cap \pd \Omega_{fp}$, the free-fluid free surface $\Gamma_{ff} := \pd \Omega_{ff} \backslash \Gamma_i$, and the pore-fluid free surface $\Gamma_{fp} := \pd \Omega_{fp} \backslash \Gamma_{i}$. The fluid is assumed to be present in a single phase, where when it occupies the porous medium, it is assumed to fully saturate the pore (or void) space $\Omega_v$.

Meanwhile, the porous medium is assumed to have undergone homogenization, where the pore space can be characterized at the macroscopic level by porosity $\theta$ defined as:
\begin{eqnarray}
    \theta := \frac{\Omega_v}{\Omega}\,,
\end{eqnarray}
where $\Omega$ is defined as the total control volume. In the most general case, the porosity may vary with space and time, i.e.~$\theta=\theta(\tb x,t)$. However, since the porous medium is considered to be static and rigid, we only consider the variation of porosity in space, i.e.~$\theta=\theta(\tb x)$. In fact, it is possible to use the porosity to characterize the entire domain, where $\theta(\tb x)=1.0,$ at the non-porous region and $0.0<\theta(\tb x)<1.0$ at the porous domain. The boundary between these two regions is denoted as $\Gamma_{pm}$, where at a given time instance $\Gamma_{pm}\supset \Gamma_i$.

In the MPM \citep{Sulsky1994, Sulsky1995b}, the initial geometry of fluid body $\mathcal{B}_f$ is discretized into $N_p$ material points (see Fig.~\ref{fig:2_porous_flow_MPM} (right)), or particles, of volume $V_p$ such that:
\begin{eqnarray}
\Omega_f = \int_\Omega \theta \td \Omega \approx \sum_{p=1}^{N_p}\theta_p V_p \,.
\label{eq:mpm_quadrature}
\end{eqnarray}
Throughout the simulation, these material points carry all relevant information, such as masses, velocities, stresses, and history-dependent properties. The initial mass of each material point then can be defined as:
\begin{eqnarray}
    m_p = \theta^{t_0}_{p} \rho_p V^{t_0}_{p}\,,
\end{eqnarray}
where $\rho$ denotes the intrinsic fluid mass density, $V^{t_0}$ is the initial volume, and $\theta^{t_0}\equiv\theta(\tb x^{t_0})$ indicates the initial porosity, which can be evaluated according to the initial material point position $\tb x^{t_0}_{p}$. Over time, the fluid material points may move across the porous and non-porous boundary $\Gamma_{pm}$, where the occupied volume correspondingly expands or shrinks according to the change in porosity. In this work, the porous medium interface $\Gamma_{pm}$ is modeled as a blurred interface, where the porosity gradually changes over a certain thickness of size $2\varepsilon$. Further implementation details on the blurred interface will be given in \Cref{subsec:porosity_field}. As depicted in Fig.~\ref{fig:2_porous_flow_MPM}, in addition to the material point discretization, a \textit{background grid} is also generated and utilized to solve the discretized balance equations, as well as to approximate field variables and their gradient.

\subsection{Governing equations and weak formulation}
\label{subsec:balance_laws}

The mass and momentum balance of the fluid flow in both the porous and non-porous regions can be expressed in the Lagrangian form following the theory of porous media \cite{bandara2015coupling, baumgarten2019general}, where the porosity field $\theta$ is utilized to unify the free-surface and the seepage flow. First, the mass balance can be written as:
\begin{eqnarray}
\frac{\pd \bar{\rho}}{\pd t} + \nabla \cdot \left( \bar{\rho} \tb v \right) = 0 \,, \label{eq:fluid_mass_balance}
\end{eqnarray}
where $\tb v$ and $\bar{\rho}$ indicate the fluid velocity field and the effective mass density, i.e.~$\bar{\rho}=\theta \rho$. Following the assumption of incompressible fluids, the intrinsic fluid density is assumed to be constant, i.e.~${\td \rho}/{\td t} = 0$. Therefore, since the variation of porosity in time can be neglected, the continuity equation can be derived: 
\begin{eqnarray}
    \nabla \cdot (\theta \tb v) = 0 \,.
    \label{eq:continuity}
\end{eqnarray}

Subsequently, the linear momentum balance is expressed as:
\begin{eqnarray}
	\bar{\rho} \dot{\tb{v}} = \nabla \cdot \tb \sigma + \bar{\rho}\tb{b} + \tb{f}_b + \tb{f}_d \,,
	\label{eq:mpm_fluid_momentum_balance}
\end{eqnarray}
where $\tb \sigma$ denotes the Cauchy stress and $\tb b$ is the body force vectors per unit mass. The Cauchy stress tensor can be separated into hydrostatic and deviatoric components as follows:
\begin{eqnarray}
    \tb \sigma &=& -\theta p \tb I + \tb s \,,
    \label{eq:cauchy_stress}\\
    \tb s &=& \mu \left( \nabla \tb v + \left(\nabla \tb v \right)^T \right) = 2 \mu \dot{\tb \varepsilon}\,,
    \label{eq:molecular_shear_stress}
\end{eqnarray}
where $p$ is the fluid pressure and $\tb{s}$ is for the deviatoric stress. In the above expression, $p$ is assumed to be positive in compression, while $\tb{\sigma}$ and $\tb{s}$ are considered positive in tension. Additionally, $\tb{I}$, $\mu$, and $\dot{\tb{\varepsilon}}$ represent the identity tensor, the dynamic viscosity, and the strain rate tensor, respectively.

In the current work, Brinkman's correction \cite{brinkman1949calculation} is considered, where the fluid viscous force is non-zero inside the porous media\footnote{This assumption differs from the classical Darcy's law \cite{darcy1856fontaines} and Terzaghi's effective stress principle \cite{terzaghi1951theoretical}, which both assume negligible or zero fluid shear stress within porous media.}. The viscous shear stress $\tb s$ previously given by Eq.~\eqref{eq:molecular_shear_stress} can be slightly modified to:
\begin{eqnarray}
    \tb s = \mu_e \left( \nabla \tb v + \left(\nabla \tb v \right)^T \right) = 2 \mu_e \dot{\tb \varepsilon}\,,
    \label{eq:brinkman_shear_stress}
\end{eqnarray}
where $\mu_e$ represents the effective viscosity within a porous medium. Previous studies, including direct numerical simulations \cite{martys1994computer} and experimental investigations \cite{givler1994determination}, highlight the significance of distinguishing the effective viscosity $\mu_e$ and the fluid viscosity $\mu$. In this work, in coherence with the blurred interface method, we consider a linear porosity-dependent viscosity formulation expressed in the form of:
\begin{eqnarray}
    \widetilde{\mu}:= \widetilde{\mu}(\theta) = \left(1+\eta\,(1-\theta) \right)\mu\,,
    \label{eq:porosity_dependent_visc}
\end{eqnarray}
which mathematical form is borrowed from the rheology of suspension flow \cite{einstein1906calculation, stickel2005fluid}. Here, $\eta$ is a precomputed slope given by:
\begin{eqnarray}
    \eta = \frac{\mu_e - \mu}{\mu(1-\theta_{\mathrm{ref}})}\,,
\end{eqnarray}
where $\theta_{\mathrm{ref}}$ denotes the reference porosity, ensuring $\widetilde{\mu}=\mu_e$ when $\theta=\theta_{\mathrm{ref}}$ and $\widetilde{\mu}=\mu$ when $\theta=1$. 

Meanwhile, $\tb{f}_b$ and $\tb{f}_d$ indicate the buoyant force and the viscous drag force, respectively. The former generally takes the form given by \citet{drumheller2000theories},
\begin{eqnarray}
	\tb{f}_b = p \nabla \theta\,,
	\label{eq:buoyancy_force}
\end{eqnarray}
whereas the drag force is assumed to follow the Darcy-Forchheimer equation\cite{forchheimer1901wasserbewegung} to account for the additional drop of the hydraulic head observed at high-velocity porous flow, i.e.
\begin{eqnarray}
	\tb{f}_d = - \theta^2 \mu {\tb k}^{-1} \tb{v} - \beta_F \theta^3 \rho ||\tb v|| \tb v\,.
	\label{eq:drag_force}
\end{eqnarray}
Here, $\beta_F$ [$L^{-1}$] is the non-Darcy flow coefficient which contains the information on the structure and tortuosity of the porous media. Furthermore, $\tb k \equiv k_{ij}$ $[L^2]$ is the intrinsic permeability tensor, which can be simplified to a scalar $k$ for isotropic conditions; this will be considered in the current work. 

The Kozeny-Carman equation \cite{carman1937fluid} is used to compute the permeability of the porous media:
\begin{eqnarray}
k =  \dfrac{\theta^3}{\left( 1- \theta \right)^2} \kappa\,.
\label{eq:kozeny-carman}
\end{eqnarray}
Here, $\kappa$ is the reference permeability, which can be determined experimentally by considering grain size distribution, shape, and roughness. According to \citet{ergun1952fluid}, $\kappa$ and $\beta_F$ can be computed as,
\begin{eqnarray}
    \kappa = \frac{d^2}{A}\,, \qquad \beta_F = \frac{B}{\sqrt{A}\,\theta^{3/2}\sqrt{k}}\,,
    \label{eq:ergun}
\end{eqnarray}
where $d$ is the effective grain size parameter, while $A$ and $B$ are dimensionless constants that can be set according to the flow type or type of porous materials. \citet{ergun1952fluid} utilizes $A=150$ and $B=1.75$, whereas setting $B=0$ recovers the Darcian flow as the nonlinear term vanishes. As suggested by many previous works \cite{liu1999numerical, del2012three, akbari2013moving}, $A$ and $B$ often require further calibrations for specific porous materials. Substituting Eq.~\eqref{eq:ergun}$_2$ into Eq.~\eqref{eq:drag_force}, we can rewrite the drag force equation after certain mathematical rearrangements as:
\begin{eqnarray}
	\tb{f}_d = - \theta \rho \Biggl( \underbrace{\frac{\theta \mu}{\rho k}}_{\widetilde{A}} + \underbrace{\frac{B \sqrt{\theta}}{\sqrt{Ak}}}_{\widetilde{B}} ||\tb v|| \Biggl) \tb v\,.
	\label{eq:drag_force_simplifed}
\end{eqnarray}
Here, to simplify the notation, two coefficients $\widetilde{A}$ and $\widetilde{B}$ are introduced as given in Eq.~\eqref{eq:drag_force_simplifed}. Other equivalent expressions to compute $\widetilde{A}$ and $\widetilde{B}$ according to different possible input parameters are given in \ref{app:drag_force_coeff}.

Having Eqs.~\eqref{eq:cauchy_stress}, \eqref{eq:buoyancy_force}, and \eqref{eq:drag_force_simplifed} in hand, the momentum balance given in Eq.~\eqref{eq:mpm_fluid_momentum_balance} can be rewritten as:
\begin{eqnarray}
	\theta \rho \left( \dot{\tb{v}} - \tb b\right) + \theta \nabla p - \nabla \cdot \tb s + \theta \rho \left( \widetilde{A} + \widetilde{B} ||\tb v|| \right) \tb v = 0 \,.
	\label{eq:equilibrium}
\end{eqnarray}
Here, notice that if the porosity is $\theta=1$, which characterizes the free fluid region, the mass and momentum balance expressed in Eq.~\eqref{eq:continuity} and \eqref{eq:equilibrium} will recover the incompressible Navier-Stokes equation. Furthermore, when the inertia term is omitted, Eq.~\eqref{eq:equilibrium} leads to the Darcy-Brinkman-Forchheimer equation \cite{forchheimer1901wasserbewegung, brinkman1949calculation}. 

The balance equations described above are to be numerically solved within domain $\Omega \subseteq \mathcal{R}^{3}$ over the time interval $t \in [0,T]$, by taking into account the following Dirichlet and Neumann boundary conditions (BCs):
\begin{eqnarray}
\tb{u} = \overline{\tb{u}} &\qquad& \mathrm{on} \, \Gamma_D \,, \label{eq:dir_BC} \\
\tb{\sigma}\cdot\tb{n} = \overline{\tb{t}} &\qquad&  \mathrm{on} \, \Gamma_N \,, \label{eq:neu_BC}
\end{eqnarray}
where $\tb{n}$ is the outward unit normal, while $\overline{\tb{u}}$ and $\overline{\tb{t}}$ are the corresponding Dirichlet and Neumann boundary displacement and traction, respectively. Meanwhile, $\tb{u}$ represents the displacement field, which can be related to the velocity and acceleration fields using material time derivatives, specifically $\tb{v} = \dot{\tb{u}}$ and $\tb{a} = \dot{\tb{v}} = \ddot{\tb{u}}$.

By performing the standard Galerkin approximation, the weak form of the momentum balance equation can be derived as:
\begin{eqnarray}
\begin{split}
R^{\mathrm{mom}}(\tb u, p,  \delta \tb u) =& 
\int_{\Omega}\theta \rho \left(\ddot{\tb u} - \tb b\right) \cdot \delta \tb u\, \td{\rm \Omega}
-\int_{\Omega} \theta p \left( \nabla \cdot \delta \tb u \right)\, \td\Omega
+\int_{\Omega} \tb s \colon \nabla \delta \tb u\, \td\Omega 
-\int_{\Omega} p \nabla\theta \cdot \delta \tb u \, \td\Omega\\
& \quad -\int_{\Gamma_{N}}\overline{\tb t} \cdot \delta  \tb u\, \td\Gamma + \int_{\Omega} \theta \rho \left( \widetilde{A} + \widetilde{B}||\dot{\tb u}||
\right) \dot{\tb u} \cdot \delta \tb u\, \td \Omega =0
\,.
\label{eq:weak_form_momentum}
\end{split}
\end{eqnarray}
Here, $\delta \tb u$ is an arbitrary test function, which is lying on the solution space of $\tb u$, denoted as $\mathcal{V}\subset \mathcal{R}^{3}$, such that $\delta \tb u = \lbrace \delta \tb u \in \mathcal{V}\,|\,\delta \tb u = \tb 0$ on $\Gamma_D \rbrace$. Similarly, the weak form of the mass balance can be derived by performing an $L_2$ inner product with a scalar test function $\delta p\in \mathcal{Q}$, where $\mathcal{Q}\subset \mathcal{R}$ is the space of virtual pressure as:
\begin{eqnarray}
    R^{\mathrm{mass}}(\tb u, \delta p) = 
    \int_{\Omega} \left( \theta \nabla \cdot \dot{\tb u} + \dot{\tb u} \cdot \nabla \theta \right) \, \delta p \, \td{\rm \Omega} = 0\,.
\label{eq:weak_form_mass}
\end{eqnarray}

\subsection{Spatial discretization}
\label{subsec:spatial_disc}

The standard FE interpolation is utilized to compute the displacement and pressure, along with their respective test functions, at each material point $p$. This can be expressed as:
\begin{subequations}
\begin{eqnarray}
\tb u_{p}:=\tb u (\tb x_p)=\sum_{I=1}^{n_n}S_I(\tb x_p) \tb u_{I}\,, \qquad 
\delta \tb u_{p}:=\delta \tb u (\tb x_p)=\sum_{I=1}^{n_n}S_I(\tb x_p) \delta  \tb u_{I}\,,
\label{eq:shape function_u}\\
p_{p}:=p (\tb x_p)=\sum_{I=1}^{n_n}S_I(\tb x_p) p_{I}\,, \qquad 
\delta p_{p}:=\delta p (\tb x_p)=\sum_{I=1}^{n_n} S_I(\tb x_p) \delta p_{I}\,,
\label{eq:shape function_p}
\end{eqnarray}
\label{eq:fe_interpolation}%
\end{subequations}
where $S_I(\tb x_p)$ denotes the generic notation of the nodal basis function. In the current work, equal order basis functions for $\tb u$ and $p$ are considered, where the notation $S_{Ip}\equiv S_I(\tb x^n_p)$ will be used onward for brevity. In the expressions above, $n_n$ represents the total number of the computational node (or \textit{connectivity})  linked with particle $p$, which includes the nodes of neighboring cells when utilizing higher-order basis functions. 

In this study, the quadratic B-Spline basis function introduced by \citet{steffen2008analysis} is utilized and integrated with the weighted-least-square (WLS) kernel correction method proposed by \citet{nakamura2023taylor} to address interpolation issues near the domain boundaries. It has been recently discovered, however, that the WLS kernel correction method does not guarantee non-negative values of basis functions. This condition is especially problematic for MPM, as it may lead to nodes with negative mass \cite{andersen2007material}. To address this issue, we utilize an iterative kernel correction strategy, where detailed information can be found in \ref{app:basis_funct_it_kc}. 

\subsection{Particle-grid transfer}
\label{subsec:particle_grid_transfer}

In the MPM, given that no variables are saved at the nodes, data transfer from particles to background nodes should be performed at the onset of each time step. Additionally, at the end of every time step, the kinematic variables and pressure at the nodes are interpolated back to the material points, after which the mesh is restored to its original configuration. In the current work, the transfers from particle to the grid (P2G) and from the grid to particles (G2P) are done following the FLIP (fluid-implicit-particle)\cite{brackbill1986flip}, APIC (affine-particle-in-cell)\cite{jiang2015affine}, and TPIC (Taylor-particle-in-cell)\cite{nakamura2023taylor} approach for mixed MPM as discussed in our previous work\cite{chandra2023stabilized}. For the sake of brevity, the explanation of the P2G and G2P transfers is skipped in this paper. Interested readers are recommended to refer to \cite{chandra2023stabilized} for more details.

\subsection{Modeling porosity variation with blurred interface method}
\label{subsec:porosity_field}

\begin{figure}[h!]
    \centering
     \begin{subfigure}[b]{0.3\textwidth}
         \centering
         \includegraphics[width=\textwidth]{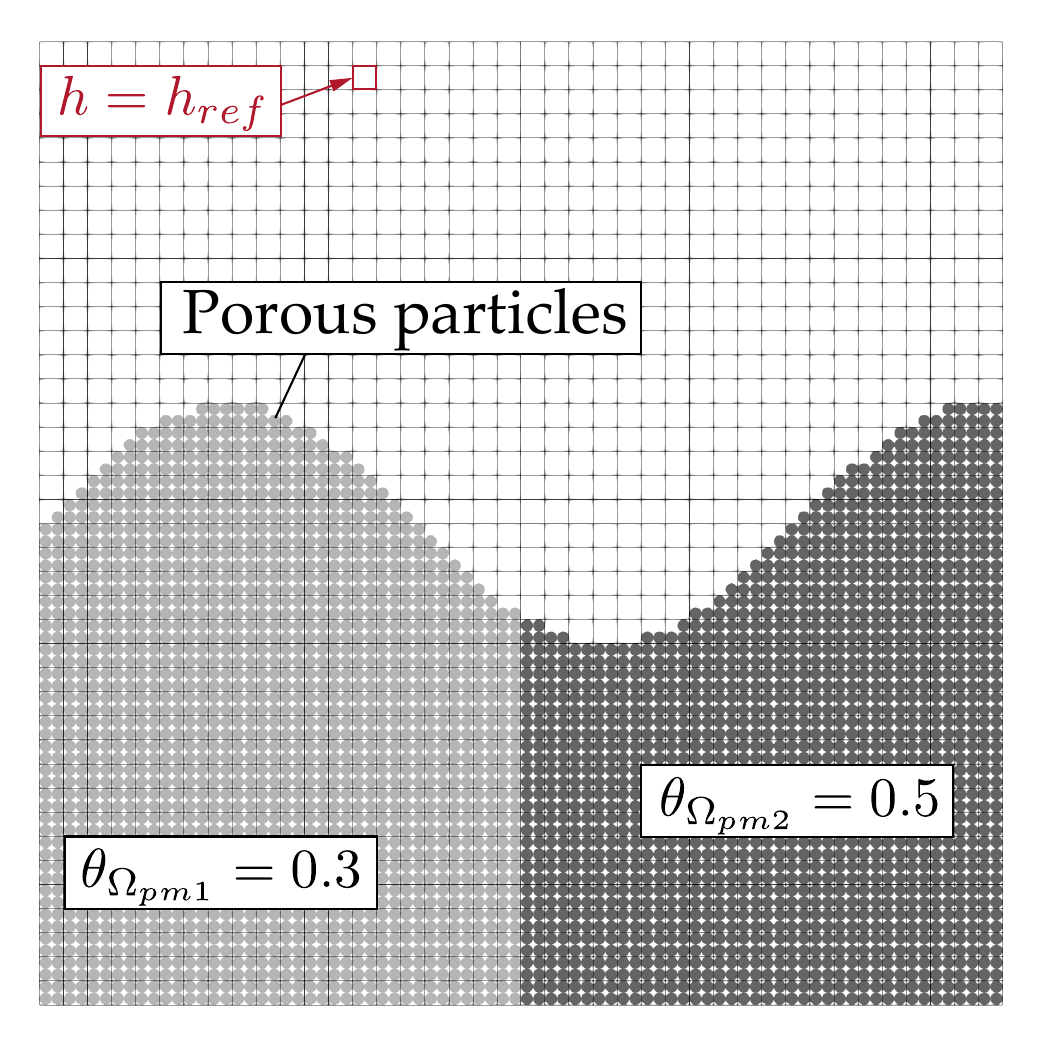}
         \caption{}
         \label{fig:2_mesh_particle}
     \end{subfigure}
     \begin{subfigure}[b]{0.3\textwidth}
         \centering
         \includegraphics[width=\textwidth]{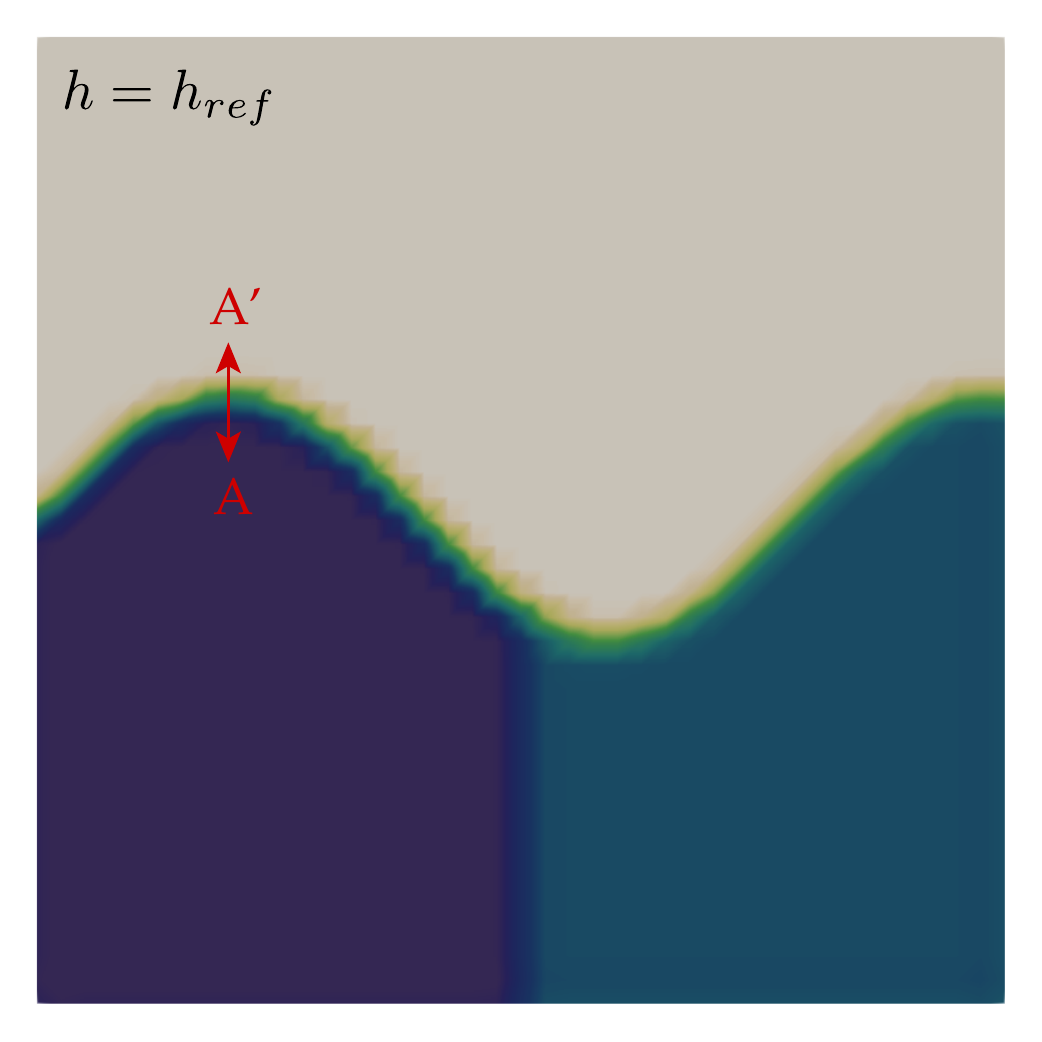}
         \caption{}
         \label{fig:2_diffuse_interface_coarse}
     \end{subfigure}
     \begin{subfigure}[b]{0.36\textwidth}
         \centering
         \includegraphics[width=\textwidth]{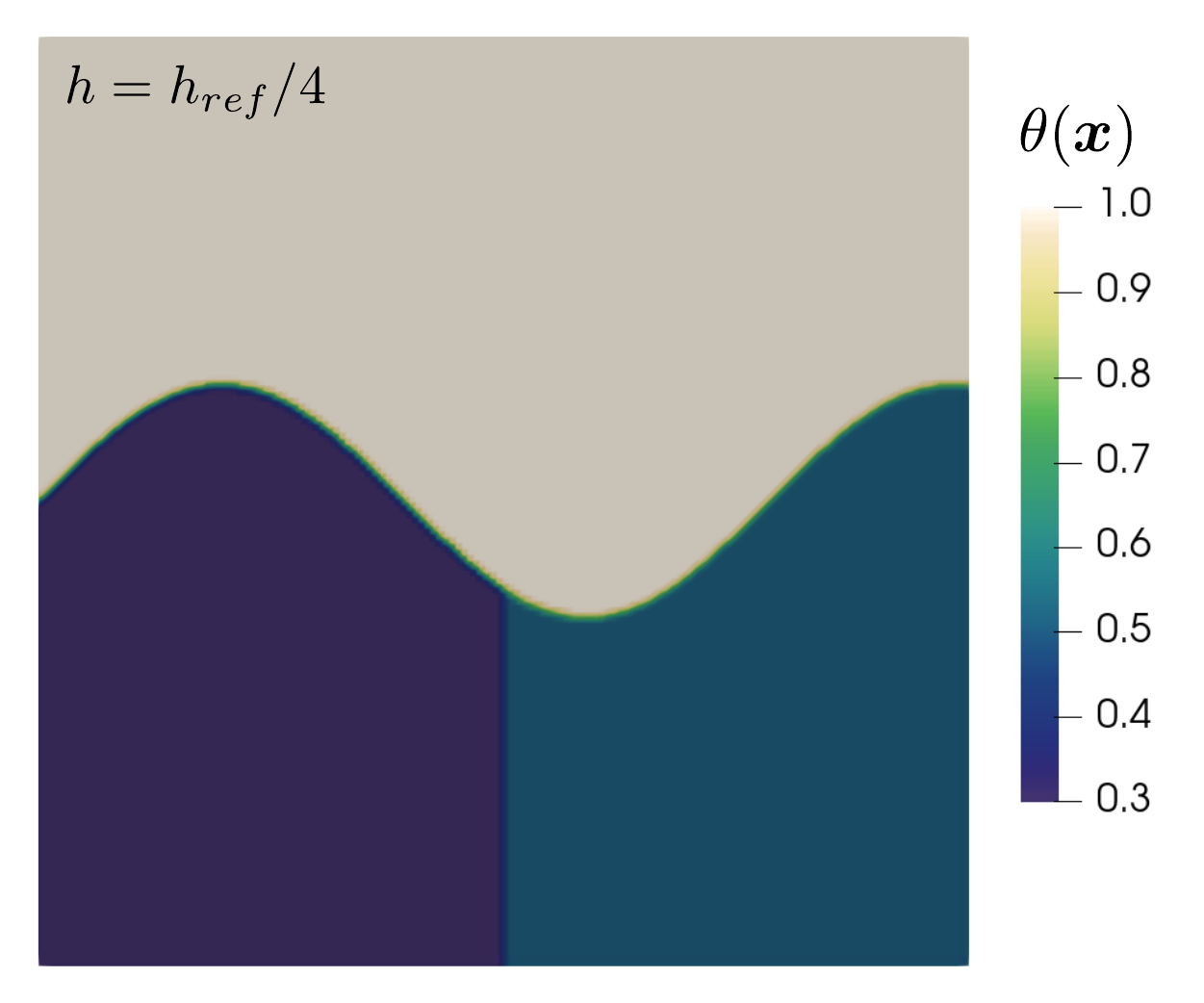}
         \caption{}
         \label{fig:2_diffuse_interface_fine}
     \end{subfigure}\\
     \vspace{0.3cm}
     \begin{subfigure}[b]{0.6\textwidth}
         \centering
         \includegraphics[width=\textwidth]{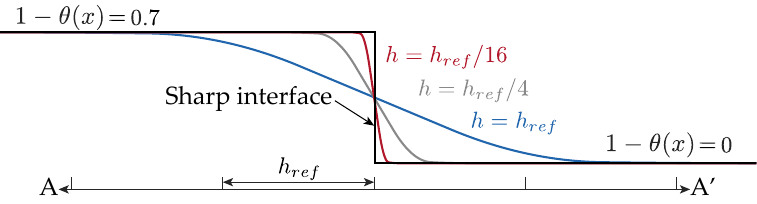}
         \caption{}
         \label{fig:2_blurred_interface}
     \end{subfigure}\\
     \vspace{0.3cm}
     \begin{subfigure}[b]{0.6\textwidth}
         \centering
         \includegraphics[width=\textwidth]{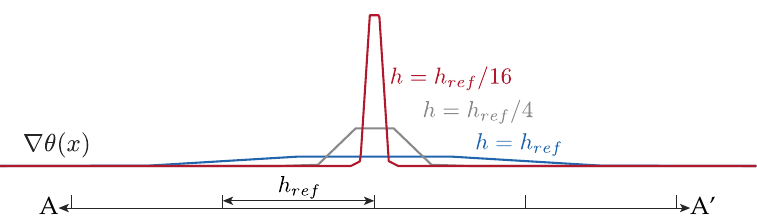}
         \caption{}
         \label{fig:2_blurred_interface_grad}
     \end{subfigure}
    \caption{Example of a blurred porous interface with two different porosity values, $\theta=0.3$ and 0.5 with the shape of a sine function: (a) discretization into porous particles and background grid with a cell size of $h_{\mathrm{ref}}=0.2$. The reconstructed nodal porosity fields in different cell sizes are plotted in (b) $h=h_{\mathrm{ref}}$ and (c) $h=h_{\mathrm{ref}}/4$. The one-dimensional porosity profile and its gradient over sections A-A' are blown up in (d) and (e), respectively, for different cell sizes.}
    \label{fig:2_diffuse_interace}
\end{figure}

In the current work, the blurred interface assumption is considered to model the porosity jump between the porous and non-porous domains and between two different porous domains. This approach has been previously considered in many prior two-phase double-point MPM works, e.g.~\cite{zhang2008material, martinelli2016soil, yamaguchi2020solid}, such that the numerical instabilities with respect to the sudden jump in porosity can be avoided. Following their approach, the porous domain, which is considered to be static and rigid, is discretized into $N_{pp}$ number of porous particles with a certain porosity value $\theta=\theta_{\Omega_{pm\alpha}}$ (see Fig.~\ref{fig:2_mesh_particle}). Here, the subscript $\alpha$ indicates the distinct porous media index with possibly different porosity and permeability values. The discretization of the porous media into porous particles is a natural choice in the MPM, particularly since this assumption can be seamlessly extended into double-point formulations \cite{bandara2015coupling, yamaguchi2020solid}, where the porous solid may also move and deform.

At the beginning of the simulation, we can construct a porosity field represented by the nodal scalar value $\theta_I$ by performing a P2G mapping. Here, we follow the approach proposed by \citet{zhang2008material} to generate a gradual profile of porosity, i.e.:  
\begin{eqnarray}
    \theta_I = \frac{1}{V_I}\sum_{pp=1}^{n_{pp}}S_{I}(\tb x_{pp}) \theta_{pp} V_{pp}\,,
    \label{eq:nodal_porosity}
\end{eqnarray}
where $V_I$ is the volume domain represented by each node $I$ which can be computed by performing an analytical integration or by Gauss quadrature as:
\begin{eqnarray}
    \quad V_I = \int S_{I}(\tb x)\, \td \Omega\,.
\end{eqnarray}
In the notation above $n_{pp}$ is the total number of static porous particles associated with node $I$. 

It is worth mentioning that, even though it is possible to model the blur zone by using distinct interpolation functions, with possibly higher orders of continuity, e.g.~\cite{stoter2017diffuse, rycroft2020reference}, this work assumes an isoparametric interpolation function for simplicity. Since a $C^1$-continuous basis function is utilized, the generated porosity gradient profile is in linear order (see Fig.~\ref{fig:2_blurred_interface_grad}). Furthermore, the width of the blur zone is $\varepsilon=1.5h$ as dictated by the kernel size of the function. Unlike the sharp interface approach, this approach introduces discretization errors around the interface $\Gamma_{i}$, which should vanish along with mesh refinement. The constructed blurred interface is depicted in 2D and 1D through Figs.~\ref{fig:2_diffuse_interface_coarse}-\ref{fig:2_blurred_interface} for different cell sizes.

Once the porosity field is constructed, at every time step, the porosity at each fluid material point $p$ can be interpolated following the standard FEM interpolation as:
\begin{eqnarray}
    \theta^n_p &:=& \theta(\tb x_p^n) = \sum_{I=1}^{n_{n}}S_{Ip} \theta_I\,,\\
    \nabla \theta^n_p &:=& \nabla \theta(\tb x_p^n) = \sum_{I=1}^{n_{n}}\nabla S_{Ip} \theta_I\,.
\end{eqnarray}
Upon traveling across domains with different porosity fields, the mass of each material point, $m_p$, should be kept constant, and thus, its volume should expand or contract correspondingly. The material point's effective density and occupied volume can be scaled by using the interpolated porosity as \cite{bandara2015coupling}:
\begin{eqnarray}
    \bar{\rho}_p^n=\theta^n_p \rho_p\,, \qquad V^n_p = \theta_p^{t_0} V_p^{t_0} / \theta^n_p\,.
\end{eqnarray}

The interpolated particle porosity will then be used to numerically solve the momentum balance equation which governs flows in both free-fluid and porous domains. In order to account for porous media with different permeability values, a similar mapping process is performed to blur the transition of the permeability parameters previously given in Eq.~\eqref{eq:drag_force_simplifed}, e.g.~the intrinsic permeability $k$, as well as the constant $A$ and $B$, if necessary. The permeability parameters are to be mapped according to:
\begin{eqnarray}
    \tb \Pi_I = \frac{1}{V_I}\sum_{pp=1}^{n_{pp}}S_{I}(\tb x_{pp}) V_{pp} \tb \Pi_{pp}\,,\qquad \mathrm{where,} \quad \tb \Pi = \{k, A, B\}\,.
\end{eqnarray}
A similar strategy can be formulated for other permeability variables, such as the hydraulic conductivity $K\,[LT^{-1}]$ or the effective grain size $d\,[L]$, to smoothly transition those parameters across different materials.

\begin{remark}
\label{rem:smoothing_type}
As pointed out by \citet{yamaguchi2020solid}, the mapping method presented in Eq.~\eqref{eq:nodal_porosity} is preferable to the typical Shepard interpolation method, e.g.~as used by \citet{bandara2015coupling}, 
\begin{eqnarray}
    \theta_I = \frac{\sum_{pp=1}^{n_{pp}} S_I(\tb x_{pp}) \theta_{pp} V_{pp}}{\sum_{pp=1}^{n_{pp}} S_I(\tb x_{pp}) V_{pp}}\,,
    \label{eq:nodal_porosity_NG}
\end{eqnarray}
for the purpose of improving accuracy and stability. As illustrated in Fig.~\ref{fig:2_blurred_interface_alternative}, the commonly used Shepard interpolation fails to produce a porosity profile that aligns with the shape of the porous interface. Furthermore, it extends the porous regions, resulting in an increased number of areas where drag forces are applied; this may lead to excessive energy dissipation. In contrast, the profile produced by Eq.~\eqref{eq:nodal_porosity} is centered at the porous interface. This configuration ensures that the domains subject to higher and lower applied drag forces are equivalent, effectively canceling out the overall numerical errors.
    
\begin{figure}[h!]
    \centering
    \includegraphics[width=0.8\textwidth]{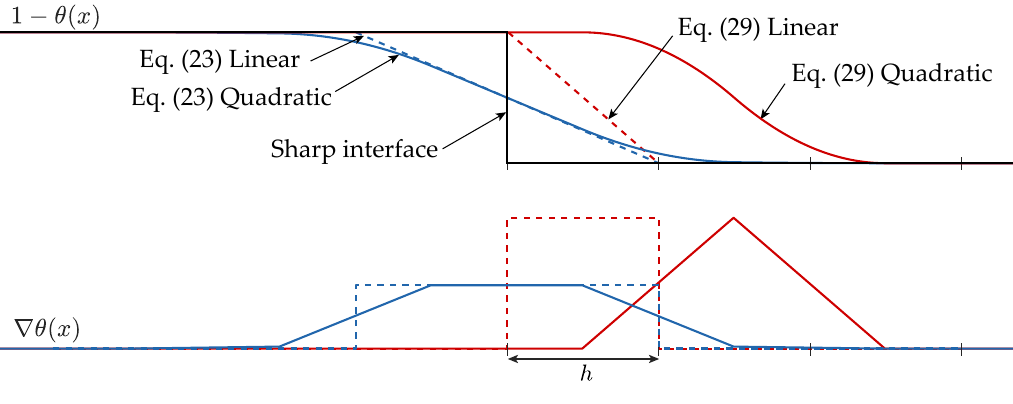}
    \caption{One-dimensional porosity profile and its gradient generated from two different formulations: Eqs.~\eqref{eq:nodal_porosity} and \eqref{eq:nodal_porosity_NG}, assuming linear and quadratic B-Spline basis functions.}
    \label{fig:2_blurred_interface_alternative}
\end{figure}
\end{remark}

\subsection{Discretized form and solution strategy}
\label{subsec:discrete_form_sol_strategy}

Given the spatial discretization described in \Cref{subsec:spatial_disc} and the MPM quadrature approximation provided in Eq.~\eqref{eq:mpm_quadrature}, and by invoking the arbitrariness of the nodal test functions, $\delta \tb u_{I}$ and $\delta p_{I}$, the discretized versions of the momentum and mass balances, Eqs.~\eqref{eq:weak_form_momentum} and \eqref{eq:weak_form_mass}, at the nodal level can be formulated as follows:
\begin{linenomath}
\begin{align}
\begin{split}
\tb R^{\mathrm{mom}}_I =&  m^n_I \Ddot{\tb{u}}_I
+
\sum_{p=1}^{n_{p}} \left(\tb B_{Ip}^T \tb{s}_p - 
\theta_p \rho_p  \tb{b} S_{Ip}\right) V_p - \sum_{J=1}^{n_n} \left(\sum_{p=1}^{n_{p}} \left(  \theta_p \nabla S_{Ip} S_{Jp} + \nabla \theta_p S_{Ip} S_{Jp}\right) V_p \right) p_J \\
& \quad -\int_{\Gamma_{N}}\overline{\tb t} S_{I}(\tb x_{\Gamma_N})\, \td\Gamma  + \Biggl( \sum_{p=1}^{n_{p}} \left( \widetilde{A}_p + \widetilde{B}_p \left\Vert \dot{\tb u}_p \right \Vert
\right) S_{Ip} m_p\Biggl) \dot{\tb u}_I
\,,
\label{eq:momentum_balance_discrete_residual}
\end{split}\\
R^{\mathrm{mass}}_I =& \sum_{J=1}^{n_n} \left( \sum_{p=1}^{n_{p}} \left( \theta_p  S_{Ip} \nabla S^T_{Jp} + \left(\nabla \theta_{p}\right)^T S_{Ip}  S_{Jp} \right) V_p \right) \dot{\tb u}_J\,.
\label{eq:mass_balance_discrete_residual}
\end{align}
\end{linenomath}
In the above expressions, $\tb B_{Ip}$ denotes the nodal deformation matrix commonly employed in the standard FE formulation \cite{hughes1986new}. Following the common practice in MPM, we adopt the lumped grid nodal mass, $m_I^n$, owing to its practicality and efficiency. It's worth noting that in MPM, applying traction forces over the Neumann boundary $\Gamma_N$ requires modeling efforts as the position of the boundary, $\tb x_{\Gamma_N}$, typically does not conform to the background grid \cite{chandra2021nonconforming, liang2023imposition}. For brevity, these terms will be omitted in the subsequent formulations.

Following the conventional nonlinear FE procedure, the iterative solution to the balance equations can be obtained by employing Newton-Raphson's method, utilizing directional derivatives. To perform the linearization with respect to the primary variables, the Newmark-$\beta$ time integration scheme is utilized, where the updated nodal velocity and acceleration vectors, $\tb v_I^k$ and $\tb a_I^k$, can also be computed. Following the Newmark-$\beta$ scheme, the nodal velocity and acceleration can be linearized as:
\begin{eqnarray}
\delta_{\tb u_I}\left(\tb v_I\right) = \frac{\gamma_N}{\beta_N \Delta t}\,, \qquad \delta_{\tb u_I}  \left(\tb a_I\right) = \frac{1}{\beta_N \Delta t^2}\,,
\label{eq:linearization_kinematics}
\end{eqnarray}
where $\delta_{\square}(\square)$ is the linearization operator with subscripts denoting the variables linearized. Meanwhile, $\gamma_N$ and $\beta_N$ are the Newmark's coefficients and $\Delta t$ is the considered time increment. In this study, values of $\gamma_N=1.0$ and $\beta_N=0.5$ are adopted to introduce algorithmic dissipation and improve numerical stability.

Employing Taylor's expansion and Eq.~\eqref{eq:linearization_kinematics}, the linearization of the weak form, Eqs.~\eqref{eq:weak_form_momentum} and \eqref{eq:weak_form_mass}, yields the following mixed linear systems of equations to be solved at every Newton-Raphson's iteration $k$:
\begin{eqnarray}
\begin{bmatrix}
    \tb K^k_{uu} & \tb K^k_{up} \\
    \tb K^k_{pu} & \tb 0
\end{bmatrix}
\begin{Bmatrix}
    \Delta \tb{u}^{k+1} \\
    \Delta {p}^{k+1}
\end{Bmatrix} = 
- \begin{Bmatrix}
    \tb R^{\mathrm{mom},k} \\
    R^{\mathrm{mass},k}
\end{Bmatrix}\,,
\label{eq:lin_sys_equation}
\end{eqnarray}
where the block-wise Jacobian matrices given in the discrete form read:
\begin{linenomath}
\begin{align}
\begin{split}
\left(\tb{K}^k_{uu}\right)_{IJ}
=&  \frac{1}{\beta_N \Delta t^2} \tb I_{IJ} m_J^n + \sum_{p=1}^{n_{p}} \left(  \frac{\gamma_N}{\beta_N \Delta t}\tb B_{Ip}^T \mathbb{C} \tb B_{Jp}\right) V^k_p + \frac{\gamma_N}{\beta_N \Delta t} \tb I_{IJ} \left( \sum_{p=1}^{n_{p}} \left( \widetilde{A}_p + \widetilde{B}_p \left\Vert \dot{\tb u}_p^k \right \Vert
\right) S_{Jp} m_p \right)\\
& +  \frac{\gamma_N}{\beta_N \Delta t} \left( \sum_{p=1}^{n_{p}} \widetilde{B}_p 
{\tb w}_p^k \otimes \dot{\tb u}_p^k S_{Ip}S_{Jp} m_p \right)\,,
\label{eq:LHS_matrix_uu}
\end{split}\\
\left(\tb{K}^k_{up}\right)_{IJ}
=&- \sum_{p=1}^{n_{p}} \left( \theta_p^k \nabla S_{Ip} S_{Jp} + \nabla \theta_p^n S_{Ip} S_{Jp}\right)  V^k_p
\,,
\label{eq:LHS_matrix_up}\\
\left(\tb{K}^k_{pu}\right)_{IJ}
=& \sum_{p=1}^{n_{p}} \frac{\gamma_N}{\beta_N \Delta t} \left( \theta_p^k  S_{Ip} \nabla S^T_{Jp} + \left(\nabla \theta^{n}_{p}\right)^T S_{Ip}  S_{Jp} \right)  V^k_p 
\,.
\label{eq:LHS_matrix_pu}
\end{align}
\end{linenomath}
Here, $\tb K_{uu}$ is the Jacobian matrix corresponding to the displacement field, which includes the dynamic, material, and drag components. The constitutive tensor $\mathbb{C}$ is a fourth-order tensor containing the porosity-dependent viscosity given previously in Eq.~\eqref{eq:porosity_dependent_visc}. It can be written in the Voigt notation as:
\begin{eqnarray}
    \mathbb{C} = \mathrm{diag}\left(2,2,2,1,1,1\right)\widetilde{\mu}^k\,.
\end{eqnarray}
Furthermore, notice that the last term in the expression of $\boldsymbol{K}_{uu}$ appears as a consequence of linearizing the velocity norm found in the nonlinear drag force term. In this context, we represent the normalized velocity vector as $\tb w_p = \dot{\tb u}_p/||\dot{\tb u}_p||$. Meanwhile, $\tb K_{up}$ and $\tb K_{pu}$ are the off-diagonal coupling matrices, whose values are different, thus, yielding to an asymmetric global Jacobian matrix. However, since $\tb K_{pu}$ can be expressed as a linear combination of $\tb K_{up}$, i.e.~$\tb K_{pu}=-\gamma_N/(\beta_N \Delta t)\tb K_{up}^T$, the assembly process can be greatly simplified. 

When evaluating the Jacobian matrices and residual vectors, the porosity field should be iteratively updated following the fluid deformation. While it is convenient to assume the porosity to remain unchanged through the nonlinear iterations, i.e.~$\theta^k=\theta^n$, this assumption is only accurate in three conditions: (i) when the fluid remains in either the non-porous or the porous regions, (ii) the thickness of the blurred interface is sufficiently large such that the porosity gradient is small, i.e.~$\varepsilon \gg h$, and (iii) the deformation at every time step is reasonably small, i.e.~$||\tb u|| \ll h$. Since these conditions may not always be guaranteed, a first-order approximation by Taylor expansion can be considered to update the porosity at every iteration $k$, i.e.:
\begin{eqnarray}
    \theta^k_p &=& \theta^n_p + \nabla \theta^n_p \cdot \tb u^k_p + \mathcal{O}\left((\tb u^k_p)^2\right)\,.
    \label{eq:iterative_porosity}
\end{eqnarray}

Here, to avoid discretization of higher-order derivatives, the porosity gradient is assumed to follow the previous time step configuration, $\nabla \theta^n_p$. However, because of this assumption, the approximation given by Eq.~\eqref{eq:iterative_porosity} may induce an over-projection error of porosity value. It is important to ensure the updated porosity magnitude lies within the physical bound:
\begin{eqnarray}
    0\leq\theta^{\min}\leq\theta^k_p\leq\theta^{\max}\leq1\,,
\end{eqnarray}
where $\theta^{\min}:=\min_I\left(\theta_I\right)$ and $\theta^{\max}:=\max_I\left(\theta_I\right)$ denote the minimum and maximum value of nodal porosity within the connectivity of the material point $p$, respectively. 

Once the material point porosity is updated, the material point volume and viscosity should be updated within the iteration as:
\begin{eqnarray}
    V_p^k &=& V_p^n \theta_p^n/\theta_p^k\,,\\
    \widetilde{\mu}^k &=&  \left(1+\eta\,(1-\theta^k_p) \right)\mu\,.
\end{eqnarray}

\begin{remark}
    It is worth mentioning that, even though the first-order approximation of porosity is considered according to Eq.~\eqref{eq:iterative_porosity}, the linearization of the porosity term with respect to the displacement field is omitted when deriving the Jacobian matrix for practicality, i.e.~$\delta_{\tb u_I} \left( \theta \right)\approx 0$.
\end{remark}

\subsection{Stabilized formulation by the variational multiscale (VMS) method}

Given that the primary variables are discretized using equal-order interpolation functions (see \Cref{subsec:spatial_disc}), the use of a stabilization technique becomes necessary to satisfy the discrete inf-sup conditions  \cite{hughes1986new}. Following our previous work \cite{chandra2023stabilized}, the VMS stabilization is derived and implemented for the unified balance equations. For detailed mathematical formulation and theoretical proof of the VMS method in mixed FEM, interested readers are strongly recommended to refer to \cite{hughes1995multiscale, hughes1998variational, codina2001stabilized, codina2002stabilized, codina2007time}.

The basis of the VMS method hinges on the division of the solution fields into two distinct scales:
\begin{subequations}
\begin{eqnarray}
    \dot{\tb u} = \dot{\tb u}_h + \dot{\tb u}_s\,, \label{eq:velocity_vms}\\
    p = p_h + p_s\,.
\end{eqnarray}
\label{eq:vms_fields}%
\end{subequations}
Here, $\dot{\tb{u}}_h$ and $p_h$ denote the resolvable scales, while $\dot{\tb{u}}_s$ and $p_s$ represent the subscale fluctuations, which cannot be resolved by the considered discretization. As given in Eq.~\eqref{eq:velocity_vms}, the scale separation is executed within the velocity field, aligning with the typical Eulerian CFD approach. It is worth emphasizing that the VMS method considered in this study does not alter the flow constitutive relation, unlike approaches in turbulence modeling; rather, its primary purpose is to resolve instabilities related to the compatibility conditions of the interpolation. Past research has explored utilizing the VMS method for turbulence modeling, given its close similarities with the large eddy simulation (LES) formulation, e.g.~\cite{hughes2000large}.

Eq.~\eqref{eq:vms_fields} is subsequently substituted to the weak form of the momentum and mass balances, Eqs.~\eqref{eq:weak_form_momentum}-\eqref{eq:weak_form_mass}, which yield to the following modified combined weak form: 
\begin{eqnarray}
\begin{split}
\widehat{R}({\tb u}_h, p_h,{\tb u}_s, p_s,  \delta \tb u_h, \delta p_h ) =& 
\int_{\Omega}\theta \rho \left(\ddot{\tb u}_h + \ddot{\tb u}_s - \tb b\right) \cdot \delta \tb u_h\, \td{\rm \Omega}
-\int_{\Omega} \left( p_h + p_s\right) \left(\theta \left( \nabla \cdot \delta \tb u_h \right) + \nabla\theta \cdot \delta \tb u_h \right)\, \td\Omega\\
& +\int_{\Omega} \mathbb{C}:\nabla^s (\dot{\tb u}_h + \dot{\tb u}_s) \colon \nabla^s \delta \tb u_h\, \td\Omega \\
& + \int_{\Omega} \theta \rho \left( \widetilde{A} + \widetilde{B}||\dot{\tb u}_h + \dot{\tb u}_s||
\right) \left(\dot{\tb u}_h+ \dot{\tb u}_s\right) \cdot \delta \tb u_h\, \td \Omega \\
& + \int_{\Omega} \nabla \cdot \left( \theta\left(\dot{\tb u}_h + \dot{\tb u}_s \right) \right) \, \delta p_h \, \td{\rm \Omega} = 0
\,.
\label{eq:weak_form_vms_1}
\end{split}
\end{eqnarray}
In this context, $\nabla^s\square$ represents the symmetric gradient operator, while $\delta \tb{u}_h$ and $\delta p_h$ stand for the resolved trial functions for the corresponding solution fields. It is important to highlight that both the subscale variables $\dot{\tb{u}}_s$ and $p_s$ are assumed to have zero boundary integrals. Additionally, the traction term related to the resolvable scales is excluded, as done previously.

In the current work, we follow the quasi-static Algebraic Sub-Grid Scales (ASGS) \cite{codina2001stabilized, codina2007time} approach where the subscale variables are modeled by projecting the FE residuals as:
\begin{eqnarray}
\dot{\tb u}_s = \tau_1 \widetilde{\tb R}^{\mathrm{mom}}(\dot{\tb u}_h, p_h)\,,&\qquad 
\ddot{\tb u}_s = \tb 0 \,,\\
\label{eq:vms_subscale_displacement}
p_s = \tau_2 \widetilde{R}^{\mathrm{mass}}(\dot{\tb u}_h)\,, &\qquad \dot{p}_s = 0\,.
\label{eq:vms_subscale_pressure}
\end{eqnarray}
Here, $\widetilde{\tb R}^{\mathrm{mom}}$ and $\widetilde{R}^{\mathrm{mass}}$ are the strong form residual of the momentum and mass balance given as:
\begin{subequations}
\begin{eqnarray}
    \widetilde{\tb R}^{\mathrm{mom}}(\dot{\tb u}_h, p_h) &=& \theta \rho \tb b - \theta \rho \ddot{\tb u}_h - \theta \nabla p_h + \nabla \cdot \left(\mathbb{C} : \nabla^s \dot{\tb u}_h \right) - \theta \rho \left( \widetilde{A} + \widetilde{B} ||\dot{\tb u}_h|| \right) \dot{\tb u}_h\,,\\
    \widetilde{R}^{\mathrm{mass}}(\dot{\tb u}_h) &=& -\nabla \cdot \left(\theta\dot{\tb u}_h \right)\,.
\end{eqnarray}
\label{eq:vms_subscale_residuals}%
\end{subequations}
Furthermore, the stabilization coefficients $\tau_1$ and $\tau_2$ are given in the form proposed by \citet{zorrilla2019modified} though with a modification to include the effect of drag force coefficients as suggested by \citet{larese2015finite}:
\begin{eqnarray}
\tau_1 = \left(\frac{\overline{\theta} \rho \tau_{dyn}}{\Delta t} + \frac{c_2 \overline{\theta} \rho||\overline{\dot{\tb u}}_h||}{h_\#} + \frac{c_1 \widetilde{\mu}}{h_\#^2} + \overline{\theta} \rho \left(\overline{\widetilde{A}} + \overline{\widetilde{B}} ||\overline{\dot{\tb u}}_h|| \right) \right)^{-1}\,, \quad
\tau_2 = \frac{h_\#^2}{c_1 \tau_1}\,,
\label{eq:stabilization_coefficients}
\end{eqnarray}
where $c_1=4.0$, $c_2=2.0$, and $\tau_{dyn}=1.0$. Here, $\widetilde{\mu}$ is the porosity-dependent fluid viscosity, whereas $||\overline{\dot{\tb u}}_h||$ and $h_\#$ are the cell-wise averaged velocity norm and the equivalent cell length, respectively. Furthermore, since $\tau_1$ and $\tau_2$ are assumed to be constant per cell, the porosity and permeability parameters, $\overline{\theta}$, $\overline{\widetilde{A}}$ and $\overline{\widetilde{B}}$, are to be evaluated at the cell-center in a similar way as the averaged velocity field.

With the aforementioned assumptions and by omitting the nonlinear and higher-order terms appearing in the drag and viscous terms, the modified weak formulation can be simplified to:
\begin{eqnarray}
\begin{split}
\widehat{R}(\tb u_h, p_h,  \delta \tb u_h, \delta p_h ) =&  R^{\mathrm{mom}}(\tb u_h, p_h, \delta \tb u_h) + R^{\mathrm{mass}}(\tb u_h, \delta p_h) \\
& + \int_\Omega \tau_2 \left( \theta \nabla \cdot \dot{\tb u}_h + \dot{\tb u}_h \cdot \nabla\theta \right) \left( \theta \left(\nabla \cdot \delta \tb u_h \right) + \nabla \theta \cdot \delta \tb u_h \right)\, \td \Omega\\
& + \int_{\Omega} \tau_1 \theta^2 \rho \nabla \delta p_h \cdot \left( \ddot{\tb u}_h-\tb b + \frac{\nabla p_h}{\rho} + \left( \widetilde{A} + \widetilde{B} ||\dot{\tb u}_h|| \right) \dot{\tb u}_h \right)   \, \td{\rm \Omega}
= 0
\,,
\label{eq:weak_form_vms_2}
\end{split}
\end{eqnarray}
where $R^{\mathrm{mom}}$ and $R^{\mathrm{mass}}$ are the standard weak-form of momentum and mass balance residuals given previously in Eqs.~\eqref{eq:weak_form_momentum} and \eqref{eq:weak_form_mass}, respectively. Moreover, by performing a linearization, as well as the spatiotemporal discretization, the modified system of equations can be obtained as:
\begin{linenomath}
\begin{eqnarray}
\begin{bmatrix}
    \tb K^k_{uu} + \widehat{\tb K}^k_{uu} & \tb K^k_{up} \\
    \tb K^k_{pu} + \widehat{\tb K}^k_{pu} & \widehat{\tb K}^k_{pp}
\end{bmatrix}
\begin{Bmatrix}
    \Delta \tb{u}^{k+1} \\
    \Delta {p}^{k+1}
\end{Bmatrix} = 
- \begin{Bmatrix}
    \tb R^{\mathrm{mom},k} + \widehat{\tb R}^{\mathrm{mom},k} \\
    R^{\mathrm{mass},k} + \widehat{R}^{\mathrm{mass},k}
\end{Bmatrix}\,,
\label{eq:mod_lin_sys_equation}
\end{eqnarray}
where the additional stabilization matrices and vectors are denoted by $\widehat{\square}$ notation. The additional momentum and mass residuals are first given as:
\begin{align}
\begin{split}
\widehat{\tb R}_I^{\mathrm{mom},k} =& \sum_{p=1}^{n_p} \tau_2^k  \Biggl( \left(\theta^k_p ({\dot{\varepsilon}_v})^k_p  + \left(\dot{\tb u}^k_p \cdot \nabla \theta^n_p\right) \right) \left( \theta^k_p \nabla S_{Ip} + \nabla \theta^n_p S_{Ip} \right) \Biggl) V_p^k 
\,, \label{eq:vms_mom_res} 
\end{split} \\
\begin{split}
\widehat{R}_I^{\mathrm{mass},k} =&  \sum_{p=1}^{n_p} \tau_1^k \theta^k_p \nabla S^T_{Ip} \cdot \left(\ddot{\tb u}_p^k - \tb b + \left( \widetilde{A}_p + \widetilde{B}_p ||\dot{\tb u}^k_p|| \right)\dot{\tb u}_p^k\right) m_p  \\
& + \sum_{J=1}^{n_n} \left( \sum_{p=1}^{n_p}  \tau_1^k \left(\theta^k_p\right)^2 \nabla S^T_{Ip} \cdot \nabla S_{Jp} V_p^k \right) p^k_J\,,
\label{eq:vms_mass_res}
\end{split}
\end{align}
\end{linenomath}
while the additional Jacobian matrices read:
\begin{linenomath}
\begin{align}
\begin{split}
\left(\widehat{\tb{K}}_{uu}^k\right)_{IJ}
=& \frac{\gamma_N}{\beta_N \Delta t} \sum_{p=1}^{n_p} \tau_2^k \Biggl(\left(\theta^k_p\right)^2 \nabla S_{Ip} \nabla S^T_{Jp} \\
& \qquad + \theta^k_p \nabla S_{Ip} S_{Jp} \left(\nabla \theta^n_p\right)^T + \theta^k_p S_{Ip} \nabla \theta^n_p \nabla  S^T_{Jp} + S_{Ip} S_{Jp} \nabla \theta^n_p \left(\nabla \theta^n_p\right)^T\Biggl)  V_p^k
\,,\label{eq:vms_kuu}
\end{split}\\
\left(\widehat{\tb{K}}_{pu}^k\right)_{IJ}
=& \sum_{p=1}^{n_p}  \tau_1^k \theta^k_p \nabla S^T_{Ip} \left(\frac{1}{\beta_N \Delta t^2} + \frac{\gamma_N}{\beta_N \Delta t}\left( \widetilde{A}_p + \widetilde{B}_p ||\dot{\tb u}_p^k|| + \widetilde{B}_p {\tb w}_p^k \otimes \dot{\tb u}_p^k \right) \right) S_{Jp} m_p  \,, \label{eq:vms_kpu}\\
\left(\widehat{\tb K}_{pp}^k\right)_{IJ} 
=& \sum_{p=1}^{n_p} \tau_1^k \left(\theta^k_p\right)^2 \nabla S^T_{Ip} \cdot \nabla S_{Jp} V_p^k\,. \label{eq:vms_kpp}
\end{align}
\end{linenomath}
Here, $\left( {\dot{\varepsilon}_v} \right)^k_p = \sum_{I=1}^{n_n} \nabla S^T_{Ip} \cdot \dot{\tb u}^k_I$ is the divergence of velocity field evaluated at the material point. It's important to note that the VMS formulation derived above will converge to the VMS formulation for incompressible Newtonian fluids in the limit of a free-surface flow, i.e.~when $\theta\rightarrow 1$, which correspondingly implies that $\widetilde{A}\rightarrow 0$, $\widetilde{B}\rightarrow 0$, and $\nabla \theta\rightarrow \tb 0$.

\begin{remark}
For the sake of practicality, the linearization of the stabilization coefficients with respect to the displacement field is omitted when deriving the stabilization matrices, i.e.~$\delta_{\tb{u}_I}\left(\tau_1\right)=0$ and $\delta_{\tb{u}_I}\left(\tau_2\right)=0$. We make this assumption because the change in the averaged velocity field, $\overline{\dot{\tb{u}}}_p$, is observed to be significantly small during the Newton iteration, and hence, the stabilization coefficients remain relatively constant, i.e.~$\tau^k_1\approx\tau^n_1$ and $\tau^k_2\approx\tau^n_2$.
\end{remark}

\subsection{Solving systems of nonlinear equations}

To solve systems of nonlinear equations using the Newton-Raphson method, Eq.~\eqref{eq:mod_lin_sys_equation} must be iteratively solved until convergence is achieved (for a detailed algorithm, c.f.~\cite{chandra2023stabilized}). In this work, the Krylov solver provided by PETSc \cite{petsc-web-page} is employed. Within each iteration, an approximate field-split Schur complement preconditioner is constructed. Meanwhile, we utilize the Jacobi and geometric algebraic multigrid (GAMG) preconditioners for the pressure and displacement fields, respectively. The linear solver of choice is the flexible generalized minimal residual (FGMRES) method. Moreover, the mixed MPM solver has been developed with hybrid shared and distributed parallelization features, incorporating the graph-partitioning domain decomposition library, KaHIP\cite{MeyerhenkeSS17}. However, it is worth noting that a comprehensive assessment of the solver's efficiency and scalability is a separate research effort, which we intend to address in a separate future work.

\subsection{Treatment of quadrature errors}

As demonstrated in several recent works, e.g.~\cite{baumgarten2019general, baumgarten2023analysis}, simulating fluid flow within the MPM often faces challenges related to poor numerical integration, primarily due to particle clustering in regions with significant shearing. For instance, in fluid flows with stagnation points, particles tend to aggregate along compressing streamlines. Using these clustered particles as quadrature points leads to significant errors, potentially compromising both solution accuracy and mass and momentum conservation\cite{chandra2023stabilized}. To address these quadrature errors, a particle-shifting technique (PST) can be introduced into the MPM formulation. However, the MPM's data structure, which combines material points and background grids, presents limitations for PST based on inter-particle proximity. This is primarily due to the computational cost introduced by the nearest neighbor searching algorithm for each material point. In this study, we incorporate the $\delta$-correction scheme, as suggested by \citet{baumgarten2023analysis}, due to its compatibility with our MPM code's data structure.

\section{Numerical examples}
\label{sec:num_ex}

In this section, we evaluate the quality and robustness of the proposed formulation and stabilization by performing a series of (quasi-)1D, 2D, and 3D numerical examples. The summary of the conducted numerical examples, along with their specific purposes, is given in Table \ref{tab:numerical_examples}. In all these examples, structured background grids with quadratic B-Spline basis functions are employed, where for all the 1D and 2D problems, single-step kernel correction is used \cite{nakamura2023taylor}, while for the 3D simulations, two-step kernel correction is considered (see \ref{app:basis_funct_it_kc}). The convergence tolerances for Newton-Raphson's method are set to $\epsilon^{\mathrm{rel}}_{\rm tol}=10^{-10}$ for the relative residual criterion and $\epsilon^{\mathrm{en}}_{\mathrm{tol}}=10^{-15}$ for the energy criterion. They are given as \cite{chandra2023stabilized}:
\begin{eqnarray}
    \mathcal{C}^{\mathrm{rel}}:=\frac{||\tb R^k||}{||\tb R^0||} \leq \epsilon^{\mathrm{rel}}_{\rm tol}\,, \qquad 
    \mathcal{C}^{\mathrm{en}}:=\sqrt{-\tb R^{\mathrm{mom},{k}}\cdot \frac{\gamma_N}{\beta_N \Delta t}\Delta \tb u^{k+1} - R^{\mathrm{mass},{k}}\cdot \Delta p^{k+1}} \leq \epsilon^{\mathrm{en}}_{\mathrm{tol}}\,,
    \label{eq:convergence_criteria}
\end{eqnarray}
where $||\tb R||$ is the $L^2$-norm of the combined residual vector. The relative convergence tolerance for all the Krylov solvers used is set to $10^{-10}$. For all cases, unless specified differently, the fluid is assumed to be water with mass density $\rho=1000$ $\mathrm{kg/m^3}$ and dynamic viscosity $\mu=0.001$ $\mathrm{Pa\cdot s}$ in the non-porous region.

\begin{table}[h!]
\centering
\caption{Summary of conducted numerical examples along with their specific purposes.}
\label{tab:numerical_examples}
\small
\resizebox{\textwidth}{!}{%
\begin{threeparttable}
\begin{tabular}{||c|c|l|c|c|c||}
\hline
\# & Section & \multicolumn{1}{c|}{Verification and validation purpose}                     & Dim. & Ref. & Type \\ \hline \hline
1.1   &   \ref{subsubsec:1d_case1}     & Volume conservation, solver performance                    & Q1D (2D)  &    \cite{sun2019numerical}         & T\\ \hline
1.2   &   \ref{subsubsec:1d_case2}     & Volume conservation, solver performance                    & Q1D (2D)  &    \cite{akbari2014modified}       & T\\ \hline
1.3   &   \ref{subsubsec:1d_case3}     & Steady-state kinematics (seepage)                          & Q1D (2D)  & \ref{app:1d_case3_derivation}  & T\\ \hline
1.4   &   \ref{subsubsec:1d_case4}     & Transient kinematics and pressure (coupled)                  & 1D        & \ref{app:1d_case4_derivation}  & T\\ \hline
1.5   &   \ref{subsubsec:1d_case5}     & Transient kinematics and pressure (coupled)                   & 1D        & \ref{app:1d_case5_derivation}  & T\\ \hline
1.6   &   \ref{subsubsec:1d_case6}     & Transient kinematics and pressure (seepage, multi layers)     & 1D        & \ref{app:1d_case6_derivation}  & T\\ \hline
2   &     \ref{subsec:couette_flow_composite}    & Steady-state kinematics (coupled)                   & 2D        &    \cite{kuznetsov1998analytical}       & T\\ \hline
3   &     \ref{subsec:2d_dam_break}    & Free-surface topology             & 2D        &  \cite{liu1999numerical}         & E, N \\ \hline
4   &     \ref{subsec:3d_dam_break}    & Free-surface topology, volume conservation                 & 3D        &     \cite{del2012three}      & N, T\\ \hline
\end{tabular}%
\begin{tablenotes}
\item[Q1D] Quasi-1D.
\item[T] Theoretical/analytical verification.
\item[E] Experimental validation.
\item[N] Numerical benchmark.
\end{tablenotes}
\end{threeparttable}
}
\end{table}

\subsection{One-dimensional gravity-driven verification suite}
\label{subsec:1d_suite}

\begin{figure}[h!]
    \centering
    \includegraphics[width=0.9\textwidth]{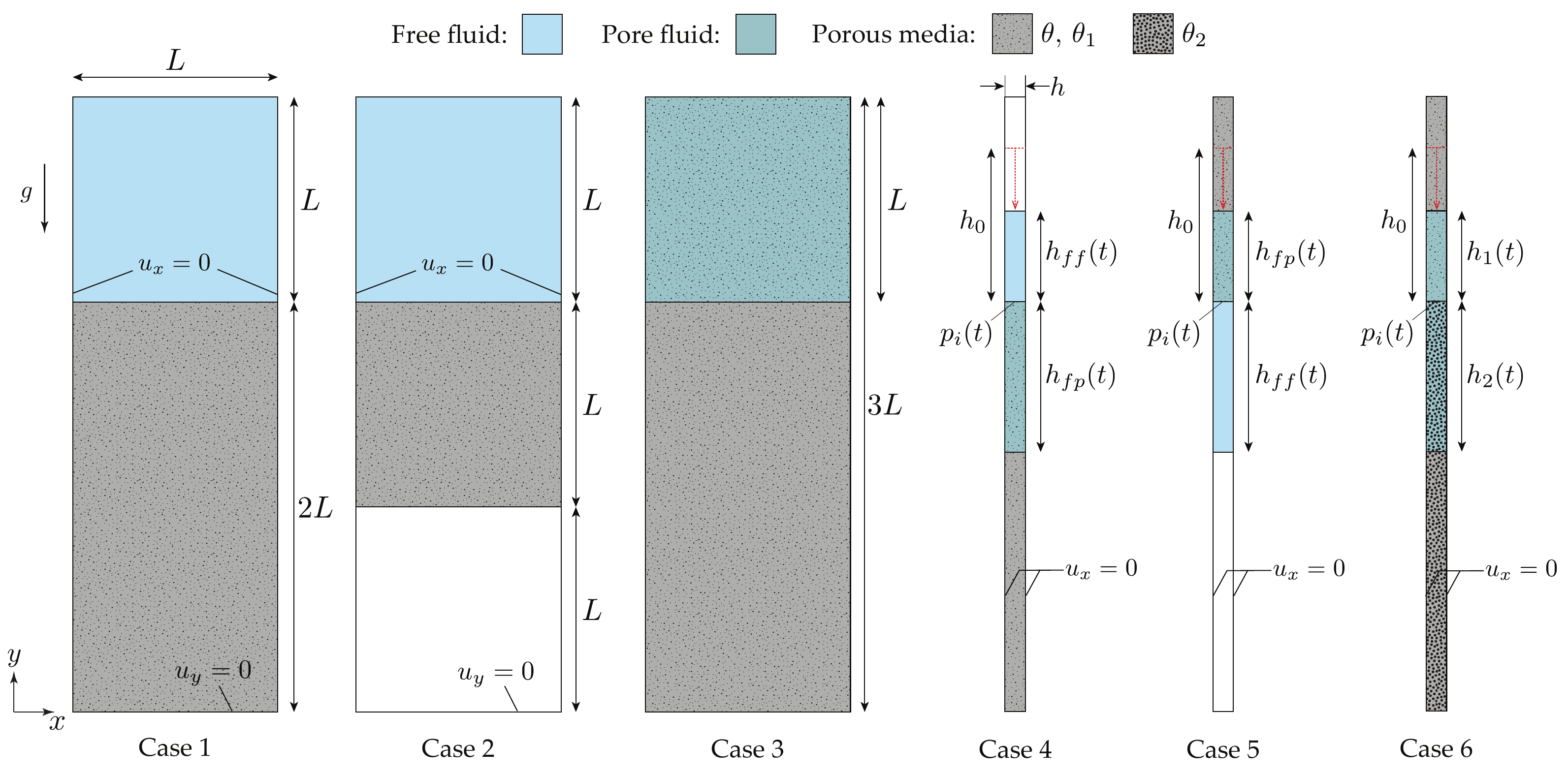}
    \caption{1D verification suite: geometrical settings of six cases and boundary conditions.}
    \label{fig:3.1_model}
\end{figure}

A series of (quasi-)1D gravity-driven ($g=9.81\,\rm{m/s^2}$) simulations are conducted to rigorously evaluate and verify the numerical accuracy, stability, and computational efficiency of the proposed method. This study encompassed six different numerical cases, each devised with specific purposes. The geometrical configurations for these cases can be categorized into two groups, as illustrated in Fig.~\ref{fig:3.1_model}. In the first set of cases, Cases 1-3, the simulations are carried out in a quasi-1D (2D) domain, with the primary gravitational force acting in the vertical direction along the $y$-axis. Here, horizontal movement is partially unconstrained, as the fluid material points at the centerline may move horizontally in the $x$-direction, though the left and right walls are subjected to homogeneous displacement BC. This setting is chosen for the purpose of quantifying macroscopic numerical properties, including volume conservation and assessing solver convergence. The second set of cases, Cases 4-6, are configured with a single element width horizontally. In this setting, horizontal motion is entirely prohibited due to enforced $x$-directional constraints on the left and right edges. These cases are designed to measure the accuracy and stability of the numerical results, e.g.~fronts propagation and interface pressure, by comparing them to the 1D analytical solutions; this will be derived and given further in \ref{app:analytical_solutions}. While not necessarily realistic, this study neglects capillary action acting at the fluid fronts for simplicity.

Inside the porous medium, the effective viscosity is considered to be the same as the fluid viscosity $\mu_e = \mu$, hence $\widetilde{\mu}=\mu$. Here, the drag force coefficients are computed using the hydraulic conductivity $K$ and by assuming Ergun's constants ($A=150$ and $B=1.75$) for the nonlinear drag components (c.f.~\ref{app:drag_force_coeff}). To compare accuracy and stability, the fractional-step approach is implemented, which details can be found in \ref{app:fractional_step}. Moreover, for this benchmark problem, only the FLIP transfer scheme is considered.

\subsubsection{Case 1 -- volume conservation of fluid entering porous domain}
\label{subsubsec:1d_case1}

In the first case, a fluid square patch with size $L=0.2$ m penetrates a porous region with porosity $\theta=0.5$. Due to the size of the porous domain, which is exactly double the size of the fluid patch, the fluid is expected to fill the entire porous domain at the end of the simulation \cite{sun2019numerical}. The MPM simulations performed for this study consider a background grid with size $h=L/10$ and 4 × 4 particles per cell (PPC), totaling 1600 material points. The time step is set as $\Delta t=0.001$ s, while the hydraulic conductivity of the porous media is considered to be considerably large, i.e.~$K= 0.1$ m/s, to reduce the number of steps needed to complete the simulation. Here, linear Darcy's law is assumed, hence, the nonlinear coefficient $\widetilde{B}$ can be set as zero.

\begin{figure}[h!]
    \centering
     \begin{subfigure}[b]{0.46\textwidth}
         \centering
         \includegraphics[width=\textwidth]{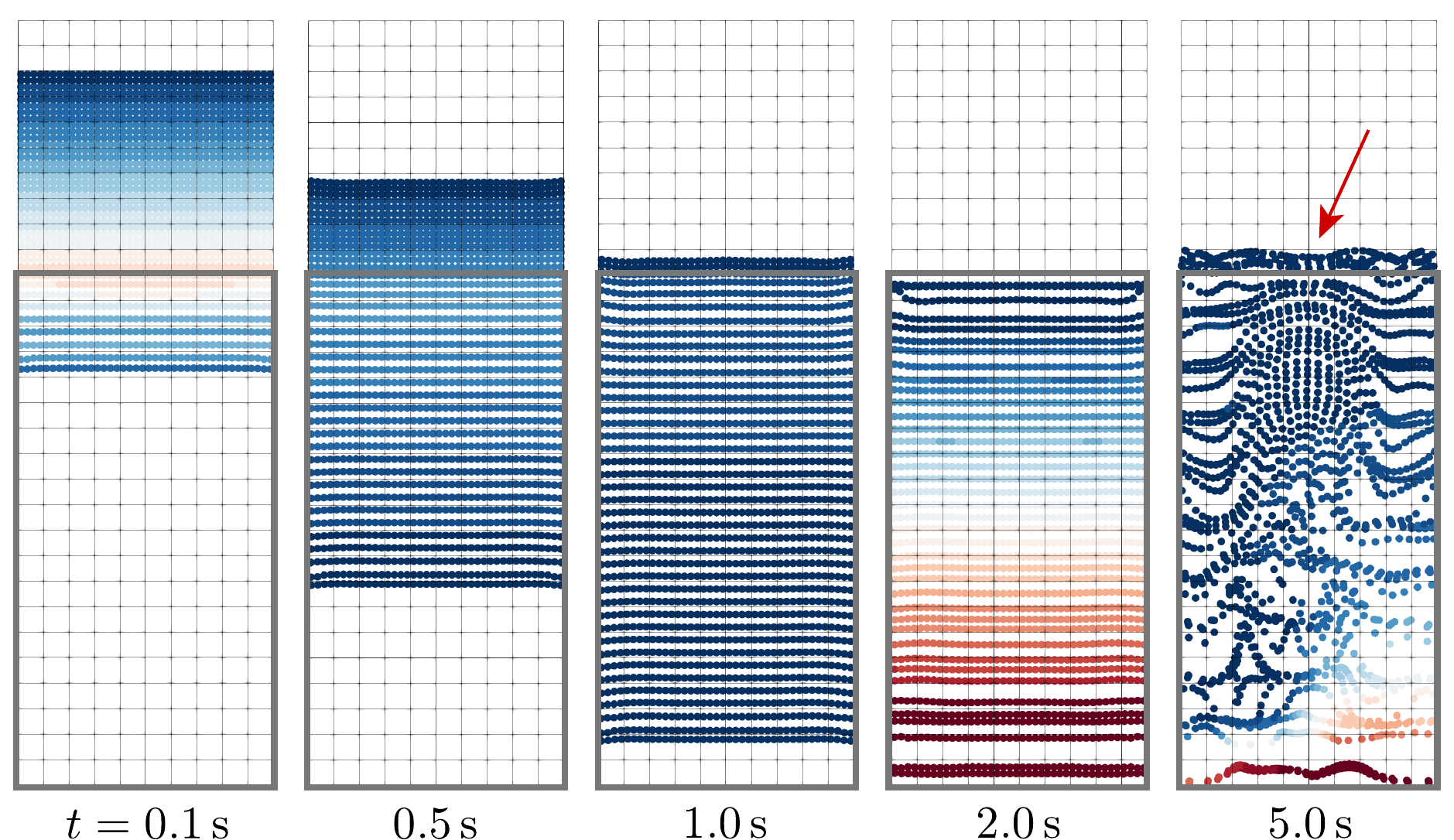}
         \caption{Fractional-step}
         \label{fig:3.1.1_screenshots_fs}
     \end{subfigure}
     \raisebox{0.4cm}{\rule{0.1pt}{4.45cm}}
     \begin{subfigure}[b]{0.46\textwidth}
         \centering
         \includegraphics[width=\textwidth]{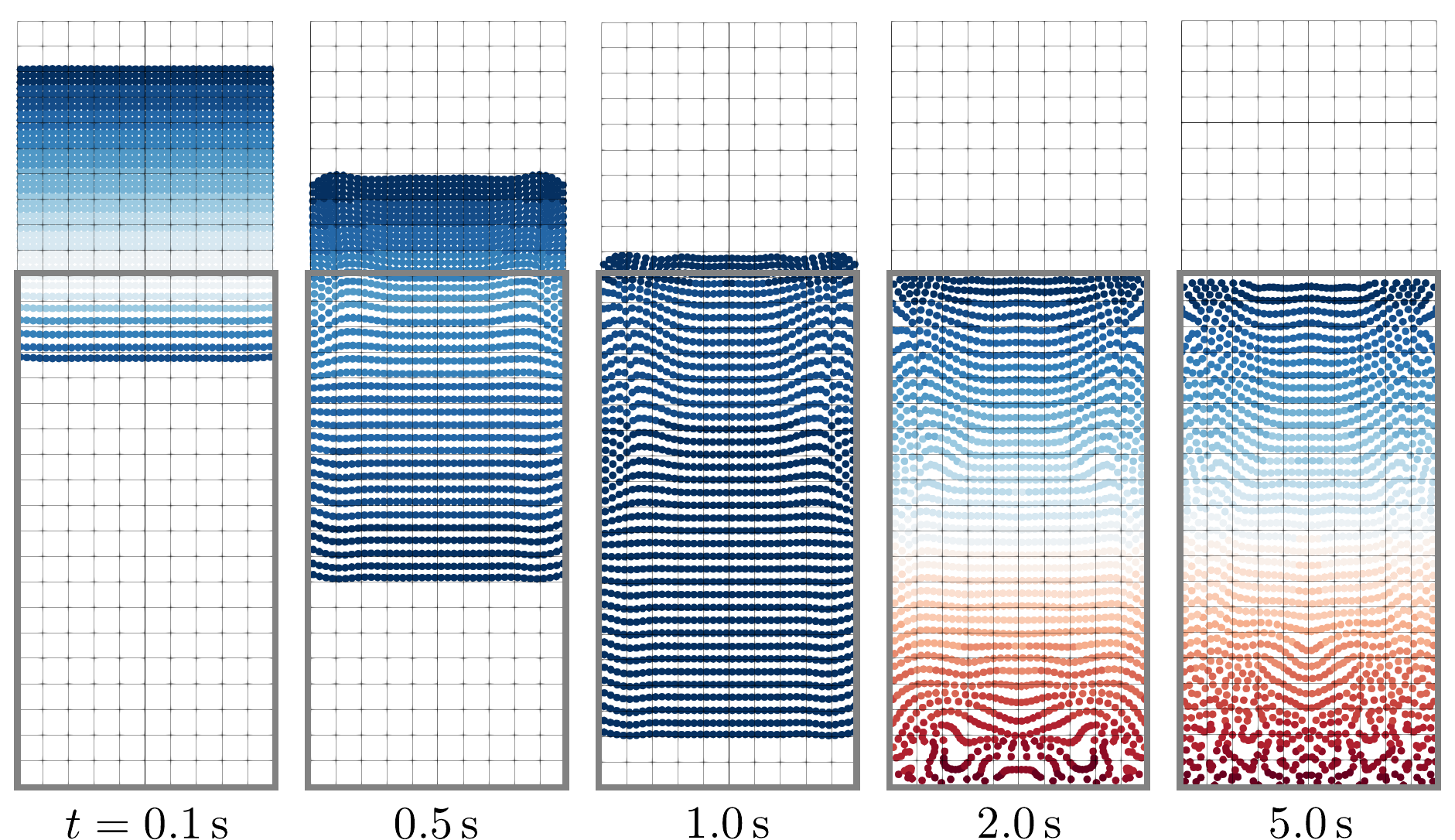}
         \caption{Mixed--VMS}
         \label{fig:3.1.1_screenshots_vms}
     \end{subfigure}
     \begin{subfigure}{0.053\textwidth}
         \centering
         \includegraphics[width=\textwidth]{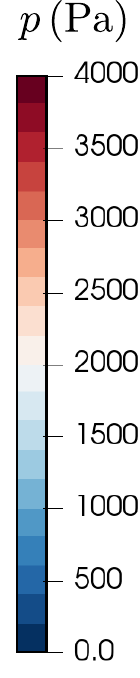}
         \caption*{}
     \end{subfigure}
    \caption{1D verification suite - Case 1: evolution of pressure field at different time snapshots $t=0.1$, 0.5, 1.0, 2.0, and 5.0 s. The red arrow highlights volumetric and pressure errors obtained by the fractional-step method.}
    \label{fig:3.1.1_results}
\end{figure}

\begin{figure}[h!]
    \centering
     \begin{subfigure}[b]{0.41\textwidth}
         \centering
         \includegraphics[width=\textwidth]{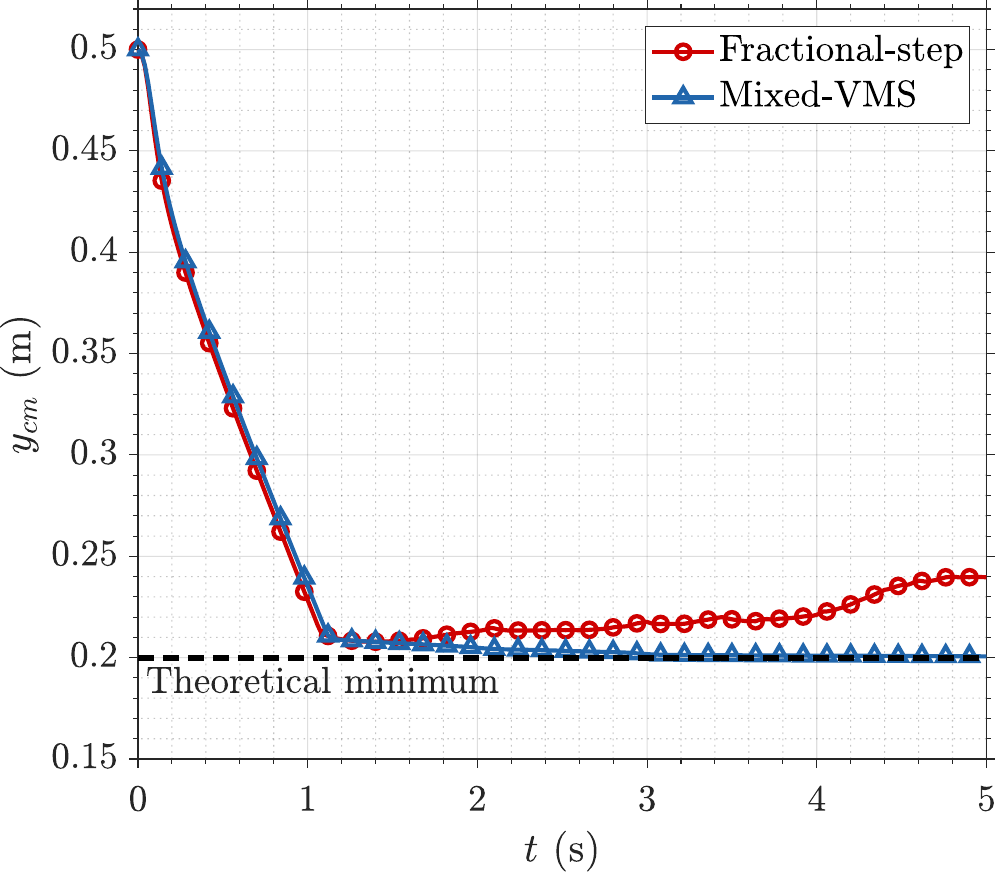}
         \caption{}
         \label{fig:3.1.1_com}
     \end{subfigure}
     \hspace{0.2cm}
     \begin{subfigure}[b]{0.4\textwidth}
         \centering
         \includegraphics[width=\textwidth]{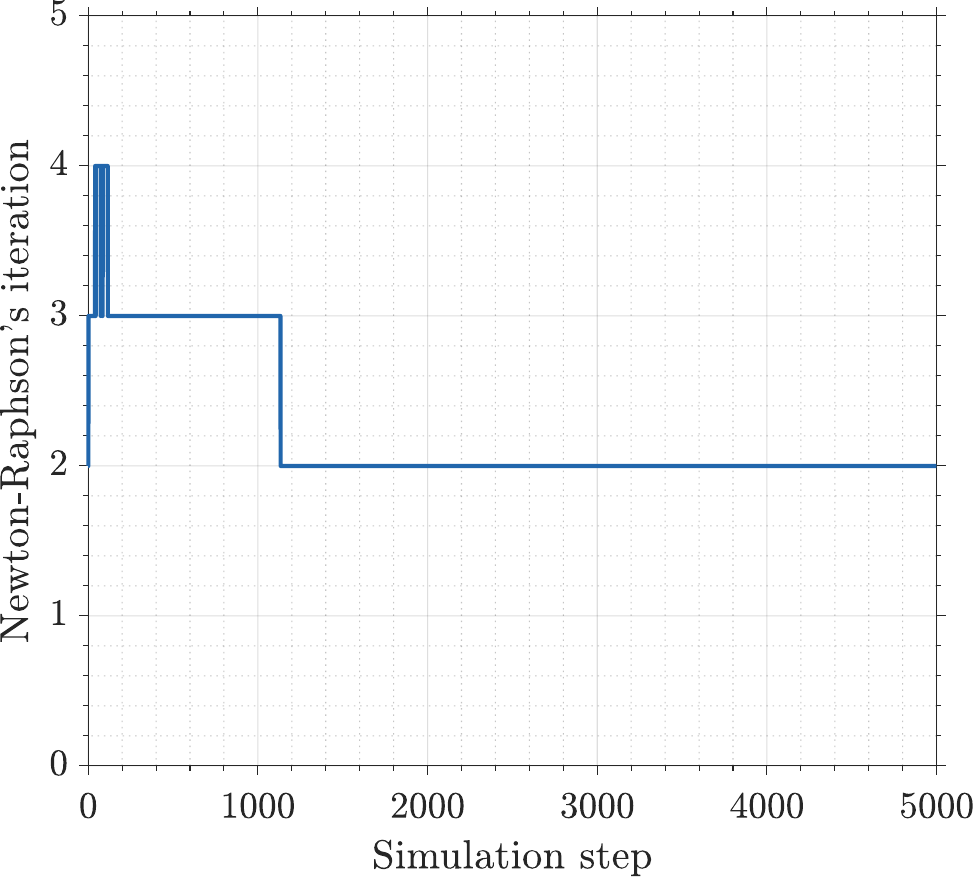}
         \caption{}
         \label{fig:3.1.1_nr_it}
     \end{subfigure}
    \caption{1D verification suite - Case 1: (a) evolution of center of mass in $y$-axis. The theoretical minimum value of the $y$-coordinate is used for comparison in place of an exact analytical solution. (b) Number of Newton-Raphson's iterations per step.}
    \label{fig:3.1.1_com_speed}
\end{figure}

The snapshots of fluid configuration and pressure field for the aforementioned discretization are given in Fig.~\ref{fig:3.1.1_results} for the proposed mixed-VMS formulation compared to the fractional-step approach. Meanwhile, the $y$-coordinate position of the center of mass $y_{cm}$ is plotted in Fig.~\ref{fig:3.1.1_com}, which can be used as a proxy to check the volume conservation. The mass center $\tb x_{cm}$ is computed through particle position averaging as:
\begin{eqnarray}
    \tb x_{cm} = \frac{1}{M_{tot}} \sum_{p=1}^{N_p} m_p \tb x_p\,, \qquad M_{tot} = \sum_{p=1}^{N_p} m_p \,,
\label{eq:center_of_mass}
\end{eqnarray}
where $M_{tot}$ indicates the total mass of the material points in the simulation domain. As can be seen from both Figs.~\ref{fig:3.1.1_results} and \ref{fig:3.1.1_com}, the mixed MPM with VMS stabilization shows superior volume conservation in the long term, compared to the fractional-step approach, which becomes unstable at the later stage of the simulation. The pressure contour and the distribution of material points inside the porous media are also obtained to be more stable and well-distributed without any spurious oscillations. 

The convergence of the proposed formulation and solver is verified for this problem by checking the number of Newton-Raphson's iterations required at each time step, which is provided in Fig.~\ref{fig:3.1.1_nr_it}. The performance of the solver is obtained to be very well, where, in general, the converged solution can be obtained in less than 5 iterations (mostly 2). Furthermore, Table \ref{tab:3.1.1_newton_krylov_iteration} presents the convergence of the relative residual norm and energy norm measured in three different time steps. As can be seen, the quadratic convergence of the residuals can be achieved, whereas the number of Krylov iterations is also obtained to be reasonably small.

\begin{table}[h!]
\centering
\caption{1D verification suite - Case 1: performance of the Newton--Krylov solver of the stabilized mixed MPM solver.}
\label{tab:3.1.1_newton_krylov_iteration}
\small
\begin{tabular}{||c|ccc|ccc|ccc||}
\hline
\multirow{2}{*}{NR it.} & \multicolumn{3}{c|}{Step 1000}                                           & \multicolumn{3}{c|}{Step 2000}                                           & \multicolumn{3}{c||}{Step 5000}                                           \\ \cline{2-10} 
                        & \multicolumn{1}{c|}{$\mathcal{C}^{\mathrm{rel}}$}    & \multicolumn{1}{c|}{$\mathcal{C}^{\mathrm{en}}$}     & Krylov it. & \multicolumn{1}{c|}{$\mathcal{C}^{\mathrm{rel}}$}    & \multicolumn{1}{c|}{$\mathcal{C}^{\mathrm{en}}$}     & Krylov it. & \multicolumn{1}{c|}{$\mathcal{C}^{\mathrm{rel}}$}    & \multicolumn{1}{c|}{$\mathcal{C}^{\mathrm{en}}$}     & Krylov it. \\ \hline \hline
0                       & \multicolumn{1}{c|}{1.0e-00} & \multicolumn{1}{c|}{-}       & 56         & \multicolumn{1}{c|}{1.0e-00} & \multicolumn{1}{c|}{-}       & 79& \multicolumn{1}{c|}{1.0e-00} & \multicolumn{1}{c|}{-}       & 85         \\ \hline
1                       & \multicolumn{1}{c|}{4.8e-05} & \multicolumn{1}{c|}{7342.9}  & 57& \multicolumn{1}{c|}{7.4e-06} & \multicolumn{1}{c|}{139.7}   & 67& \multicolumn{1}{c|}{5.6e-06} & \multicolumn{1}{c|}{8.5}     & 76\\ \hline
2                       & \multicolumn{1}{c|}{1.5e-08} & \multicolumn{1}{c|}{1.2e-04} & 56& \multicolumn{1}{c|}{1.1e-11} & \multicolumn{1}{c|}{6.9e-09} & -          & \multicolumn{1}{c|}{5.6e-13} & \multicolumn{1}{c|}{2.7e-10} & -          \\ \hline
3                       & \multicolumn{1}{c|}{4.4e-12} & \multicolumn{1}{c|}{1.2e-11} & -          & \multicolumn{1}{c|}{-}       & \multicolumn{1}{c|}{-}       & -          & \multicolumn{1}{c|}{-}       & \multicolumn{1}{c|}{-}       & -          \\ \hline
\end{tabular}
\end{table}

\subsubsection{Case 2 -- volume conservation of fluid flowing through porous domain}
\label{subsubsec:1d_case2}

\begin{figure}[h!]
    \centering
     \begin{subfigure}[b]{0.46\textwidth}
         \centering
         \includegraphics[width=\textwidth]{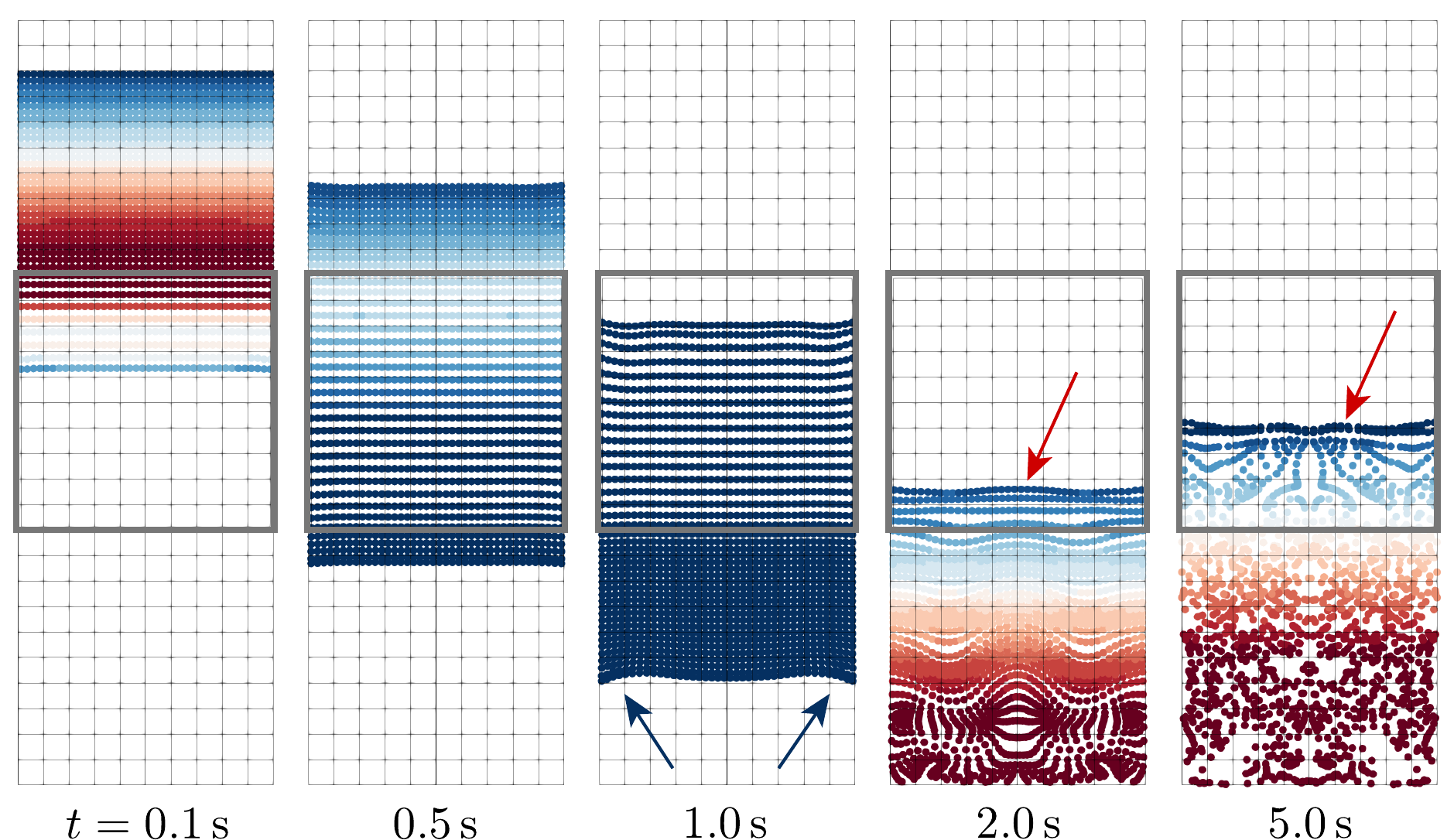}
         \caption{Fractional-step}
         \label{fig:3.1.2_screenshots_fs}
     \end{subfigure}
     \raisebox{0.4cm}{\rule{0.1pt}{4.45cm}}
     \begin{subfigure}[b]{0.46\textwidth}
         \centering
         \includegraphics[width=\textwidth]{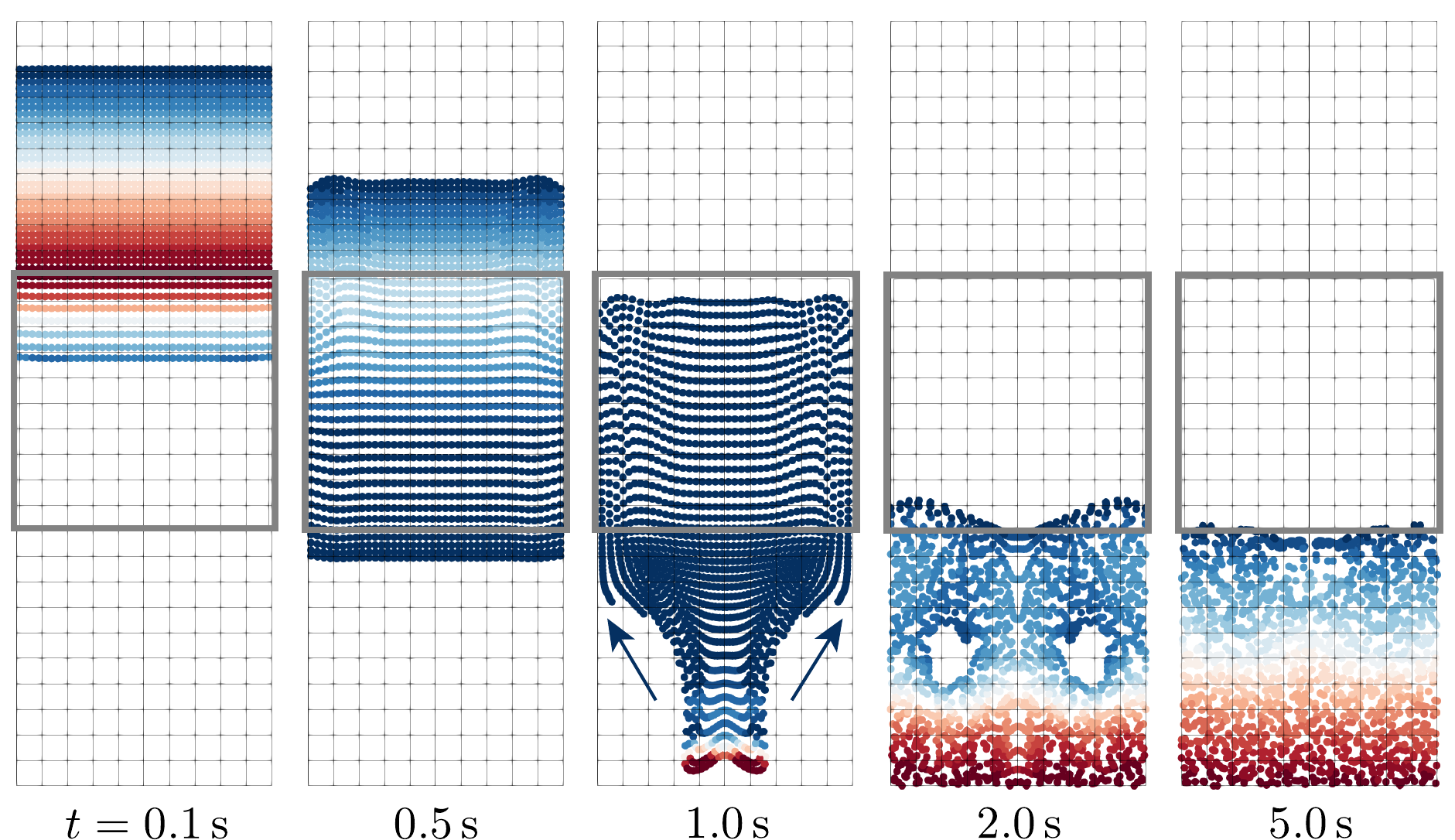}
         \caption{Mixed--VMS}
         \label{fig:3.1.2_screenshots_vms}
     \end{subfigure}
     \begin{subfigure}{0.053\textwidth}
         \centering
         \includegraphics[width=\textwidth]{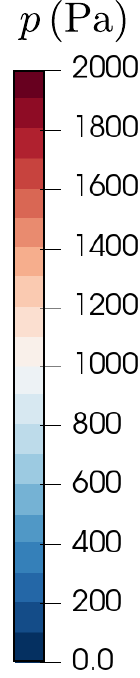}
         \caption*{}
     \end{subfigure}
    \caption{1D verification suite - Case 2: evolution of pressure field at different time snapshots $t=0.1$, 0.5, 1.0, 2.0, and 5.0 s. The red arrows highlight volumetric and pressure errors obtained by the fractional-step method, whereas the dark blue arrows highlight the over-constraint of horizontal movements by the side walls, which prohibit particle separation from the walls.}
    \label{fig:3.1.2_results}
\end{figure}

\begin{figure}[h!]
    \centering
     \begin{subfigure}[b]{0.405\textwidth}
         \centering
         \includegraphics[width=\textwidth]{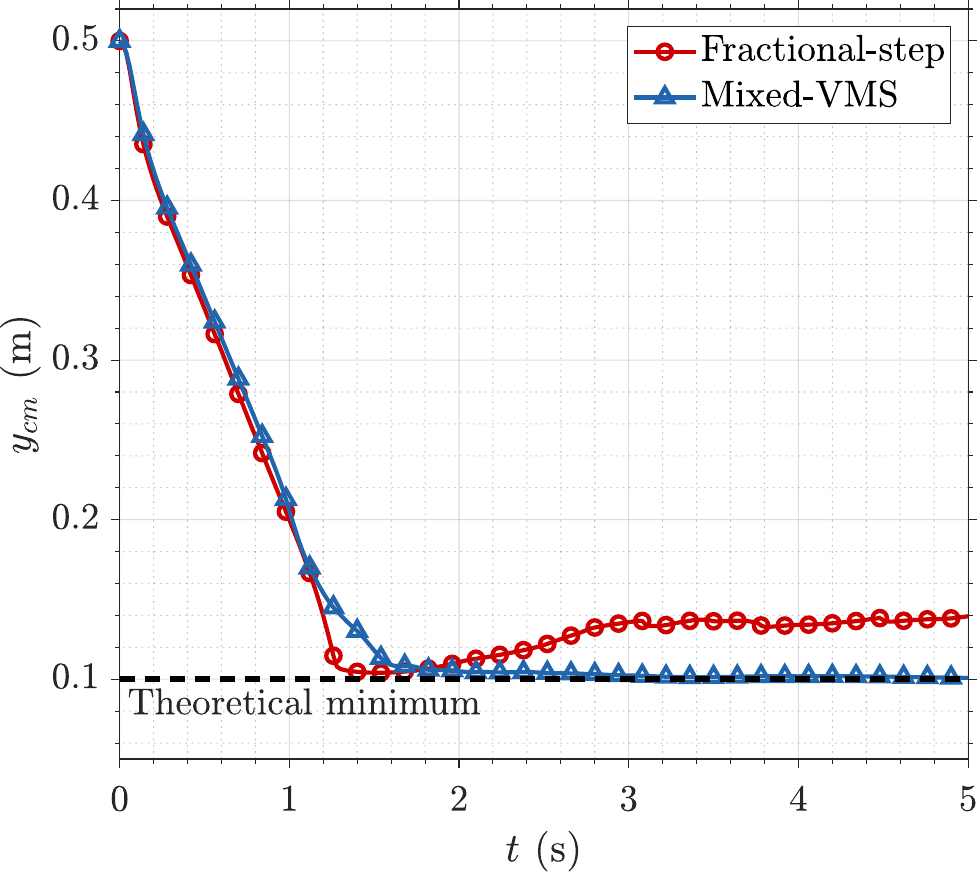}
         \caption{}
         \label{fig:3.1.2_com}
     \end{subfigure}
     \hspace{0.2cm}
     \begin{subfigure}[b]{0.4\textwidth}
         \centering
         \includegraphics[width=\textwidth]{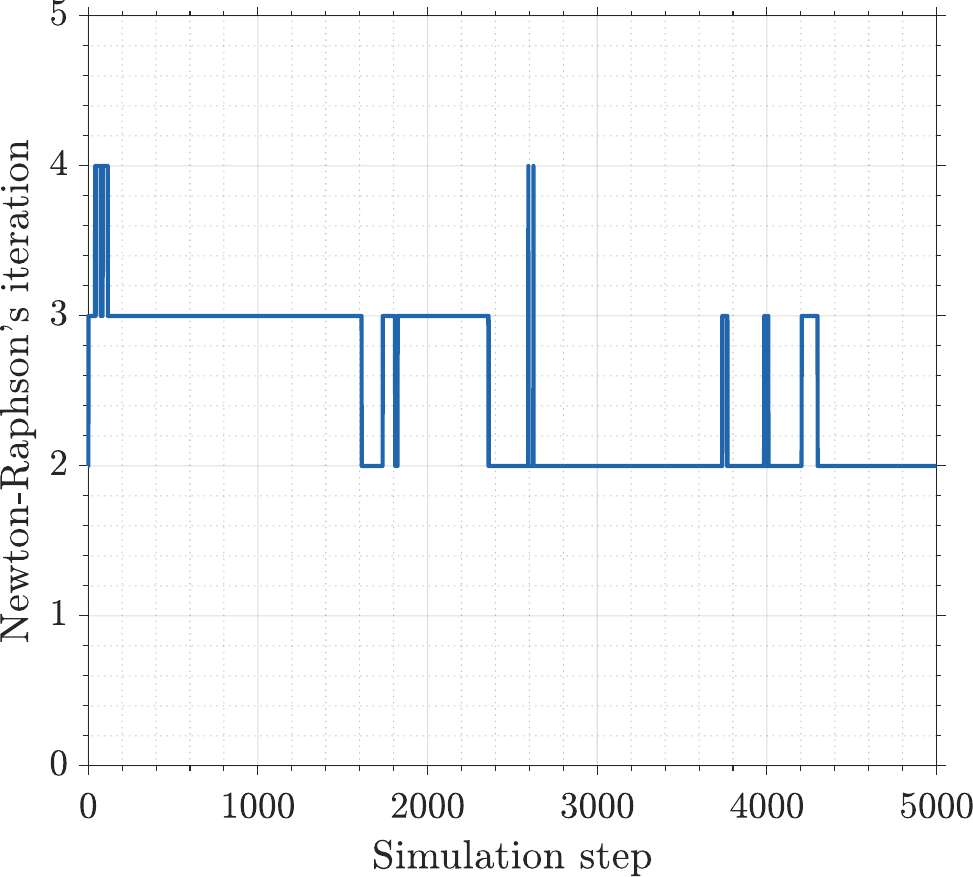}
         \caption{}
         \label{fig:3.1.2_nr_it}
     \end{subfigure}
    \caption{1D verification suite - Case 2: (a) evolution of center of mass in $y$-axis. The theoretical minimum value of the $y$-coordinate is used for comparison in place of an exact analytical solution. (b) Number of Newton-Raphson's iterations per step.}
    \label{fig:3.1.2_com_speed}
\end{figure}

To further confirm volume conservation, the height of the porous domain described in Case 1 is reduced by half \cite{akbari2014modified} (see Fig.~\ref{fig:3.1_model}). In this case, the fluid patch passing through the porous media is expected to fully occupy the free non-porous domain whose size is the same as the original size of the patch. Here, the material and numerical settings are set to be the same as in Case 1. Similar to the previous case, the proposed formulation exhibited significantly improved volume conservation and maintained a stable pressure profile, as evident in Figs.~\ref{fig:3.1.2_results} and \ref{fig:3.1.2_com}. It is worth noting that an over-constraining behavior is observed near the side walls, where particles are restricted from moving away from the walls due to fixity supports in the $x$-axis. However, this numerical artifact did not affect the accuracy of volume conservation. In addition to that, Fig.~\ref{fig:3.1.2_nr_it} and Table \ref{tab:3.1.2_newton_krylov_iteration} validate the performance of the proposed formulation, where the convergence of the nonlinear solver to the specified tolerances can be achieved in less than 5 iterations. It is worth highlighting that, in agreement with the plot presented by Fig.~\ref{fig:3.1.1_nr_it} for Case 1, there is a slight increase in the number of iterations when the fluid enters and exits the porous domain, i.e.~from 2 to 3 or 4 iterations per step.

\begin{table}[h!]
\centering
\caption{1D verification suite - Case 2: performance of the Newton--Krylov solver of the stabilized mixed MPM solver.}
\label{tab:3.1.2_newton_krylov_iteration}
\small
\begin{tabular}{||c|ccc|ccc|ccc||}
\hline
\multirow{2}{*}{NR it.} & \multicolumn{3}{c|}{Step 1000}                                           & \multicolumn{3}{c|}{Step 2000}                                           & \multicolumn{3}{c||}{Step 5000}                                           \\ \cline{2-10} 
                        & \multicolumn{1}{c|}{$\mathcal{C}^{\mathrm{rel}}$}    & \multicolumn{1}{c|}{$\mathcal{C}^{\mathrm{en}}$}     & Krylov it. & \multicolumn{1}{c|}{$\mathcal{C}^{\mathrm{rel}}$}    & \multicolumn{1}{c|}{$\mathcal{C}^{\mathrm{en}}$}     & Krylov it. & \multicolumn{1}{c|}{$\mathcal{C}^{\mathrm{rel}}$}    & \multicolumn{1}{c|}{$\mathcal{C}^{\mathrm{en}}$}     & Krylov it. \\ \hline \hline
0                       & \multicolumn{1}{c|}{1.0e-00} & \multicolumn{1}{c|}{-}       & 43& \multicolumn{1}{c|}{1.0e-00} & \multicolumn{1}{c|}{-}       & 34& \multicolumn{1}{c|}{1.0e-00} & \multicolumn{1}{c|}{-}       & 47\\ \hline
1                       & \multicolumn{1}{c|}{4.4e-05} & \multicolumn{1}{c|}{34592.6}  & 55& \multicolumn{1}{c|}{9.9e-06} & \multicolumn{1}{c|}{60721.4}   & 41& \multicolumn{1}{c|}{2.8e-06} & \multicolumn{1}{c|}{6635.4}     & 55\\ \hline
2                       & \multicolumn{1}{c|}{1.5e-08} & \multicolumn{1}{c|}{1.3e-03} & 51& \multicolumn{1}{c|}{1.2e-11} & \multicolumn{1}{c|}{4.8e-05} & 43& \multicolumn{1}{c|}{2.1e-11} & \multicolumn{1}{c|}{1.3e-07} & -          \\ \hline
3                       & \multicolumn{1}{c|}{4.2e-12} & \multicolumn{1}{c|}{1.8e-10} & -          & \multicolumn{1}{c|}{7.7e-14}       & \multicolumn{1}{c|}{5.5e-13}       & -          & \multicolumn{1}{c|}{-}       & \multicolumn{1}{c|}{-}       & -          \\ \hline
\end{tabular}
\end{table}

\subsubsection{Case 3 -- seepage velocity of Darcian and non-Darcian flow}
\label{subsubsec:1d_case3}

\begin{figure}[h!]
    \centering
     \begin{subfigure}[b]{0.4\textwidth}
         \centering
         \includegraphics[width=\textwidth]{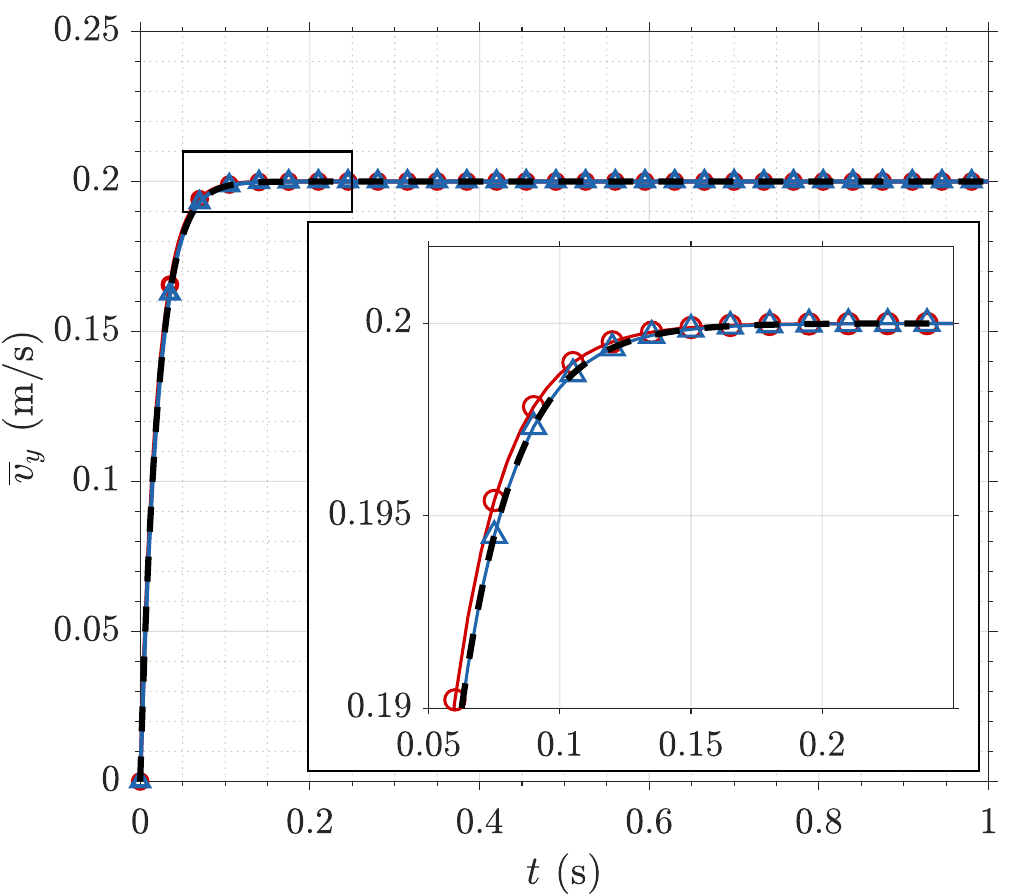}
         \caption{Darcian flow}
         \label{fig:3.1.3_linear}
     \end{subfigure}
     \hspace{0.2cm}
     \begin{subfigure}[b]{0.4\textwidth}
         \centering
         \includegraphics[width=\textwidth]{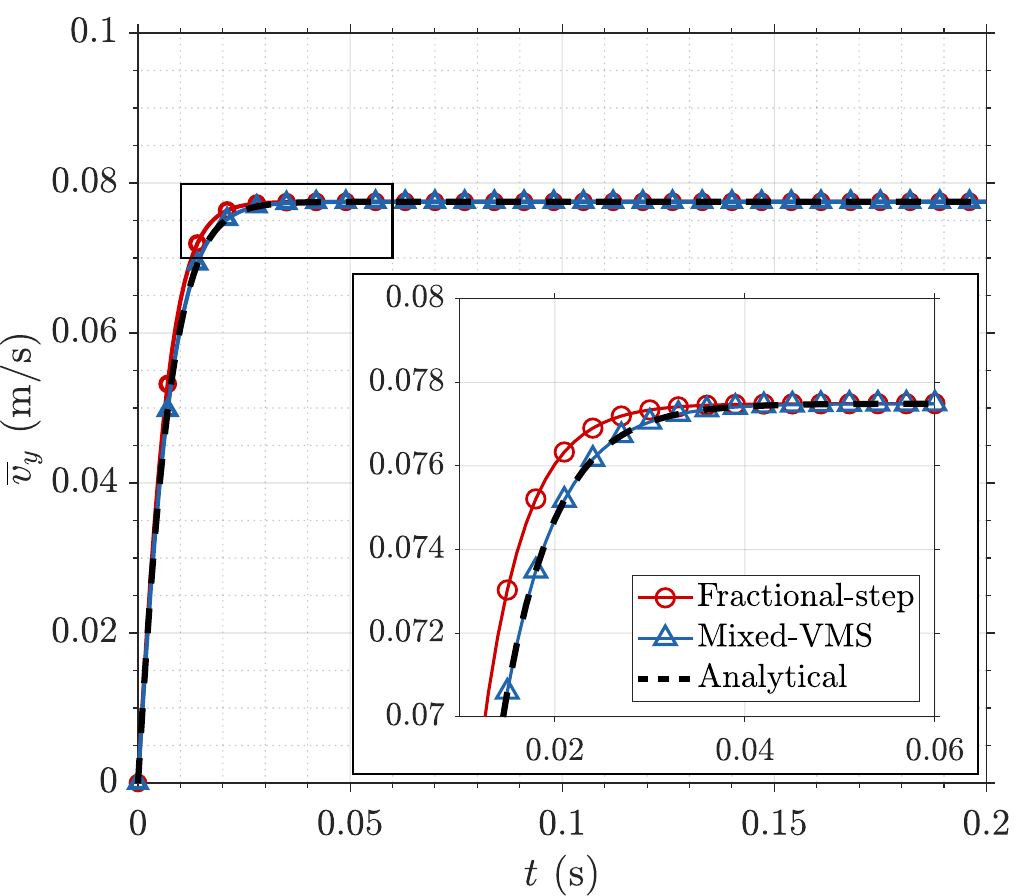}
         \caption{Non-Darcian flow}
         \label{fig:3.1.3_nonlinear}
     \end{subfigure}
    \caption{1D verification suite - Case 3: evolution of averaged vertical velocity $\overline{v}_y$ over time.}
    \label{fig:3.1.3_results}
\end{figure}

In the next verification test, the fluid material points are initialized inside the porous medium with porosity $\theta=0.5$. Here, the material and numerical settings are also the same as the ones considered in Case 1. This test aims to evaluate the kinematic accuracy of fluid motion within porous media, which is governed by the selected drag force law, specifically, the Darcy or Darcy-Forchheimer (non-Darcian) models. The obtained averaged vertical velocity $\overline{v}_y$ is plotted and compared with the semi-analytical solution (see \ref{app:1d_case3_derivation}) in Fig.~\ref{fig:3.1.3_results}. The averaged particle velocity $\overline{v}_y$ is calculated as:
\begin{eqnarray}
    \overline{v}_y = \frac{1}{M_{tot}} \sum_{p=1}^{N_p} m_p {v}_{p,y}\,.
\end{eqnarray}
As it can be observed, the fractional-step approach exhibits a faster rate of velocity increase, attributed to the explicit formulation of the drag force which considered the previous-step velocity field $\dot{\tb u}^{n}$ (c.f.~Eq.~\eqref{eq:predictor}).

\begin{figure}[h!]
    \centering
     \begin{subfigure}[b]{0.4\textwidth}
         \centering
         \includegraphics[width=\textwidth]{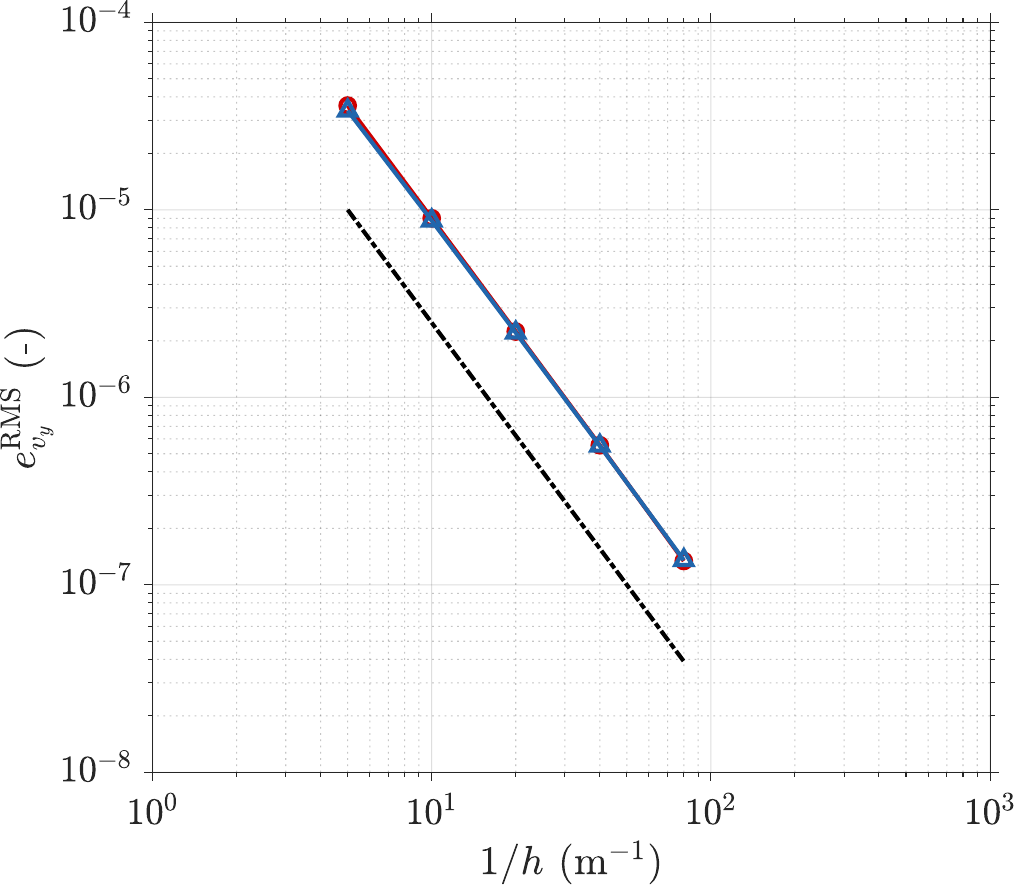}
         \caption{Darcian flow}
         \label{fig:3.1.3_convergence_linear}
     \end{subfigure}
     \hspace{0.2cm}
     \begin{subfigure}[b]{0.4\textwidth}
         \centering
         \includegraphics[width=\textwidth]{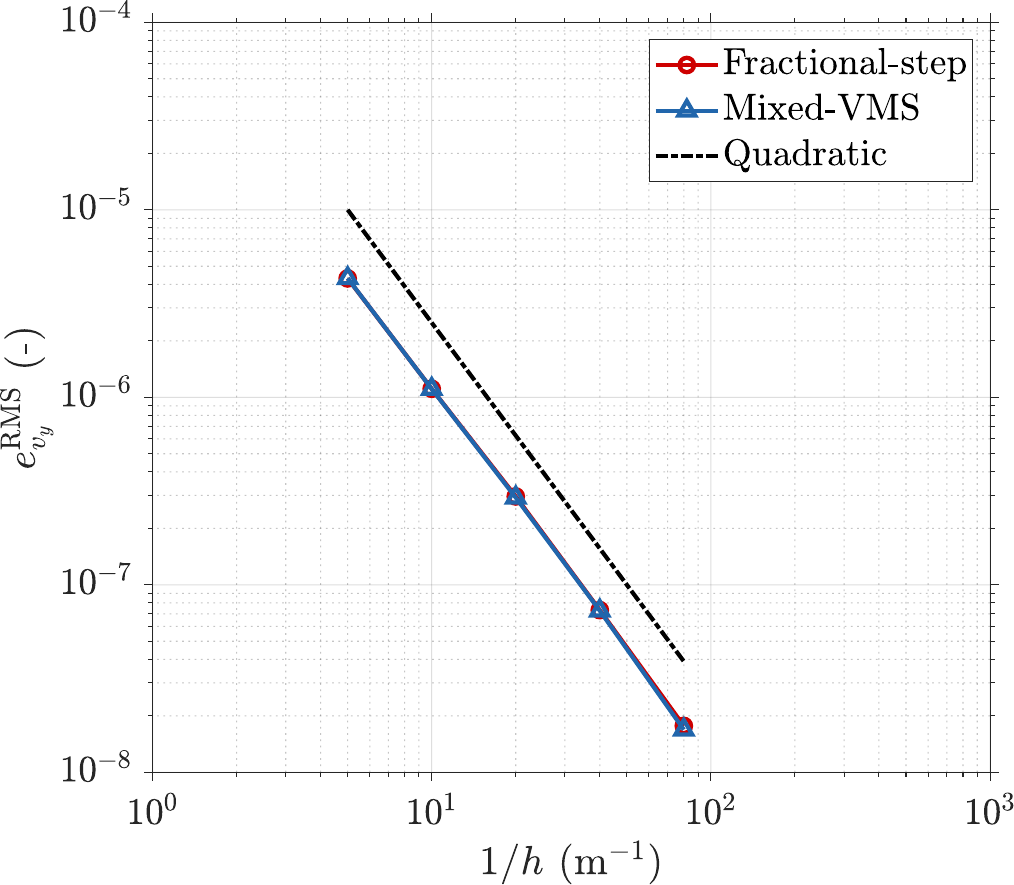}
         \caption{Non-Darcian flow}
         \label{fig:3.1.3_convergence_nonlinear}
     \end{subfigure}
    \caption{1D verification suite - Case 3: convergence plot of vertical velocity with mesh refinement measured at $t=0.2$ s.}
    \label{fig:3.1.3_convergence}
\end{figure}

 Even though the fractional-step method accelerates faster, the steady-state velocity is obtained to be equally accurate to the mixed formulation. The root-mean-squared (RMS) error of the steady-state particle velocity is plotted in Fig.~\ref{fig:3.1.3_convergence} along with mesh refinement, with $h=L/5\sim L/80$. The figures confirmed the quadratic convergence of the formulation for the seepage flow, where there is no porosity variation. The plotted RMS error can be computed as follows:
\begin{eqnarray}
    e^{\rm RMS}_{v} = \sqrt{ \frac{1}{N_p} \sum_{I=1}^{N_p} \left(\frac{v_p - v_a}{v_a} \right)^2 } \,,
\end{eqnarray}
where $N_p$ and $v_a$ are the total number of material points and analytical velocity at steady state, respectively.

\subsubsection{Case 4 -- transient analysis of fluid entering porous media}
\label{subsubsec:1d_case4}

The objective of this verification test is to analyze the transient (or unsteady) behavior of fluid entering a porous media in a one-dimensional setting. As depicted in Fig.~\ref{fig:3.1_model}, a fluid column with initial height $h_0=0.3$ m is dropped at rest right on top of a porous media with porosity $\theta=0.5$ and hydraulic conductivity $K=0.1$ m/s. As the fluid infiltrates the porous media, its height, denoted as $h_{ff}$, gradually decreases. At the same time, the pore fluid height increases with time, which is denoted as $h_{fp}$. Due to mass conservation, the lower limit of $h_{ff}$ and the upper limit of $h_{fp}$ are 0 and 0.6 m, respectively. Here, the pressure at the exact position of the sharp interface is denoted as $p_i$, which also changes with time as the fluid moves through the interface.

The MPM mesh size for this problem is set to be $h=h_0/30=0.01$ m, while the PPC and time step are selected to be $4\times 4$ and $\Delta t =0.001$ s. As mentioned earlier, to prohibit any horizontal motion that may introduce numerical errors, one element thickness is considered for this case, where the horizontal displacement of the left and right walls are fixed. In order to verify the MPM simulation results, a 1D analytical solution for this problem is derived, which will be elaborated further in \ref{app:1d_case4_derivation}.

\begin{figure}[h!]
    \centering
     \begin{subfigure}[b]{0.315\textwidth}
         \centering
         \includegraphics[width=\textwidth]{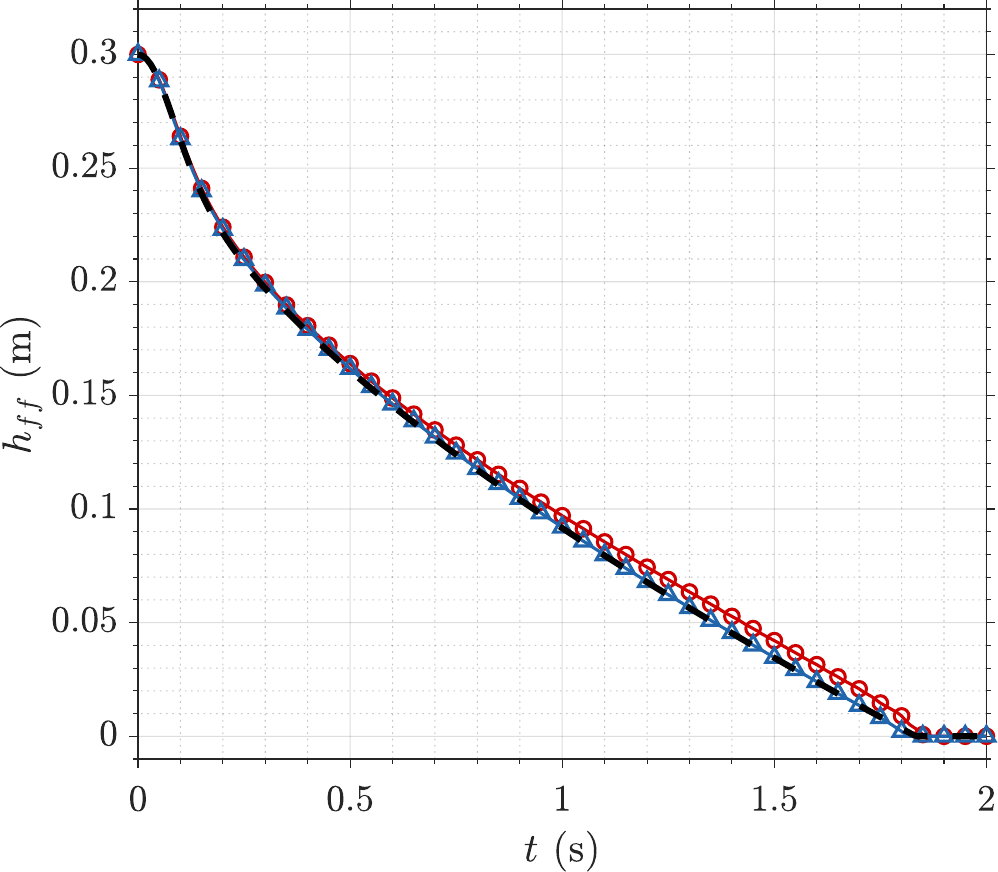}
         \caption*{}
         \label{fig:3.1.4_hff_linear}
     \end{subfigure}
     \hspace{0.1cm}
     \begin{subfigure}[b]{0.315\textwidth}
         \centering
         \includegraphics[width=\textwidth]{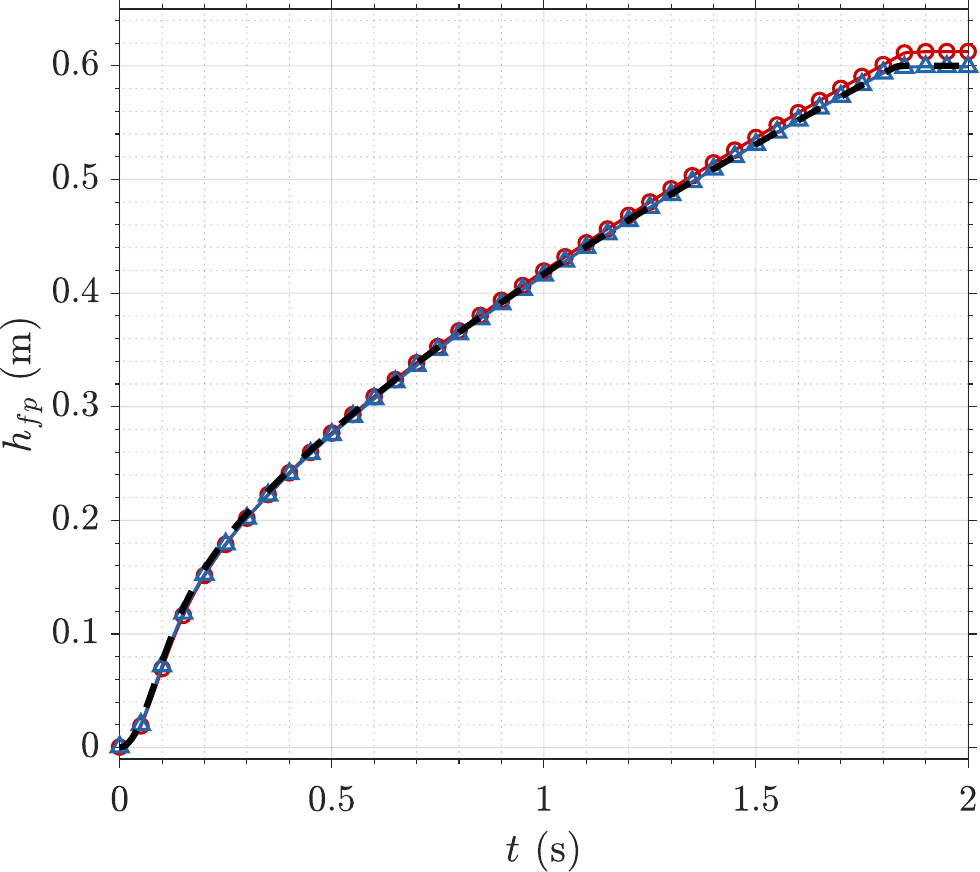}
         \caption{Darcian flow}
         \label{fig:3.1.4_hfp_linear}
     \end{subfigure}
     \hspace{0.1cm}
     \begin{subfigure}[b]{0.32\textwidth}
         \centering
         \includegraphics[width=\textwidth]{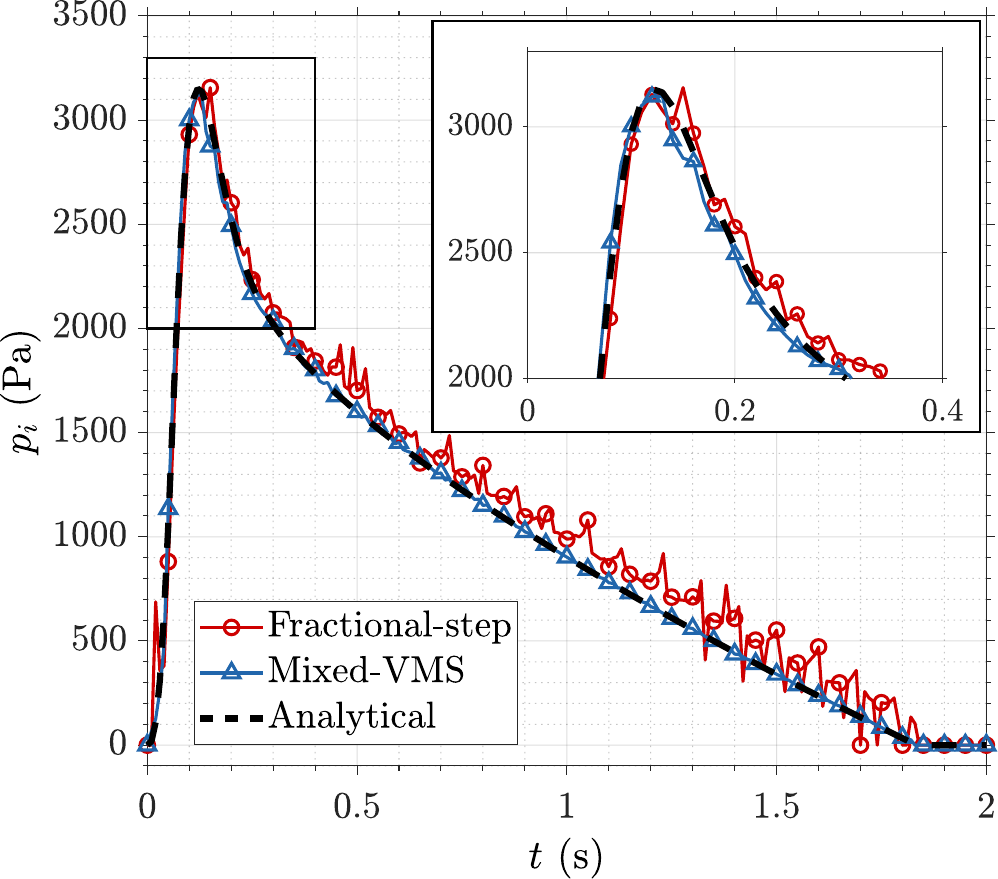}
         \caption*{}
         \label{fig:3.1.4_pi_linear}
     \end{subfigure}\\
     \vspace{0.3cm}
     \begin{subfigure}[b]{0.32\textwidth}
         \centering
         \includegraphics[width=\textwidth]{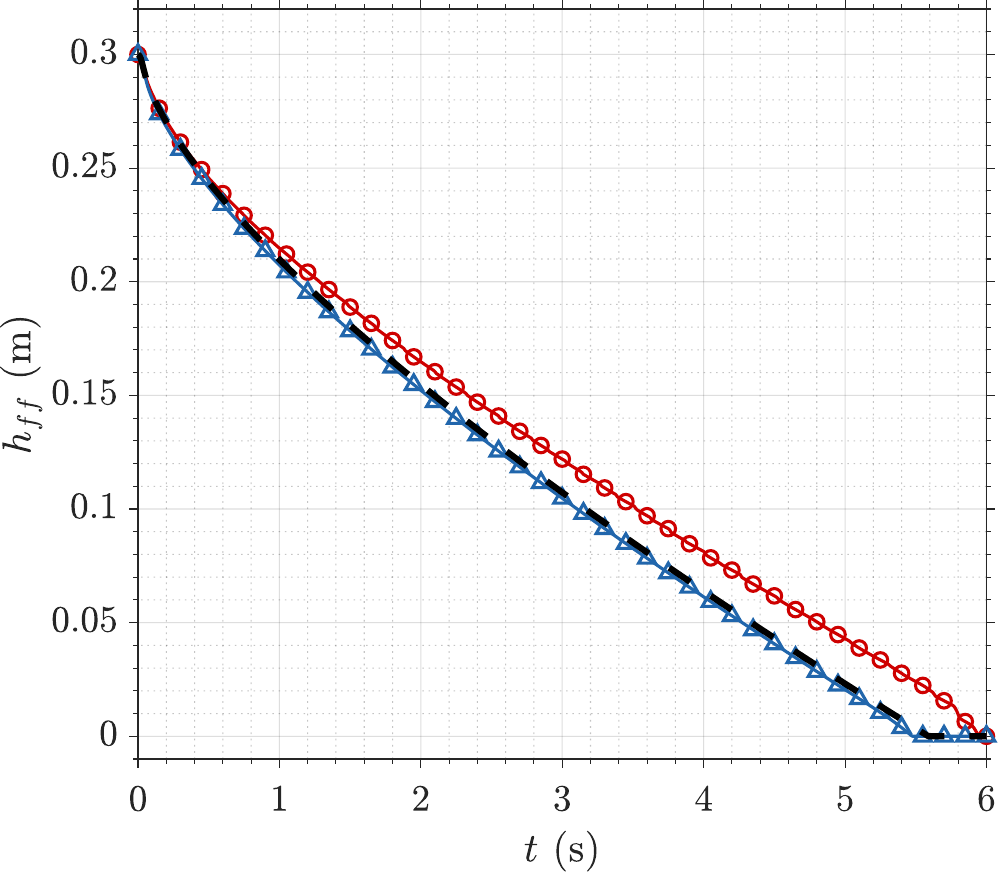}
         \caption*{}
         \label{fig:3.1.4_hff_nonlinear}
     \end{subfigure}
     \hspace{0.1cm}
     \begin{subfigure}[b]{0.315\textwidth}
         \centering
         \includegraphics[width=\textwidth]{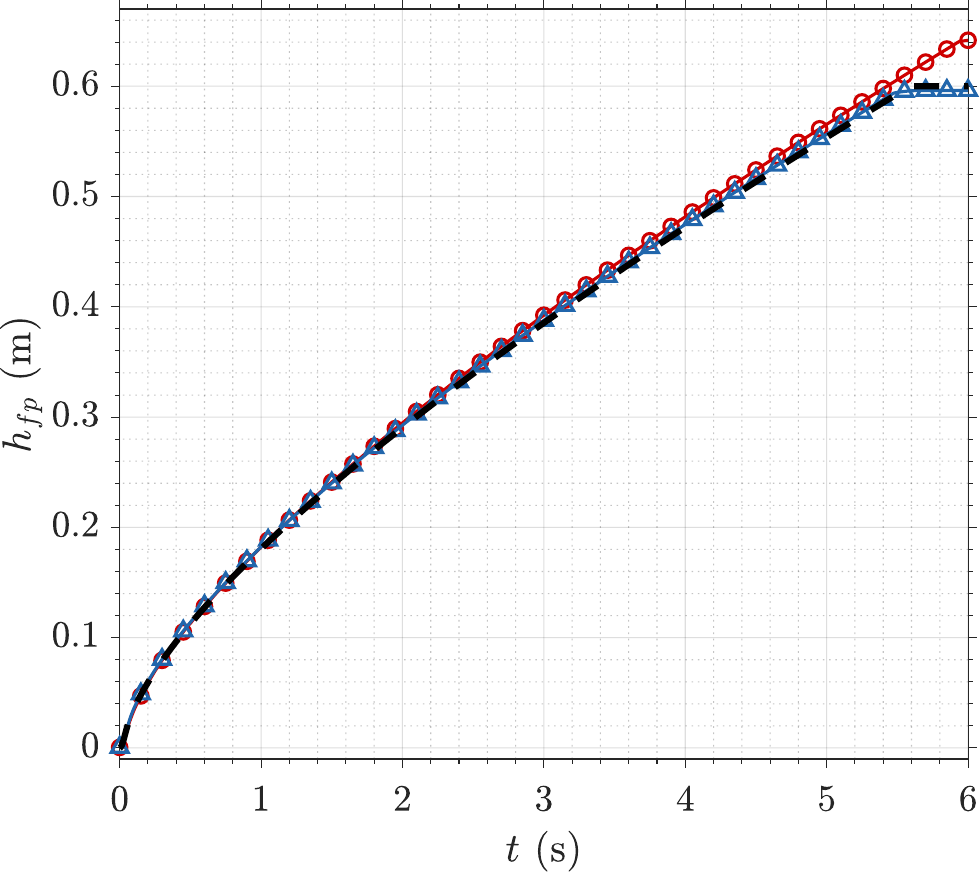}
         \caption{Non-Darcian flow}
         \label{fig:3.1.4_hfp_nonlinear}
     \end{subfigure}
     \hspace{0.1cm}
     \begin{subfigure}[b]{0.32\textwidth}
         \centering
         \includegraphics[width=\textwidth]{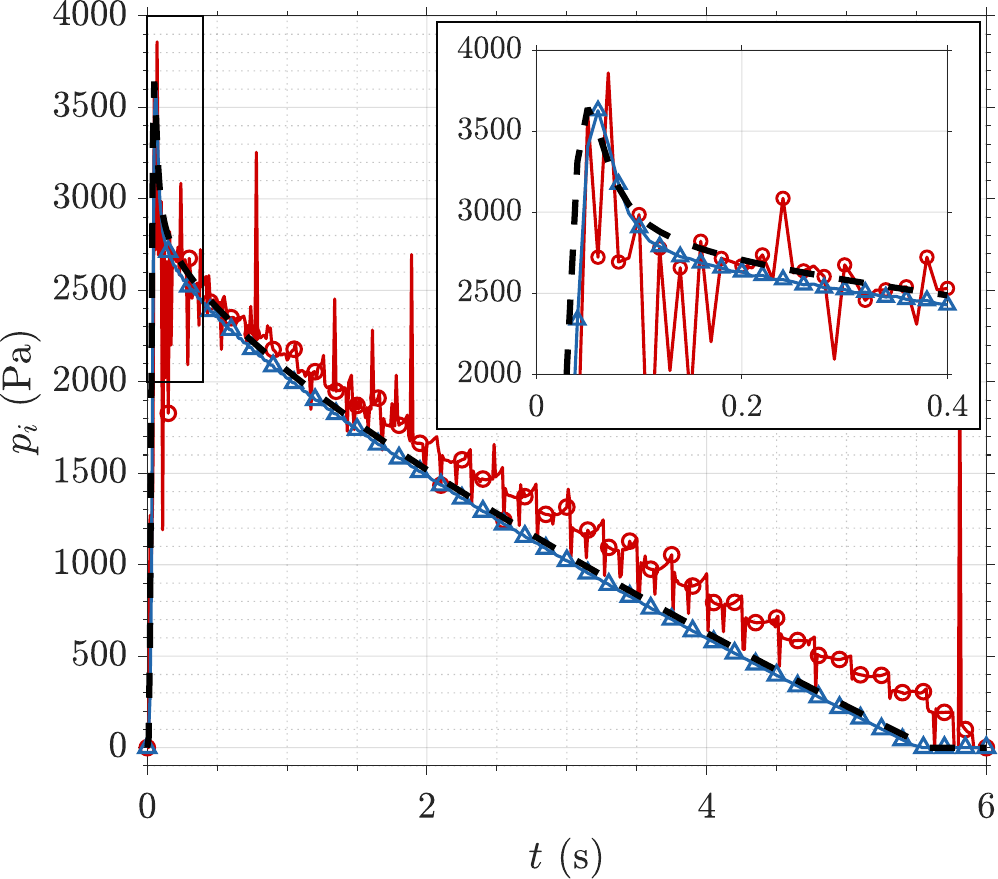}
         \caption*{}
         \label{fig:3.1.4_pi_nonlinear}
     \end{subfigure}
    \caption{1D verification suite - Case 4: evolution of free-fluid height $h_{ff}$, pore-fluid height $h_{fp}$, and interface pressure $p_i$.}
    \label{fig:3.1.4_results}
\end{figure}

The evolution of free-fluid height $h_{ff}$, pore-fluid height $h_{fp}$, and the interface pressure $p_{i}$ over time are plotted in Fig.~\ref{fig:3.1.4_results} along with the derived semi-analytical solutions. Here, the fluid heights are measured from the front and back-most material points by taking into account the particle size $h_p$ computed from particle volume. Meanwhile, the plotted interface pressure is obtained from the nodal pressure value at the exact location of the sharp porous interface. As can be seen from the results, the mixed MPM predicts the evolution of heights and interface pressure significantly more accurately compared to the fractional-step approach, where the differences between the two approaches are more apparent in the non-Darcian flow case. Specifically, the mixed MPM can reduce significantly spurious pressure oscillations that are shown by the fractional-step approach. These oscillations are primarily caused by the \textit{binary pressure imposition error}\cite{chandra2023stabilized}, inherited from the free-surface node detection and pressure imposition required by the fractional-step method.

\subsubsection{Case 5 -- transient analysis of fluid exiting porous media}
\label{subsubsec:1d_case5}

The next case examines the condition opposite to Case 4, where we analyze the transient behavior of fluid exiting one-dimensional porous media, as shown in Fig.~\ref{fig:3.1_model}. A pore-fluid column with an initial height of $h_0 = 0.3$ m is released at rest at the end of the porous domain. The porosity and hydraulic conductivity remain the same as in Case 4, with $\theta = 0.5$ and $K = 0.1$ m/s. Over time, the height of the pore fluid, denoted as $h_{fp}$, gradually decreases, while the height of the free fluid below the porous interface, $h_{ff}$, increases. Mass conservation sets the lower limit for $h_{fp}$ at 0 m and the upper limit for $h_{ff}$ at 0.15 m.

\begin{figure}[h!]
    \centering
     \begin{subfigure}[b]{0.32\textwidth}
         \centering
         \includegraphics[width=\textwidth]{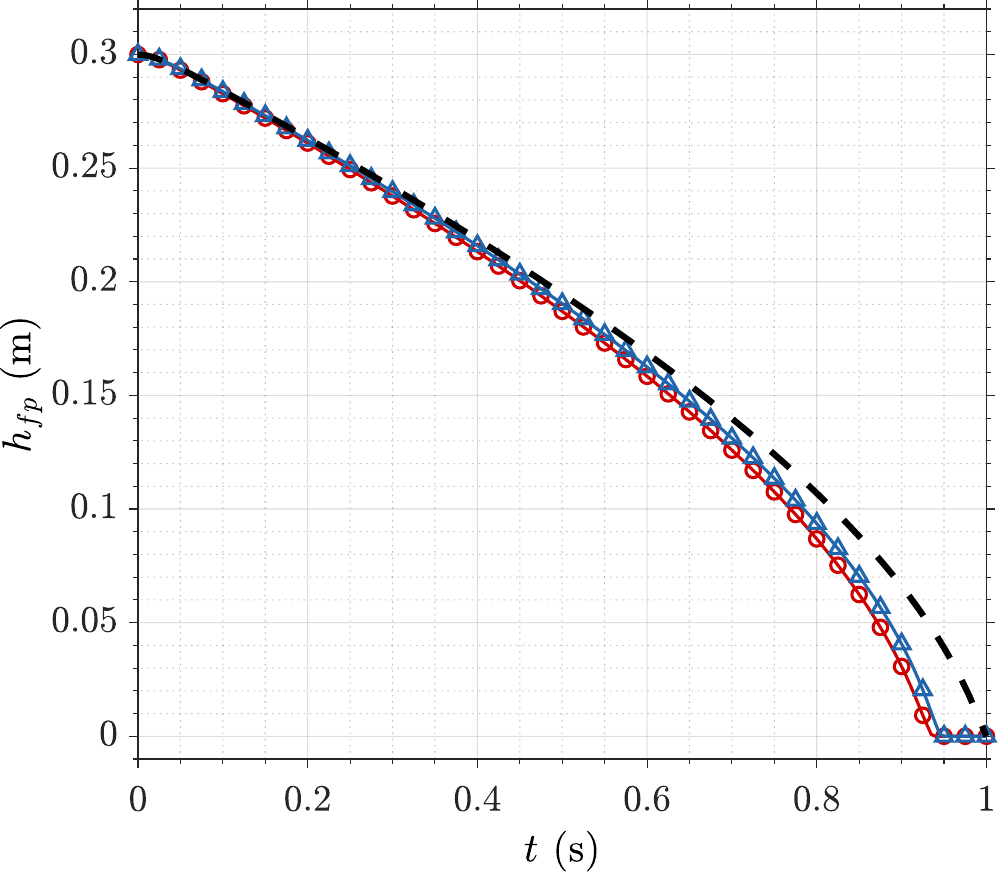}
         \caption*{}
         \label{fig:3.1.5_hff_linear}
     \end{subfigure}
     \hspace{0.1cm}
     \begin{subfigure}[b]{0.32\textwidth}
         \centering
         \includegraphics[width=\textwidth]{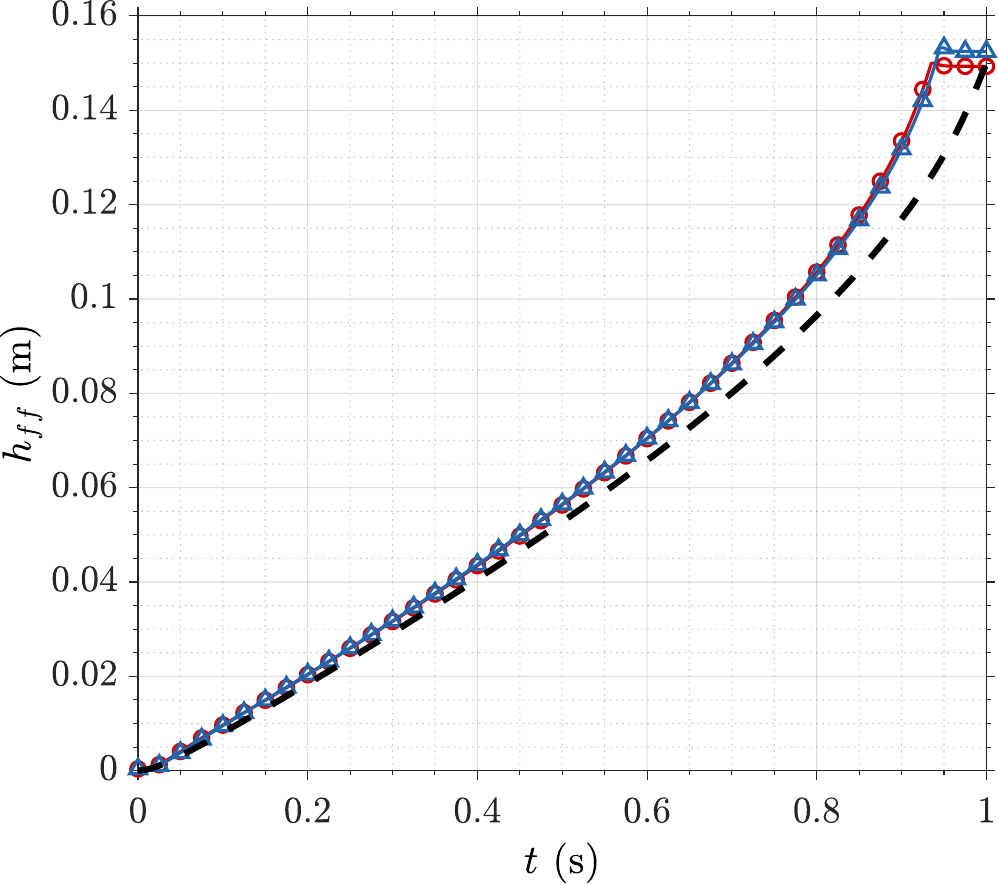}
         \caption{Darcian flow}
         \label{fig:3.1.5_hfp_linear}
     \end{subfigure}
     \hspace{0.1cm}
     \begin{subfigure}[b]{0.32\textwidth}
         \centering
         \includegraphics[width=\textwidth]{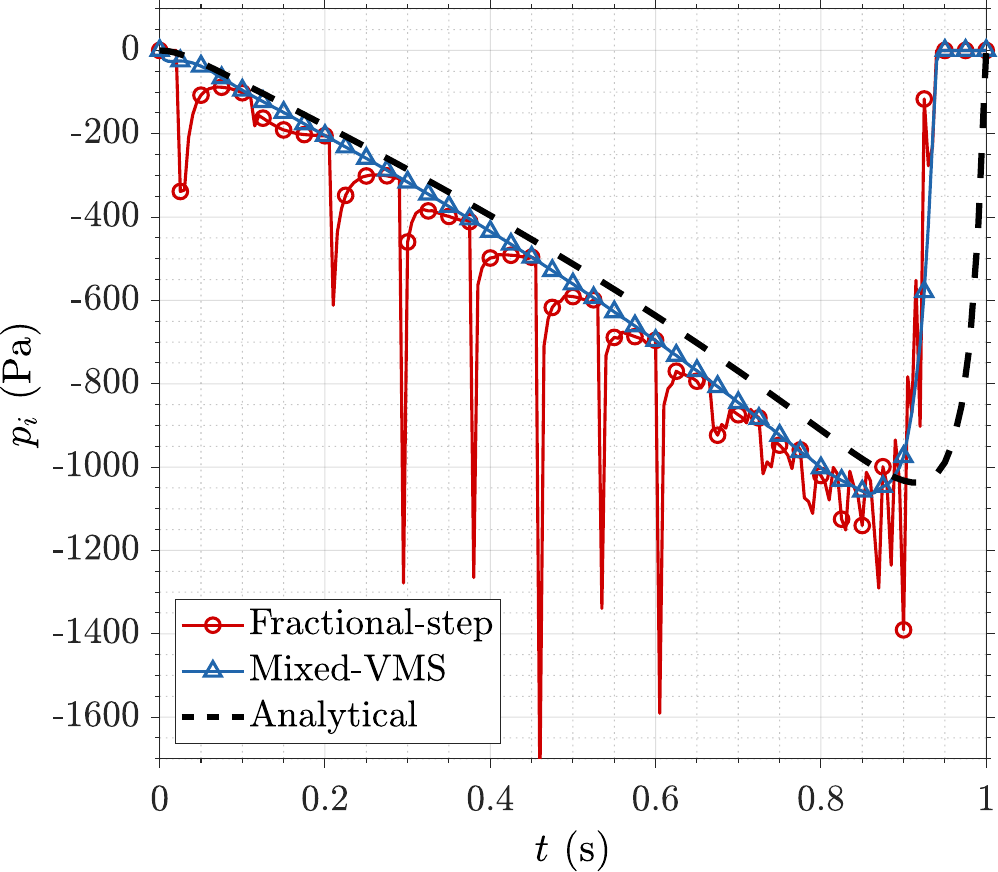}
         \caption*{}
         \label{fig:3.1.5_pi_linear}
     \end{subfigure}\\
     \vspace{0.3cm}
     \begin{subfigure}[b]{0.32\textwidth}
         \centering
         \includegraphics[width=\textwidth]{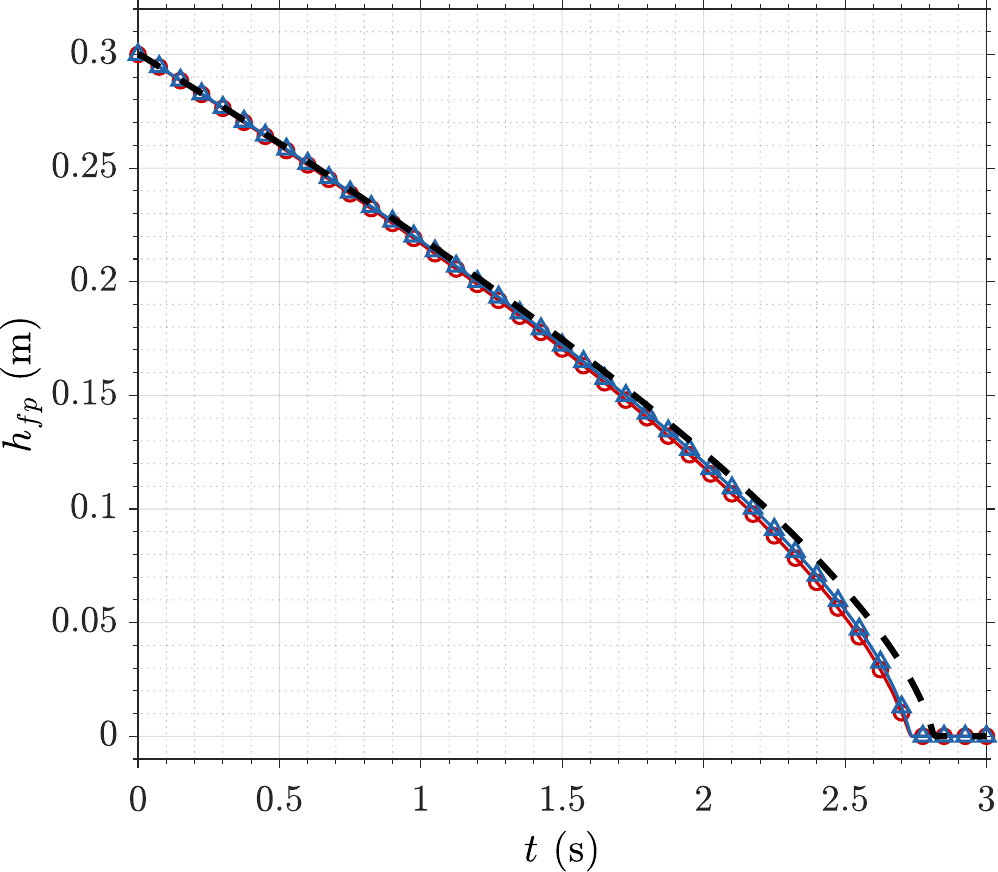}
         \caption*{}
         \label{fig:3.1.5_hff_nonlinear}
     \end{subfigure}
     \hspace{0.1cm}
     \begin{subfigure}[b]{0.32\textwidth}
         \centering
         \includegraphics[width=\textwidth]{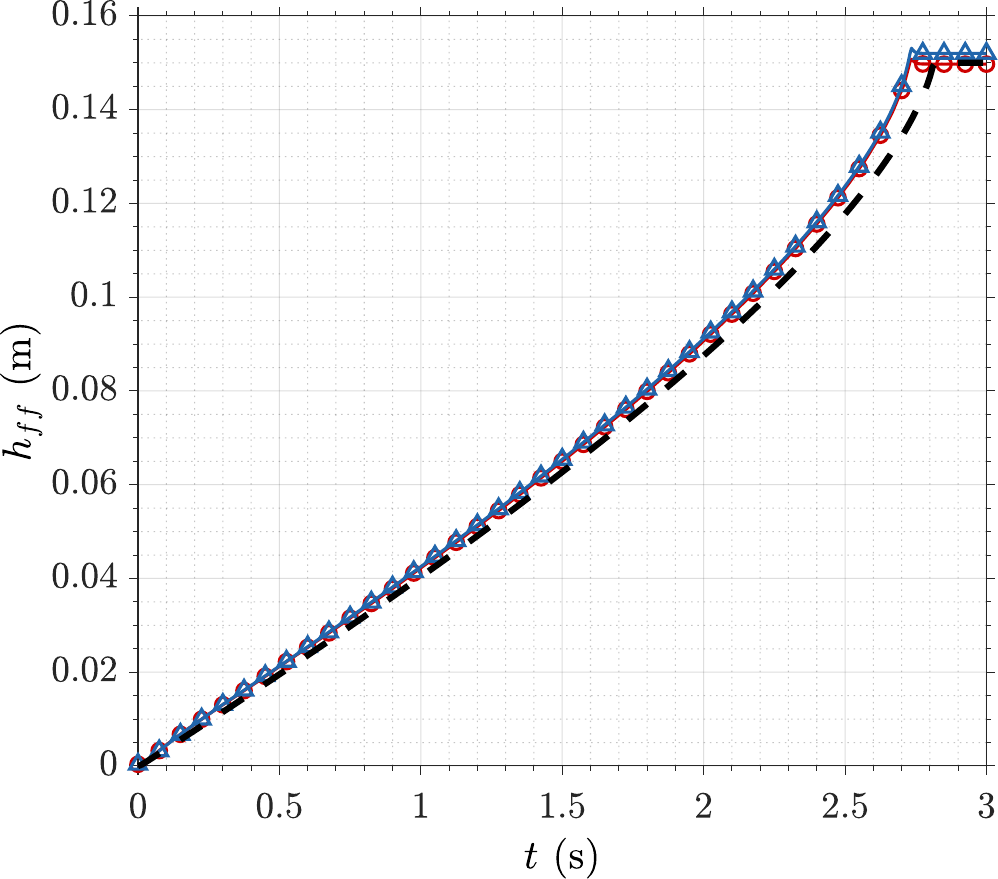}
         \caption{Non-Darcian flow}
         \label{fig:3.1.5_hfp_nonlinear}
     \end{subfigure}
     \hspace{0.1cm}
     \begin{subfigure}[b]{0.32\textwidth}
         \centering
         \includegraphics[width=\textwidth]{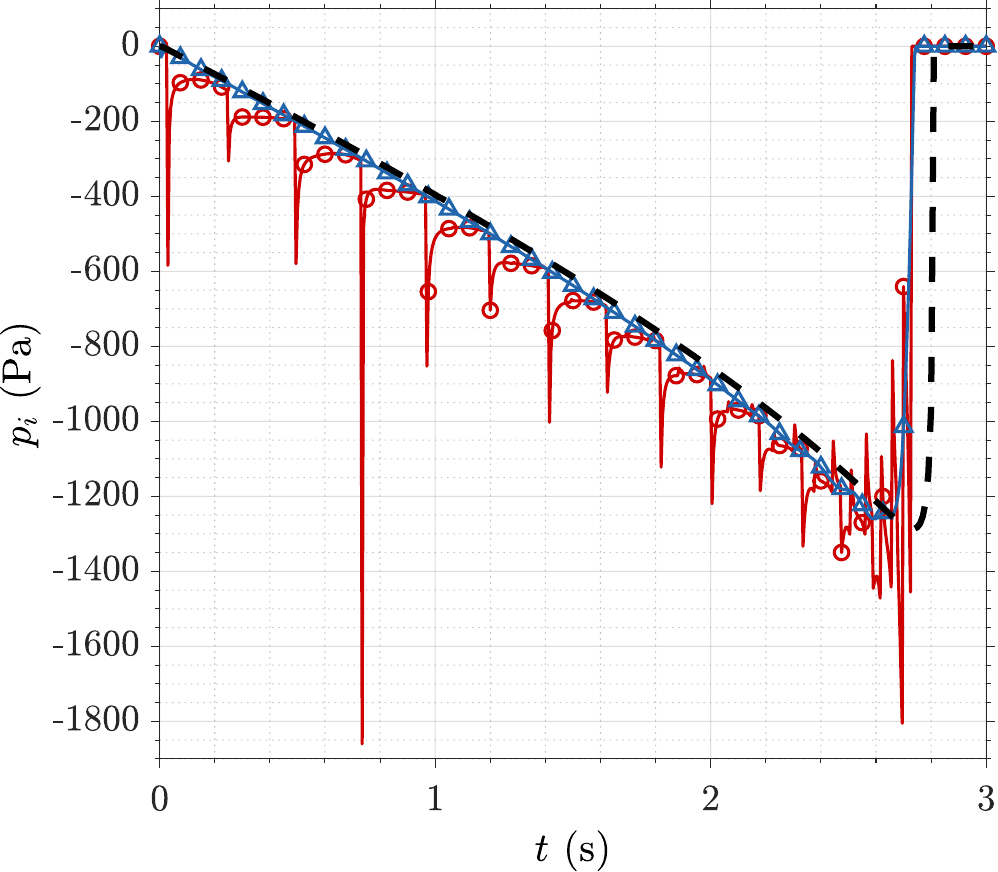}
         \caption*{}
         \label{fig:3.1.5_pi_nonlinear}
     \end{subfigure}
    \caption{1D verification suite - Case 5: evolution of pore-fluid height $h_{fp}$, free-fluid height $h_{ff}$, and interface pressure $p_i$.}
    \label{fig:3.1.5_results}
\end{figure}

Fig.~\ref{fig:3.1.5_results} plots the evolution of the fluid heights, $h_{fp}$, $h_{ff}$, and the interface pressure $p_{i}$ over time, in comparison to the fractional-step approach and the semi-analytical solution. For this test, the differential equation, as derived in \ref{app:1d_case5_derivation}, is found to be stiff, exhibiting a rapid acceleration change in a very short amount of time. In conjunction with this issue, the MPM solutions provided in Fig.~\ref{fig:3.1.5_results} exhibit noticeable inaccuracies, particularly in the later stages of the simulation. Notably, as the nonlinear Darcy flow exhibits a larger jerk (see Fig.~\ref{fig:C.3_anal_jy}), a reduction in the time step to $\Delta t=0.00025$ s, one-fourth of the time step used in the linear Darcy model, becomes necessary to obtain a considerably accurate solution. While the interface pressure obtained by the mixed MPM is significantly smoother in comparison to the fractional-step approach, it is observed that the pressure increase, from the minimum value to zero, occurs at an earlier time compared to the semi-analytical solution. This limitation is majorly caused by our time integration scheme, which employs a single-step Newmark-$\beta$ method, that is not well-suited for solving problems with stiff ODE.

\begin{figure}[h!]
    \centering
     \begin{subfigure}[b]{0.4\textwidth}
         \centering
         \includegraphics[width=\textwidth]{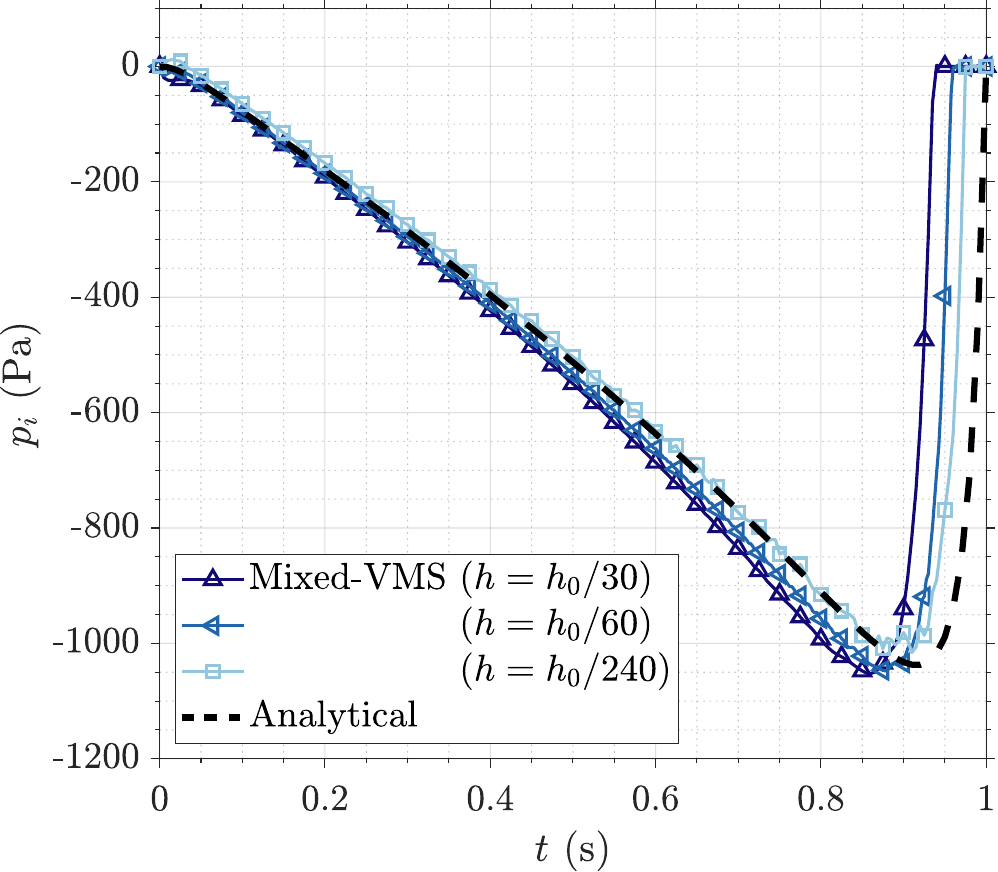}
         \caption{}
         \label{fig:3.1.5_mesh_convergence}
     \end{subfigure}
     \hspace{0.2cm}
     \begin{subfigure}[b]{0.4\textwidth}
         \centering
         \includegraphics[width=\textwidth]{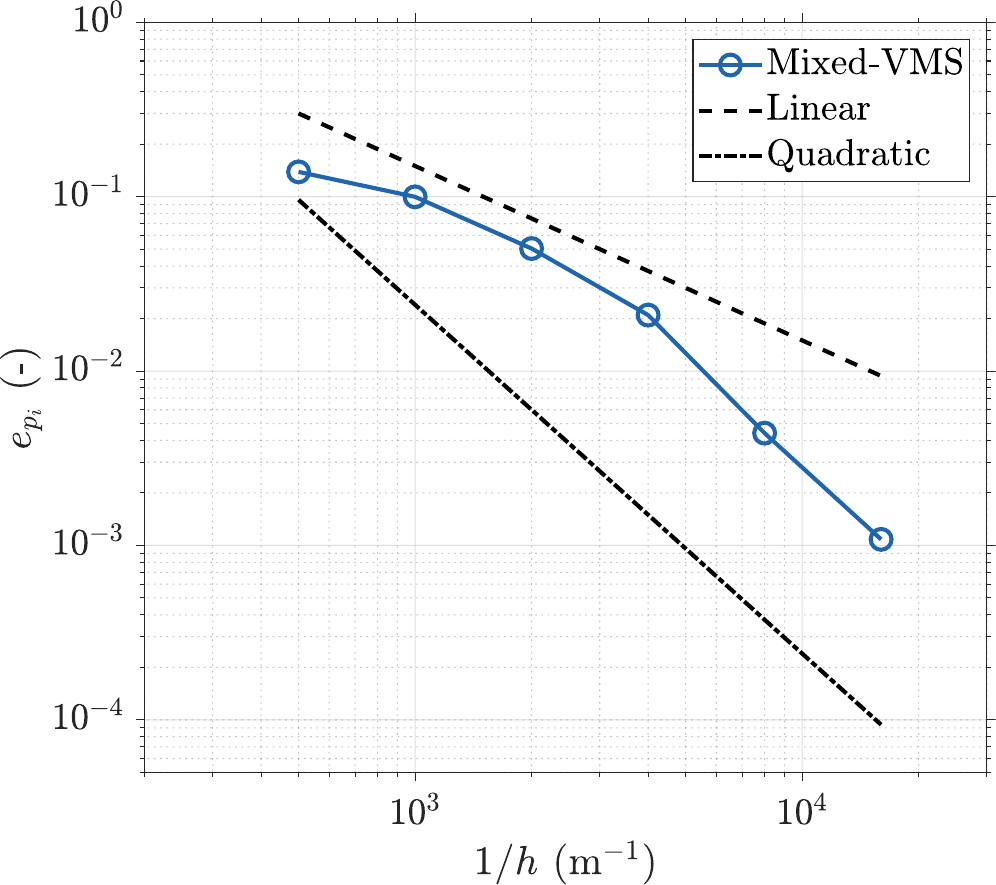}
         \caption{}
         \label{fig:3.1.5_mesh_convergence_rate}
     \end{subfigure}
    \caption{1D verification suite - Case 5: convergence of interface pressure $p_i$ with respect to mesh refinement.}
    \label{fig:3.1.5_mesh_convergence_profile}
\end{figure}

A mesh-refinement study is conducted to assess the convergence trend of the mixed MPM solution. The study is not performed for the fractional-step approach due to its tendency to exhibit spurious pressure oscillations. Consequently, there is concern that the results may not accurately represent a meaningful trend. Here, simulations with different cell sizes, ranging from $h=h_0/15\sim h_0/480$, are performed for the linear Darcy model, while keeping the PPC and ratio of $h/\Delta t$ constant. Fig.~\ref{fig:3.1.5_mesh_convergence_profile} presents the converging interface pressure $p_i$ along with mesh refinement. As highlighted in Fig.~\ref{fig:3.1.5_mesh_convergence}, the interface pressure converges towards the analytical solution as the mesh is refined. Moreover, in Fig.~\ref{fig:3.1.5_mesh_convergence_rate}, the convergence rate of the interface pressure error $e_{p_i}$ at $t=0.5$ s is presented. This error is defined as a single-point error using the formula:
\begin{eqnarray}
    e_{p_i} := \left\vert\frac{p_i - p_{i,a}}{p_{i,a}}\right\vert\,,
\end{eqnarray}
where $p_{i,a}$ represents the analytically computed interface pressure value at the same time frame. The average convergence slope is approximately 1.42. The reduction in convergence rate from quadratic is mainly attributed to the presence of the porous interface. As the blurred interface is defined within the background grids, its width decreases with mesh refinement. Around the interface, over a region of size $\mathcal{O}(\varepsilon)\simeq \mathcal{O}(h)$, the underlying equations are slightly modified, resulting in a lower convergence rate. A similar tendency has been observed previously by other researchers working on the blurred interface, e.g.~\cite{stoter2017diffuse, rycroft2020reference}. 

\subsubsection{Case 6 -- transient analysis of flow across two porous media with different porosity and permeability}
\label{subsubsec:1d_case6}

The final verification test in the 1D-gravity-driven suite involves a transient analysis of pore fluid flow through two distinct porous media. These two porous media are characterized by differing values of porosity and permeability. Two sets of material parameters are examined, which details are listed in Table \ref{tab:3.1.6_cases}. The first set of parameters represents a highly permeable porous medium, falling within a similar range as those in earlier cases. In contrast, the second set simulates realistic soil conditions, specifically sand and sandy loam for materials 1 and 2, respectively\cite{rawls1983green}. This set of porous media exhibits a lower permeability, which is chosen to assess the long-term numerical stability and accuracy of the proposed solver with a finite size of time increment. As indicated in the table, the time step is set as $\Delta t = 1$ s, which is considerably large, given a mesh size of $h = 0.01$ m. The simulated time for the Set 2 extends up to 7 hours. Here, both sets assume the Darcy-Forchheimer drag force model with Ergun's drag force parameters.

\begin{table}[h!]
\centering
\caption{1D verification suite - Case 6: different material and discretization parameters for different sets.}
\label{tab:3.1.6_cases}
\begin{tabular}{||c|c|c|c|c|c|c||}
\hline
      Set & $\theta_1$ (-) & $K_1$ (m/s) & $\theta_2$ (-) & $K_2$ (m/s) & $\Delta t$ (s) & Simulated time \\ \hline \hline
1 & 0.5             & 0.1               & 0.3              & 0.01  & 0.001 & 13 s     \\ \hline
2 & 0.437             & 3.27 $\times 10^{-5}$               & 0.453              & 3.03 $\times 10^{-6}$ & 1 & 7 h      \\ \hline
\end{tabular}
\end{table}

\begin{figure}[h!]
    \centering
     \begin{subfigure}[b]{0.32\textwidth}
         \centering
         \includegraphics[width=\textwidth]{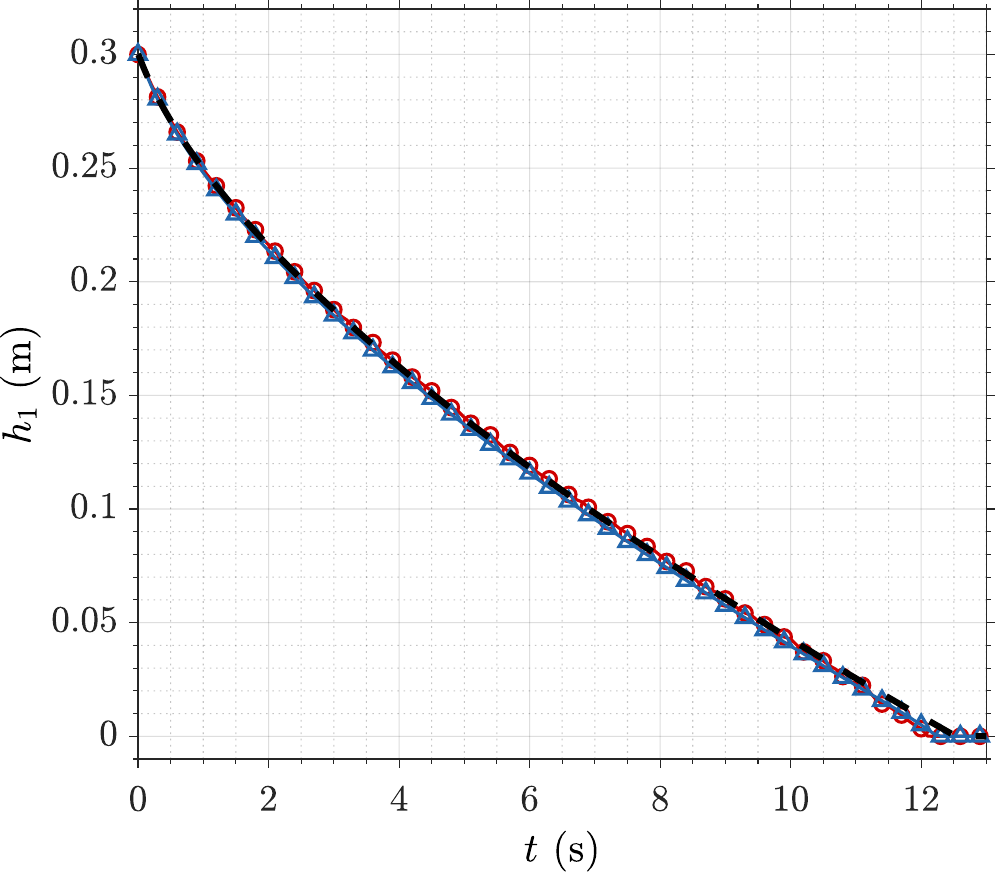}
         \caption*{}
         \label{fig:3.1.6_h1_nonlinear_case_a}
     \end{subfigure}
     \hspace{0.1cm}
     \begin{subfigure}[b]{0.315\textwidth}
         \centering
         \includegraphics[width=\textwidth]{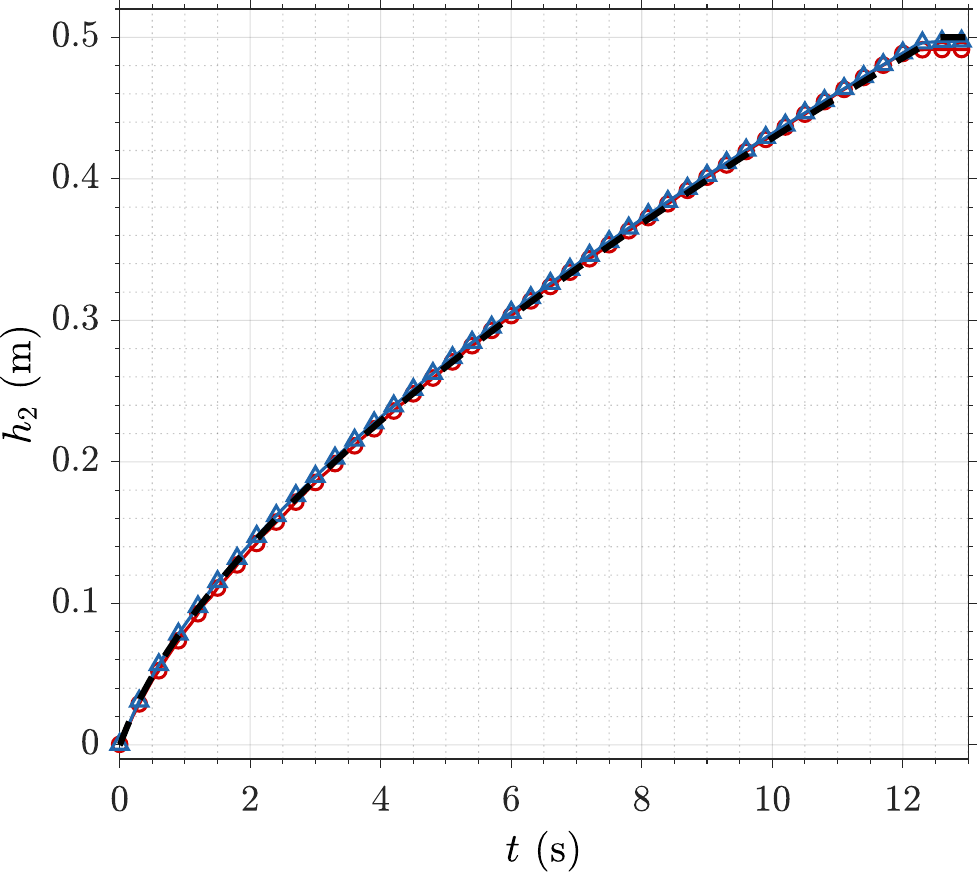}
         \caption{Set 1: high permeability porous media}
         \label{fig:3.1.6_h2_nonlinear_case_a}
     \end{subfigure}
     \hspace{0.1cm}
     \begin{subfigure}[b]{0.32\textwidth}
         \centering
         \includegraphics[width=\textwidth]{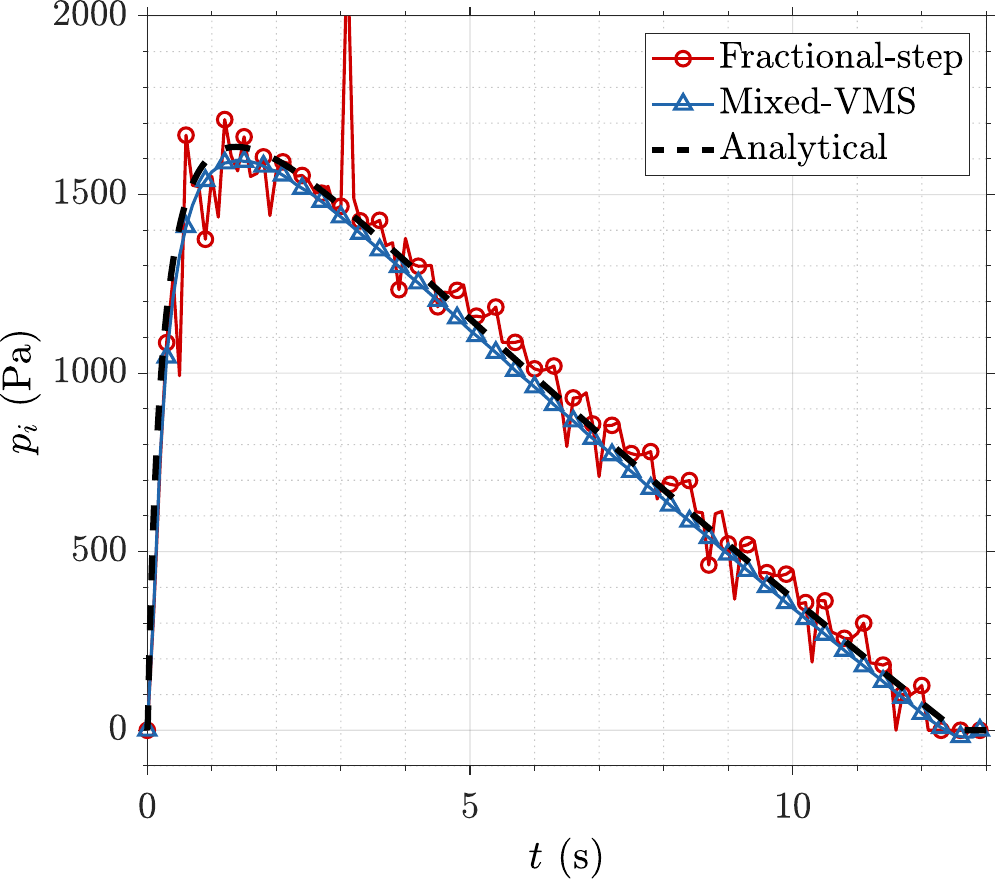}
         \caption*{}
         \label{fig:3.1.6_pi_nonlinear_case_a}
     \end{subfigure}\\
     \vspace{0.3cm}
     \begin{subfigure}[b]{0.32\textwidth}
         \centering
         \includegraphics[width=\textwidth]{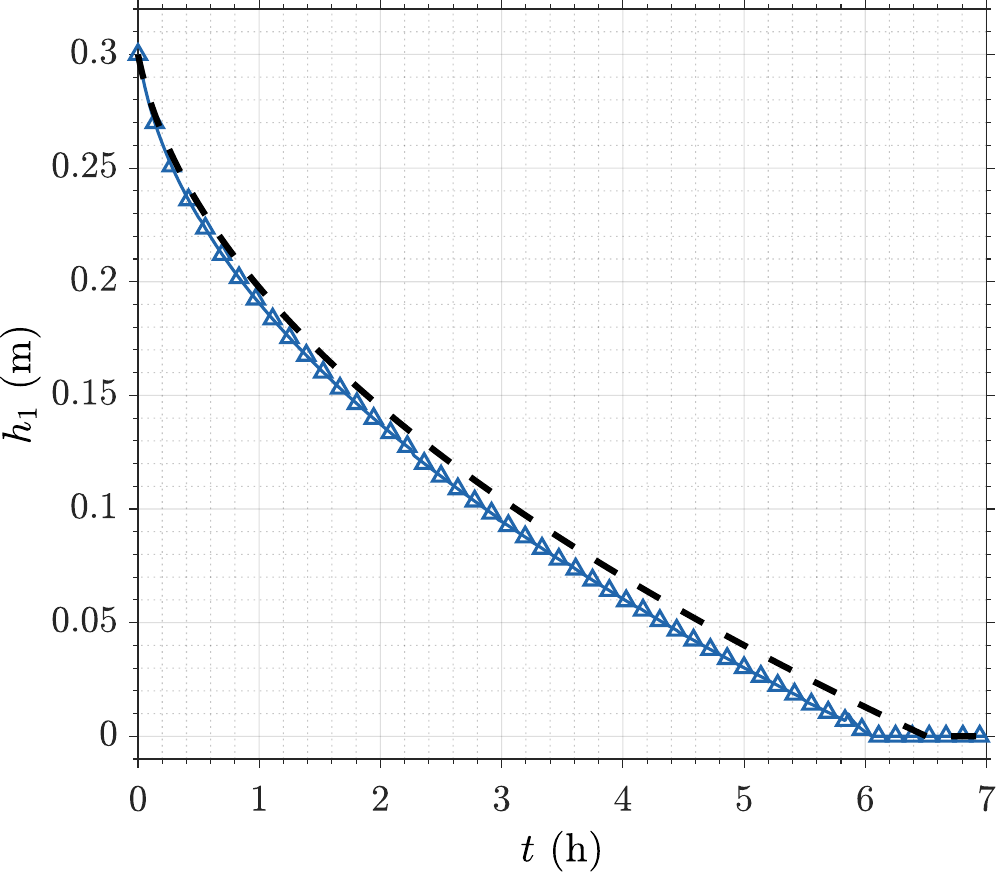}
         \caption*{}
         \label{fig:3.1.6_h1_nonlinear_case_b}
     \end{subfigure}
     \hspace{0.1cm}
     \begin{subfigure}[b]{0.32\textwidth}
         \centering
         \includegraphics[width=\textwidth]{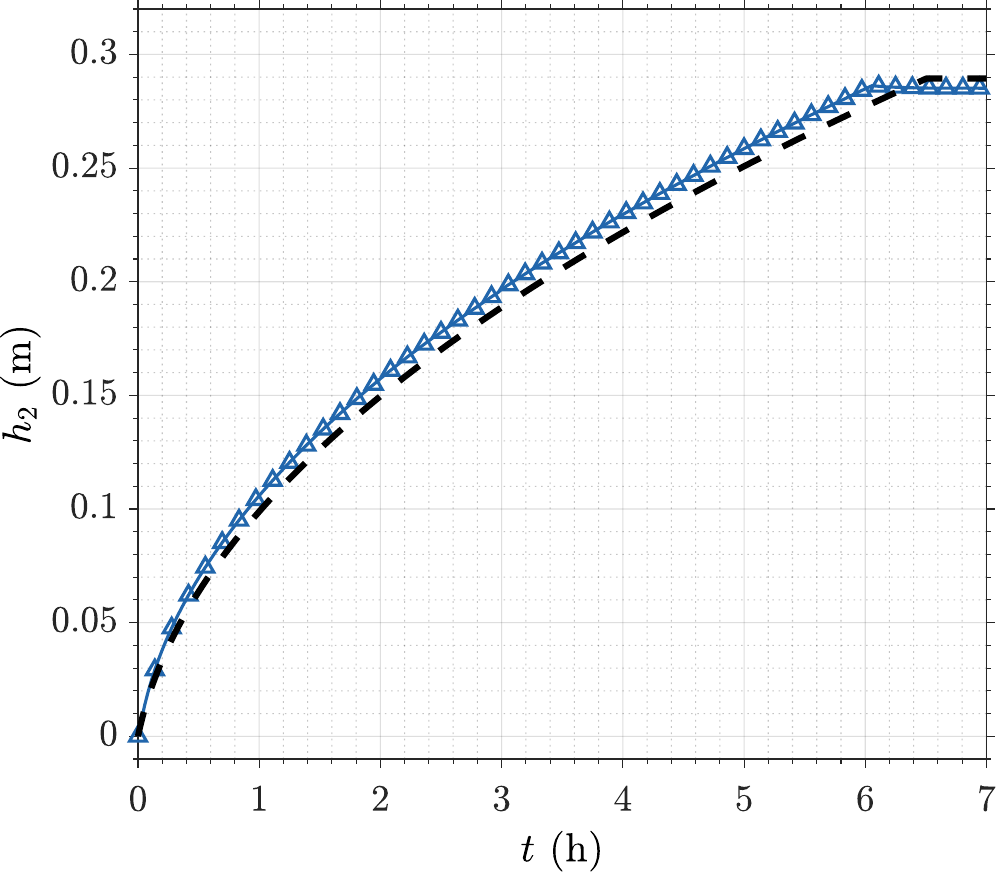}
         \caption{Set 2: low permeability porous media}
         \label{fig:3.1.6_h2_nonlinear_case_b}
     \end{subfigure}
     \hspace{0.1cm}
     \begin{subfigure}[b]{0.32\textwidth}
         \centering
         \includegraphics[width=\textwidth]{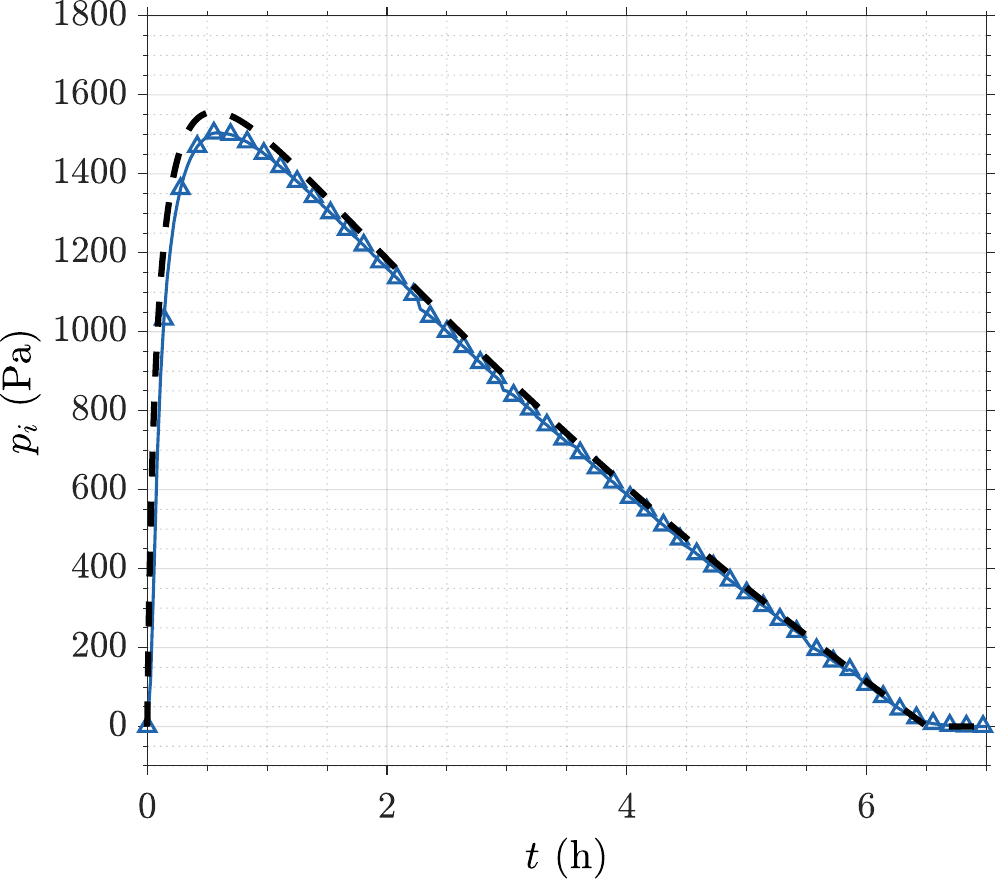}
         \caption*{}
         \label{fig:3.1.6_pi_nonlinear_case_b}
     \end{subfigure}
    \caption{1D verification suite - Case 6: evolution of pore-fluid heights in porous medium 1, $h_{1}$, porous medium 2, $h_{2}$, and interface pressure $p_i$ over time. The fractional-step approach diverges while running Set 2, hence, only the mixed-VMS results are plotted in (b).}
    \label{fig:3.1.6_results}
\end{figure}

The geometric configuration and initial conditions for this problem are illustrated in Fig.~\ref{fig:3.1_model}. Initially, a fluid column is positioned within porous material 1, with an initial height of $h_0$. Driven by gravity, the height of the fluid column inside the porous material 1, denoted as $h_1$, decreases over time, while the height of the pore fluid column in porous material 2, denoted as $h_2$, increases. Adhering to the principle of fluid-volume conservation, the lower limit of $h_1$ and the upper limit of $h_2$ are computed to be 0 m and $(\theta_1/\theta_2)h_0$, respectively. In the current work, we also present the semi-analytical solution for this problem, which is elaborated in \ref{app:1d_case6_derivation}. Additionally, the numerical settings for the MPM simulations, including the mesh size and the PPC, remain consistent with those described in Case 4.

Figs.~\ref{fig:3.1.6_h2_nonlinear_case_a} and \ref{fig:3.1.6_h2_nonlinear_case_b} present the evolution of pore-fluid heights and interface pressure for Sets 1 and 2, respectively. As can be observed from Fig.~\ref{fig:3.1.6_h2_nonlinear_case_a}, the proposed stabilized mixed formulation produces a smoother pressure profile without spurious oscillations, in contrast to the fractional-step approach. Most importantly, the proposed mixed MPM exhibits superior numerical stability, which is shown by the results from Set 2. Here, the solution does not diverge for a large value of $\Delta t$. On the other hand, the fractional-step approach blows up with the same time increment. This behavior is primarily attributed to the strong dependency of the fractional-step approach on the selected time increment, whose critical value is influenced by the magnitude of the hydraulic conductivity or permeability, see \cite{kularathna2021semi, mieremet2016numerical, morikawa2022soil}. A stable $\Delta t_{crit}$ is computed and tested to be about $10^{-5}$ s for this case, implying that about 2.5 billion steps are required to complete the same simulated time. The computation time is extrapolated to be over 6 months on an Intel Core i9-9900X desktop, in contrast to 16 minutes for the mixed MPM. It is worth acknowledging that the solution accuracy for the mixed MPM degrades with increasing $\Delta t$, which is the primary cause of the observed inaccuracies in the profiles depicted in Fig.~\ref{fig:3.1.6_h2_nonlinear_case_b}. However, given the significant improvement in numerical stability and computational efficiency, a slight reduction in accuracy may be considered to be a favorable tradeoff.

\subsection{Couette flow in a composite channel}
\label{subsec:couette_flow_composite}

\begin{figure}[h!]
    \centering  
    \includegraphics[width=0.4\textwidth]{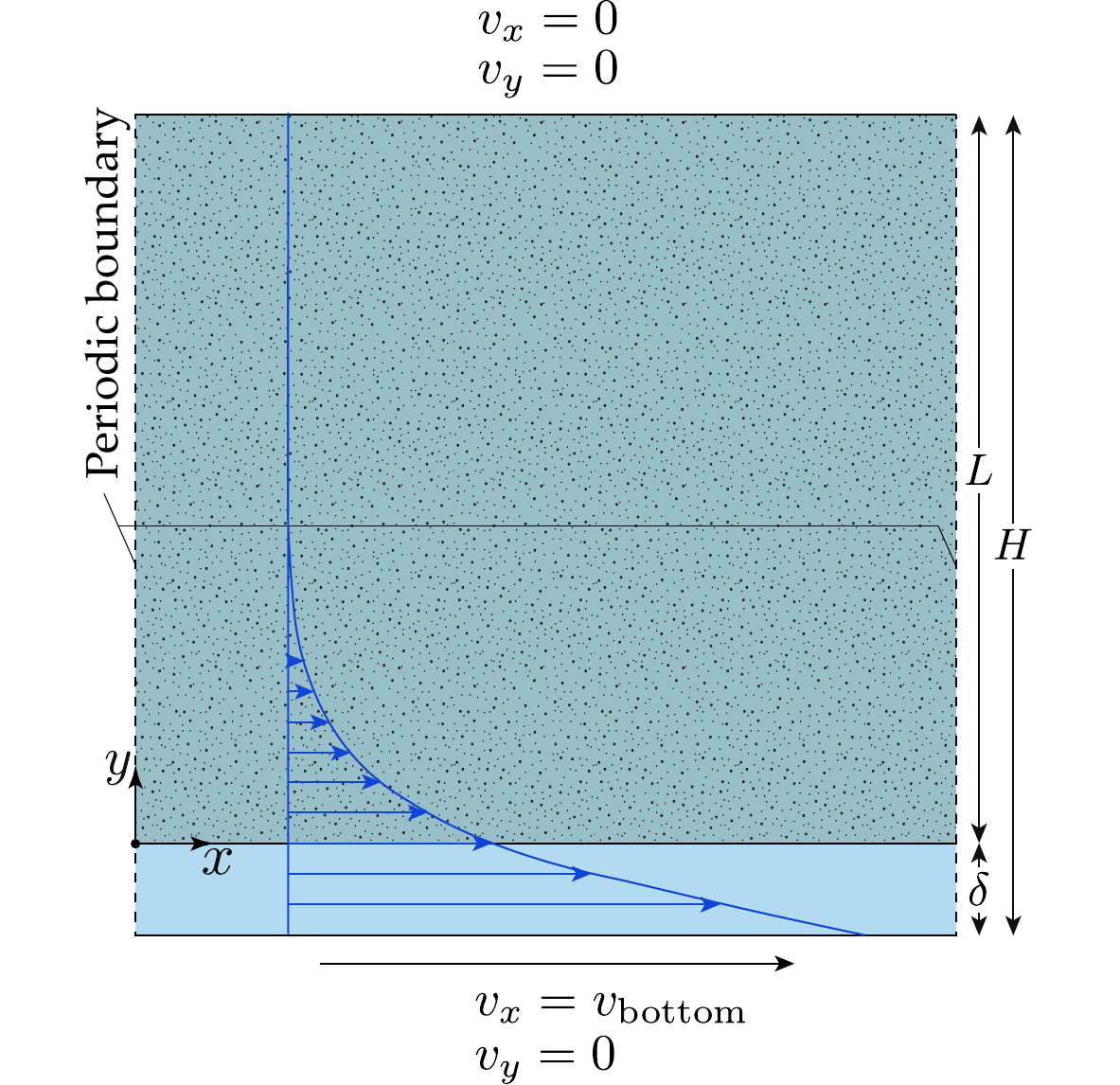}
    \caption{Composite Couette flow: geometrical settings and boundary conditions.}
    \label{fig:3.2_model}
\end{figure}

The 1D verification suite conducted in \Cref{subsec:1d_suite} confirmed the mass conservation and flow dynamics of the proposed method in simulating (quasi-)1D flow under various conditions. However, further investigation is necessary to assess the method's performance in simulating the momentum interaction between the fluid's viscous force and the non-Darcian effects in the porous medium. In the next test, the simulation of Couette flow in a composite channel is performed to analyze fluid flow tangent to a porous interface without the presence of head difference, which might also drive seepage flow. In this setup, two parallel plates are positioned at a distance $H$, with a porous medium of length $L<H$ attached to the top plate, creating a gap with a length of $\delta$ (see Fig.~\ref{fig:3.2_model}). The space between the plates is then filled with water, which saturates the porous medium. In this configuration, the bottom plate interfacing with the clear fluid is moved with velocity $v_{\text{bottom}}$, while the top plate remains fixed. Both of the plates are assumed to be non-slip.

In this analysis, the 2D simulation domain is set with $\delta/H=0.1$ and $L/H=0.9$, where the distance between the plates is considered to be $H=0.001$ m. The porous medium is assumed to have a porosity of $\theta=0.5$ and permeability in two cases: $k=2.5 \times 10^{-9}$ and $2.5 \times 10^{-11}\,\rm{m}^2$, corresponding to Darcy numbers, $Da_H=\Tilde{k}/H^2$, of $10^{-2}$ and $10^{-4}$, respectively. Here, $\Tilde{k}=k/\theta^2$ is the intrinsic permeability normalized by porosity-squared, which is introduced for mathematical convenience. The Forchheimer drag force term is considered, assuming Ergun's constants ($A=150$ and $B=1.75$). Additionally, the effective viscosity, $\mu_e$, associated with the Brinkman term is considered in two cases: (i) $\mu_e=\mu$ and (ii) $\mu_e=4\mu$. The study also focuses on two cases of Reynolds number: $Re=\rho v_{\text{bottom}} \delta/\mu=1$ and 100; thus, the horizontal velocity of the bottom plate is set as $v_{\text{bottom}}=0.01$ and $1$ m/s, respectively. The steady-state analytical solution for this problem is available as elaborated in \ref{app:couette}.

The MPM model considers a structured quadrilateral background grid of size $h=H/80$, each with $4\times4$ PPC. To obtain a fully developed steady-state flow condition, the periodic boundaries are used on both the left and right ends of the simulation domain. The periodic BC also allows us to save computational costs by reducing the number of cells horizontally. In the current work, $8$ elements are considered horizontally for the computational domain, totaling $80\times8$ cells with 10,240 material points. The time increment is set as $\Delta t = 5\times10^{-5}$ s and $1\times10^{-5}$ s for $Re=1$ and 100, respectively, where the simulations are performed until they reach the quasi-steady-state solution. To apply the velocity BC at the bottom plate, we enforce a nonhomogeneous displacement BC, $\overline{u}_{I,x}$, to achieve the desired velocity field, $v_\text{bottom}$, as: 
\begin{eqnarray}
    \overline{u}_{I,x} = \frac{\beta_N \Delta t}{\gamma_N}v_\text{bottom} + \left(1-\frac{\beta_N}{\gamma_N}\right)\Delta t {v}^n_{I,x} + \left(\frac{1}{2}-\frac{\beta_N}{\gamma_N}\right) \Delta t^2 {a}^n_{I,x}\,,
\end{eqnarray}
following the Newmark-$\beta$ time integration scheme. In addition, to model the jump in viscosity between the non-porous and porous domains, the porosity-dependent viscosity function, $\widetilde{\mu}$, as defined in Eq.~\eqref{eq:porosity_dependent_visc}, is employed with a precomputed slope of $\eta = 6$ for case (ii) with $\mu_e = 4\mu$. For this verification test, the performance of the FLIP, APIC, and TPIC transfer schemes are investigated.

\begin{figure}[h!]
    \centering
     \begin{subfigure}[b]{0.4\textwidth}
         \centering
         \includegraphics[width=\textwidth]{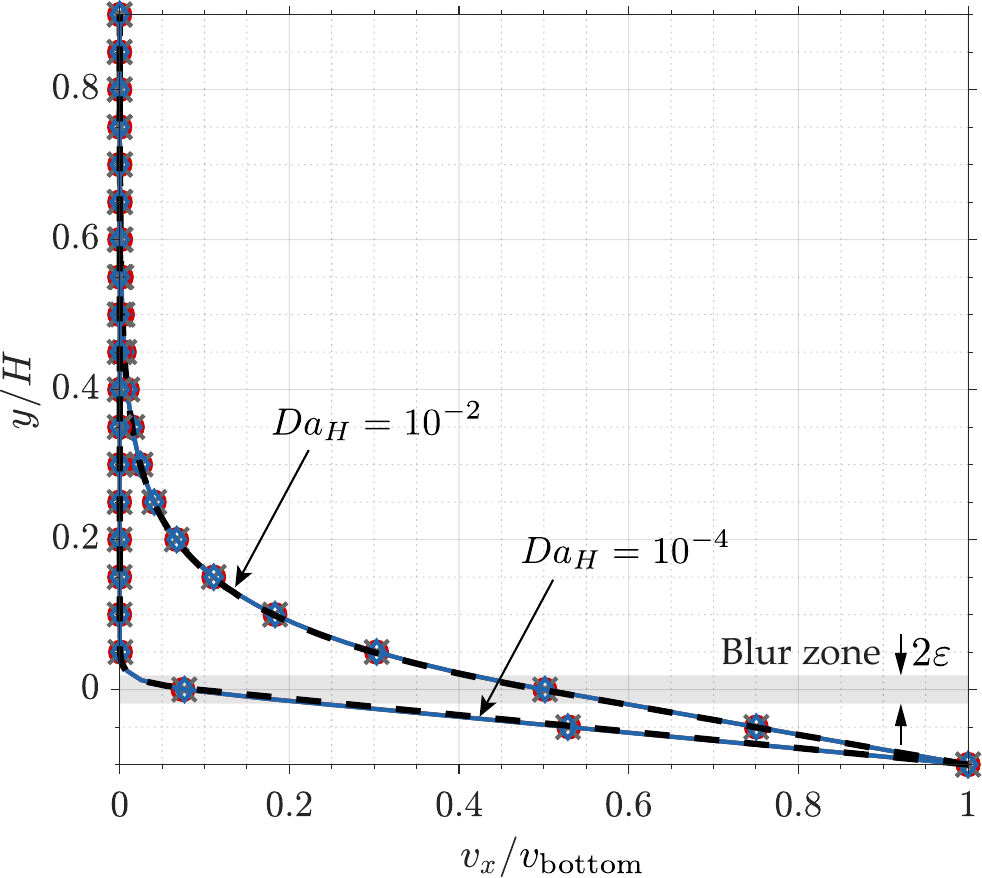}
         \caption{$Re=1,\,\mu_e=\mu$.}
         \label{fig:3.2_results_re_1_gamma_1}
     \end{subfigure}
     \hspace{0.1cm}
     \begin{subfigure}[b]{0.4\textwidth}
         \centering
         \includegraphics[width=\textwidth]{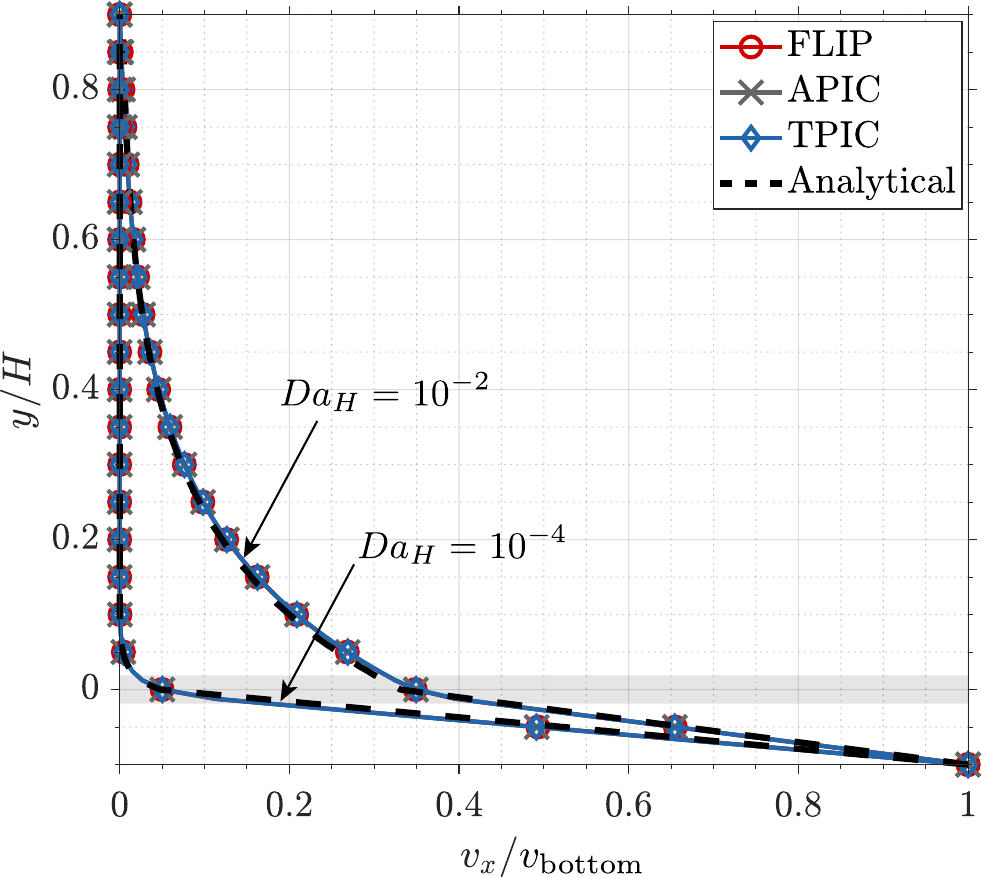}
         \caption{$Re=1,\,\mu_e=4\mu$.}
         \label{fig:3.2_results_re_1_gamma_2}
     \end{subfigure}
     \hfill\\
     \vspace{0.2cm}
      \begin{subfigure}[b]{0.4\textwidth}
         \centering
         \includegraphics[width=\textwidth]{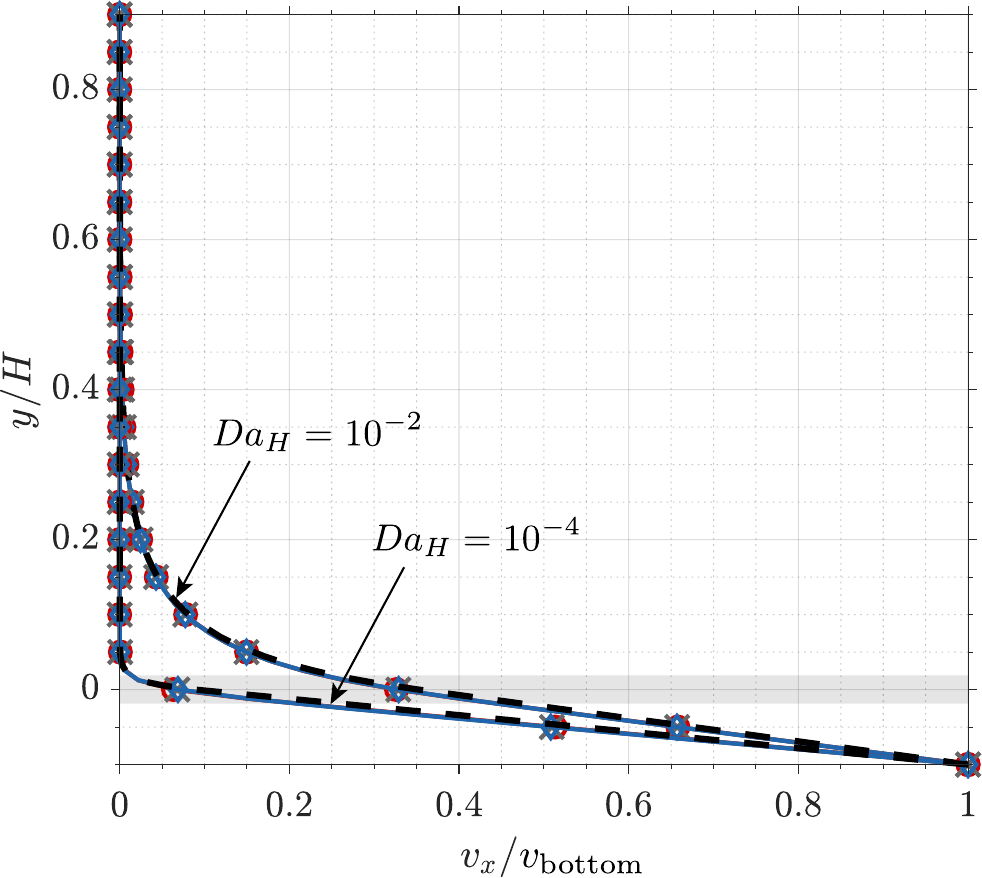}
         \caption{$Re=100,\,\mu_e=\mu$.}
         \label{fig:3.2_results_re_100_gamma_1}
     \end{subfigure}
     \hspace{0.1cm}
     \begin{subfigure}[b]{0.4\textwidth}
         \centering
         \includegraphics[width=\textwidth]{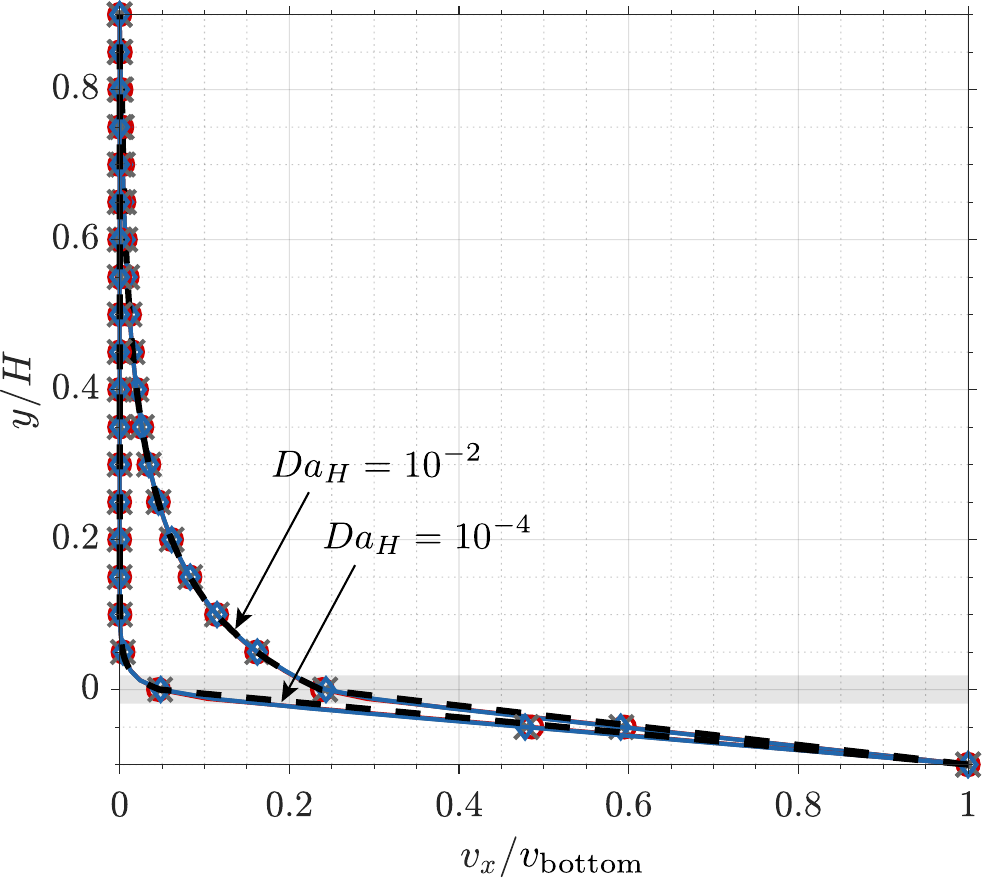}
         \caption{$Re=100,\,\mu_e=4\mu$.}
         \label{fig:3.2_results_re_100_gamma_2}
     \end{subfigure}
    \caption{Composite Couette flow: normalized horizontal velocity for different simulation conditions. The blur zone around the porous interface is highlighted in gray.}
    \label{fig:3.2_results}
\end{figure}

The obtained normalized horizontal velocity is presented in Fig.~\ref{fig:3.2_results} for different simulation settings. As can be seen, all the mixed MPM results performed with FLIP, APIC, and TPIC schemes agree very well with the steady-state analytical solutions. The flows with a larger Darcy number, $Da_H=10^{-2}$, exhibit a taller momentum boundary layer inside the porous medium as the fluid can flow relatively faster due to larger permeability. Additionally, the influence of effective viscosity $\mu_e$ on the velocity profile is observed. When $\mu_e=\mu$, the velocity profile smoothly decreases over the increasing height, whereas when $\mu_e=4\mu$, a kink in the velocity profile can be observed around the porous interface. This kink is caused by the difference in velocity gradient to satisfy the tangential stress equilibrium at the interface \cite{vafai1987analysis}. Moreover, the results of higher Reynolds numbers also show a very good match with the analytical solutions. Here, the cases with higher Reynolds numbers display a more isolated boundary layer compared to those with lower Reynolds numbers, where the velocity change is relatively more gradual with increasing height. All these flow features can be accurately captured by the three transfer schemes.

\begin{figure}[h!]
    \centering
     \begin{subfigure}[b]{0.408\textwidth}
         \centering
         \includegraphics[width=\textwidth]{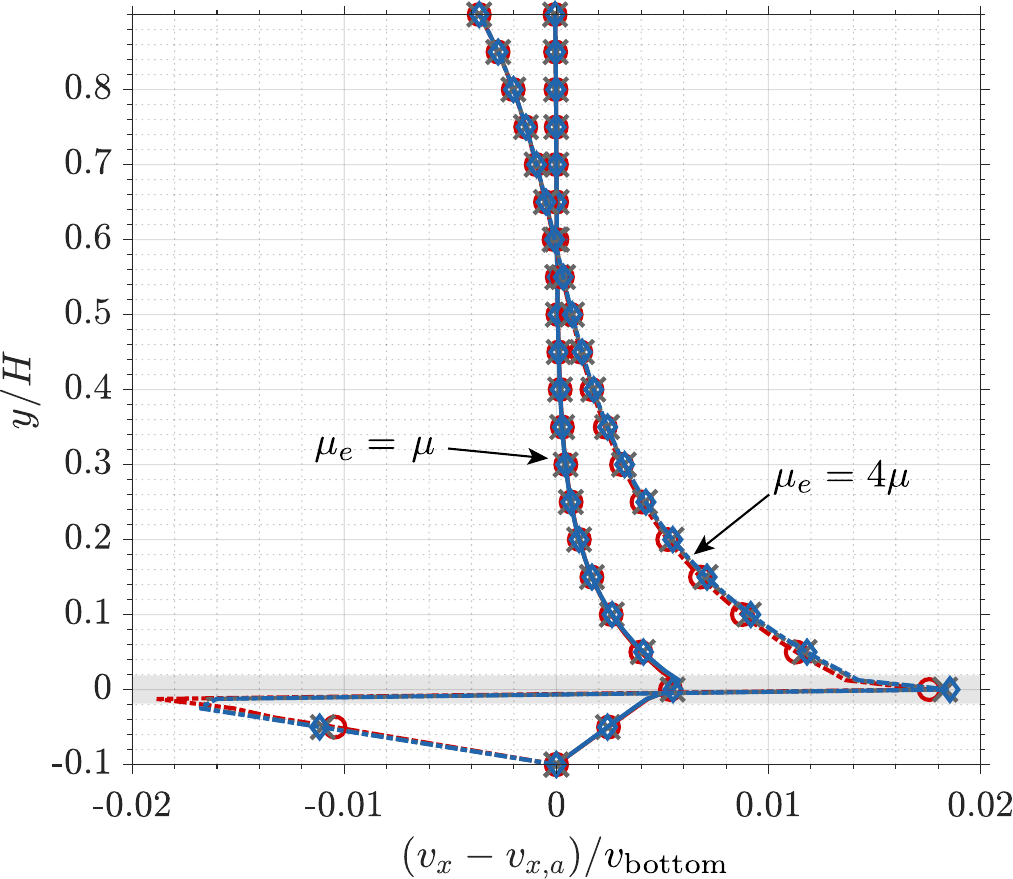}
         \caption{}
         \label{fig:3.2_error_profile_vx}
     \end{subfigure}
     \hspace{0.2cm}
     \begin{subfigure}[b]{0.4\textwidth}
         \centering
         \includegraphics[width=\textwidth]{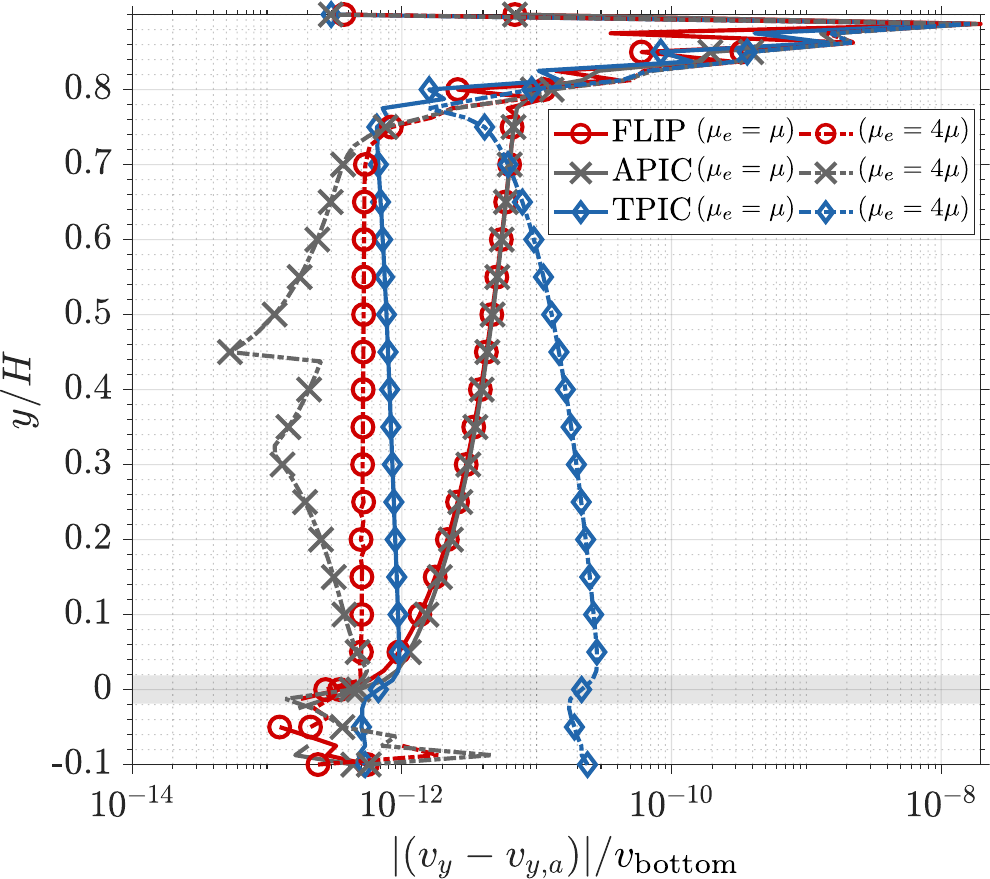}
         \caption{}
         \label{fig:3.2_error_profile_vy}
     \end{subfigure}
    \caption{Composite Couette flow: error profile across the channel height for (a) horizontal velocity and (b) vertical velocity considering $h=H/80$. The vertical velocity is plotted in a semilog plot to emphasize the order of error magnitude.}
    \label{fig:3.2_error_profile}
\end{figure}

Next, the error distribution and the convergence profile are investigated at the steady state (assumed to be at $t=0.2$ s). Here, the configuration with $Re=1$ and $Da_H=10^{-2}$ is selected for comparison, where the two cases of effective viscosity are considered: (i) $\mu_e = \mu$ and (ii) $\mu_e = 4\mu$. Fig.~\ref{fig:3.2_error_profile} depicts the error distribution of horizontal and vertical velocities across the channel's height, where the vertical velocity is expected to be zero everywhere, i.e.~$v_{y,a}=0$. Here, the error profiles of the horizontal velocity are dominated by the discretization errors near the blurred interface. Case (ii) is observed to have a larger magnitude with overestimated values around ${y=0^+}$ and underestimated values around ${y=0^-}$. Furthermore, the error near the top plate also tends to increase. This is caused by a slightly different assumption of BC considered by the analytical solution (see more on \ref{app:couette}). Meanwhile, different trends are observed for case (i). Since the momentum boundary layer height is lower in this case, the error near the top plate is much smaller. We also can observe that the continuity in viscosity reduces the error around the blurred interface. Overall, the errors for both cases are relatively small, capped below $0.6\%$ and $2\%$ for cases (i) and (ii), respectively. In contrast, the error profiles of the vertical velocity are significantly smaller, with the maximum magnitude being about $10^{-8}$ near the top plate for both cases. 

\begin{figure}[h!]
    \centering
     \begin{subfigure}[b]{0.40\textwidth}
         \centering
         \includegraphics[width=\textwidth]{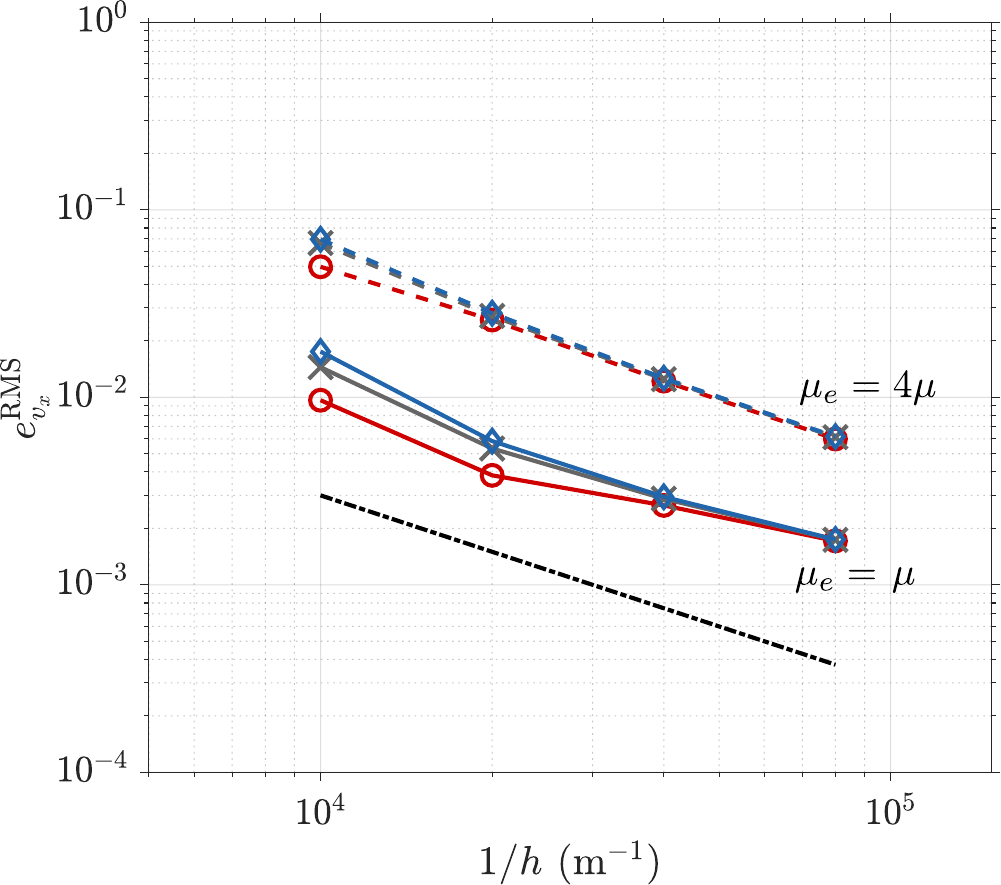}
         \caption{}
         \label{fig:3.2_error_conv_vx_rms}
     \end{subfigure}
     \hspace{0.2cm}
     \begin{subfigure}[b]{0.4\textwidth}
         \centering
         \includegraphics[width=\textwidth]{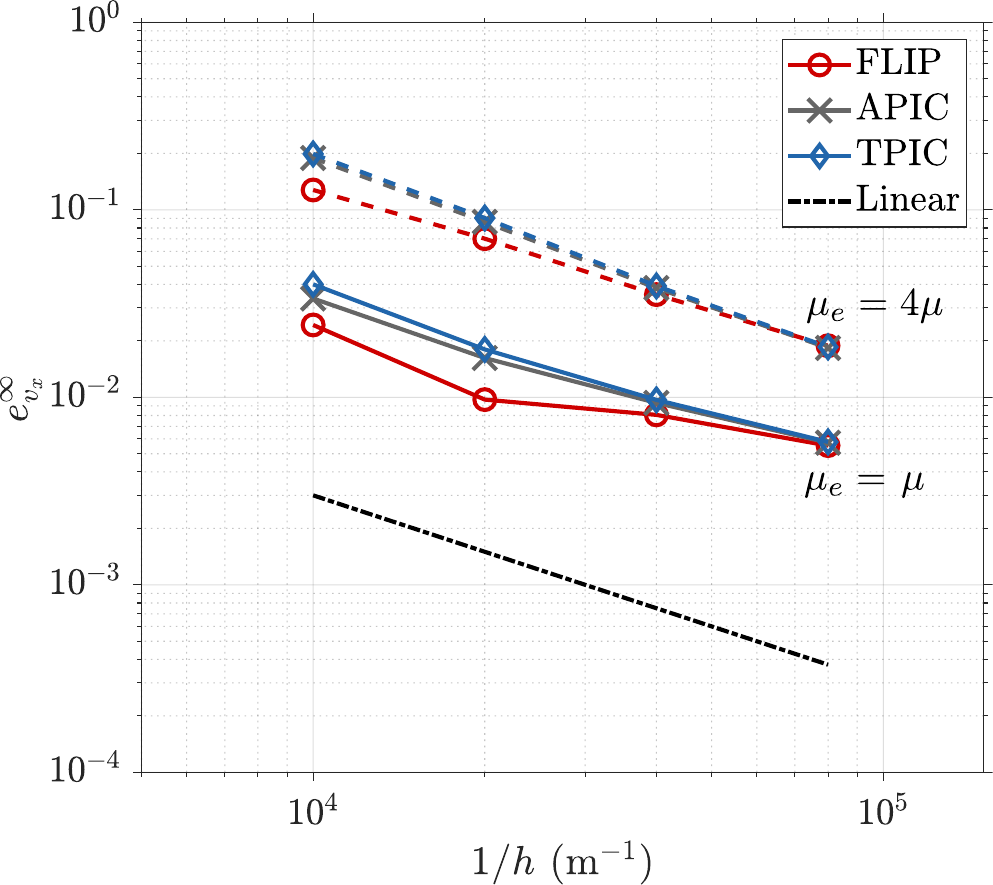}
         \caption{}
         \label{fig:3.2_error_conv_vx_infty}
     \end{subfigure}
    \caption{Composite Couette flow: convergence plot of horizontal velocity along with mesh refinement measured at $t=0.2$ s: (a) RMS error and (b) infinity error.}
    \label{fig:3.2_error_conv}
\end{figure}

Furthermore, convergence study of the horizontal velocity $v_x$ along with mesh refinement is performed with four different cell sizes, $h=H/10\sim H/80$. During the mesh refinement process, the PPC is maintained at $4\times 4$, and the ratio between mesh size and time step is kept constant at $h/\Delta t = 0.25$. The numerical error is quantified using the root-mean-square (RMS) and the infinity errors, which are defined as follows:
\begin{eqnarray}
    e^{\rm RMS}_v = \sqrt{ \frac{1}{N_n} \sum_{I=1}^{N_n} \left(\frac{v_I - v_a(\tb x_I)}{v_\text{bottom}} \right)^2 } \,, \qquad
    e^{\infty}_v = \max_I\left( \left\vert \frac{v_I - {v}_a(\tb x_I)}{v_\text{bottom}}\right\vert\right) \,,
\end{eqnarray}
where $N_n$ is the total compute node in the simulation domain. As depicted in Fig.~\ref{fig:3.2_error_conv}, the convergence trends of the two error measures are approximately linear for all schemes. As discussed earlier in \Cref{subsubsec:1d_case5}, a convergence rate lower than quadratic is expected due to the presence of a blurred porous interface. Further research on enhancing accuracy near the blurred interface, such as through the implementation of adaptive mesh refinement, is still relatively unexplored in MPM and presents an interesting avenue for future work.

\subsection{Dam break flows through rigid porous media}
\label{subsec:2d_dam_break}

\begin{figure}[h!]
    \centering
    \includegraphics[width=0.8\textwidth]{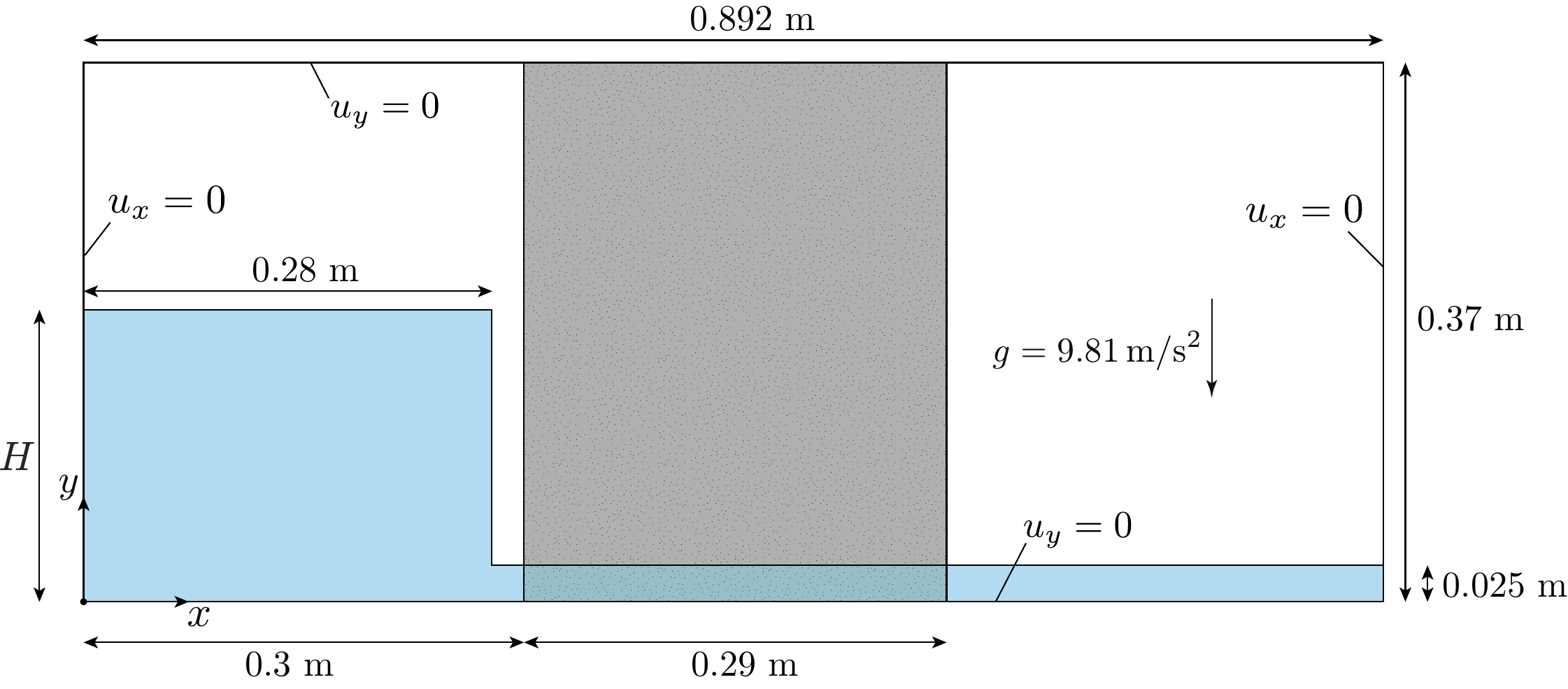}
    \caption{Dam break flows through rigid porous media: geometrical settings and boundary conditions.}
    \label{fig:3.3_model}
\end{figure}

An experimental validation of a water dam-break flow through a rigid porous column is conducted next. This validation test was first proposed by \citet{liu1999numerical} to investigate flow in porous media under two different flow regimes characterized by different ranges of Reynolds numbers \cite{del2012three}. As illustrated in Fig.~\ref{fig:3.3_model}, the experimental setup was done in a tank with a length of 0.892 m and height of 0.58 m. A porous column dam with 0.29 m length and 0.37 m height was placed at the approximate center of the tank, with a distance of 0.3 m from the left wall. In the experiment, a gate was installed at 0.02 m away from the left side of the porous column, where a water reservoir with a certain height $H$ was initially retained. Here, the effective viscosity inside the porous medium is set to be the same as fluid viscosity, i.e.~$\mu_e = \mu$ and $\widetilde{\mu}=\mu$. In the simulation, at $t=0$ s, the gate is instantaneously released and the fluid is allowed to flow under the gravity of $g=9.81$ m/s$^2$. Since the fluid mainly flows horizontally through the porous column, there is no need to model the entire height of the tank, hence, the height of the simulation domain can be reduced following the height of the porous column. Furthermore, the wall boundaries are set to be free-slip as indicated in Fig.~\ref{fig:3.3_model}.

Two types of porous media are considered in this validation test, i.e.~glass beads with a porosity of 0.39 and crushed rocks with a porosity of 0.49. Both are treated as cohered in place so they do not flow. The average grain diameter $d$ of these materials is 3 mm and 15.9 mm, respectively. To study the flow conditions inside the two different porous materials, the initial water height $H$ is set differently according to the experiments conducted by \citet{liu1999numerical}. This is correspondingly denoted as Case 1 and 2 in this study, where the height $H$ along with other material parameters are listed in \Cref{tab:3.3_cases}. Here, the same drag force constants $A$ and $B$ are utilized as suggested by \citet{liu1999numerical}, even though they calibrated their model with the drag-force model suggested by \citet{van1995porous}, which includes an extra inertial resistance term. Here, since the dam-break flow is rather non-oscillatory, we also neglected the effect of the Keulegan-Carpenter number appearing in the nonlinear drag force term suggested by \citet{van1995porous}.

\begin{table}[h!]
\centering
\caption{Dam break flows through rigid porous media: parameters of porous media for different cases.}
\label{tab:3.3_cases}
\begin{tabular}{||c|c|c|c|c|c|c|c|c||}
\hline
      Case & $H$ (m) & Porous material & $d$ (m) & $\theta$ (-) & $Re$ (-)\cite{del2012three} & $k$ (m$^2$) & $A$ (-) & B (-) \\ \hline \hline
1 & 0.14 & Glass beads & 0.003             & 0.39               & 9.6              & 7.17 $\times 10^{-9}$  & 200 & 1.1     \\ \hline
2 & 0.25 & Crushed rocks & 0.0159             & 0.49               & 325              & 1.14 $\times 10^{-7}$ & 1000 & 1.1      \\ \hline
\end{tabular}
\end{table}

The MPM simulations are performed using quadrilateral background grids with a cell size of $h=0.01$ m, resulting in a total of 89$\times$37 cells. The fluid domain is discretized as 8,712 and 13,640 material points for Cases 1 and 2, respectively, which are initially arranged in $4\times4$ PPC configuration. In this study, we set the initial material point's volume uniformly, resulting in smaller initial masses for material points within the porous domain. Additionally, the time increment is set to be $\Delta t = 0.001$ s. The simulations primarily utilize the FLIP transfer scheme, with the APIC and TPIC schemes simulated for comparison purposes.

\begin{figure}[h!]
    \centering
    \includegraphics[width=1.0\textwidth]{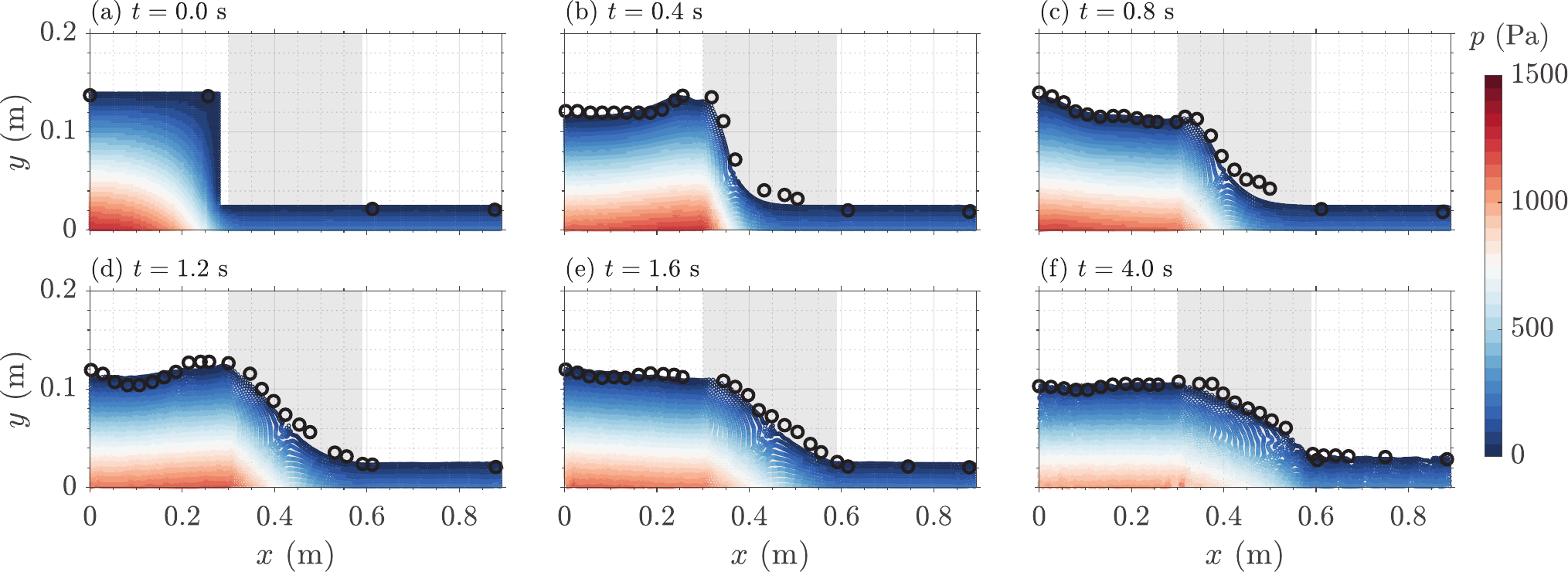}
    \caption{Dam break flows through rigid porous media - Case 1: evolution of pressure contour obtained by the FLIP scheme. The experimental data conducted by \citet{liu1999numerical} are denoted by the black circles.}
    \label{fig:3.3_case1_results}
\end{figure}

\begin{figure}[h!]
    \centering
    \includegraphics[width=1.0\textwidth]{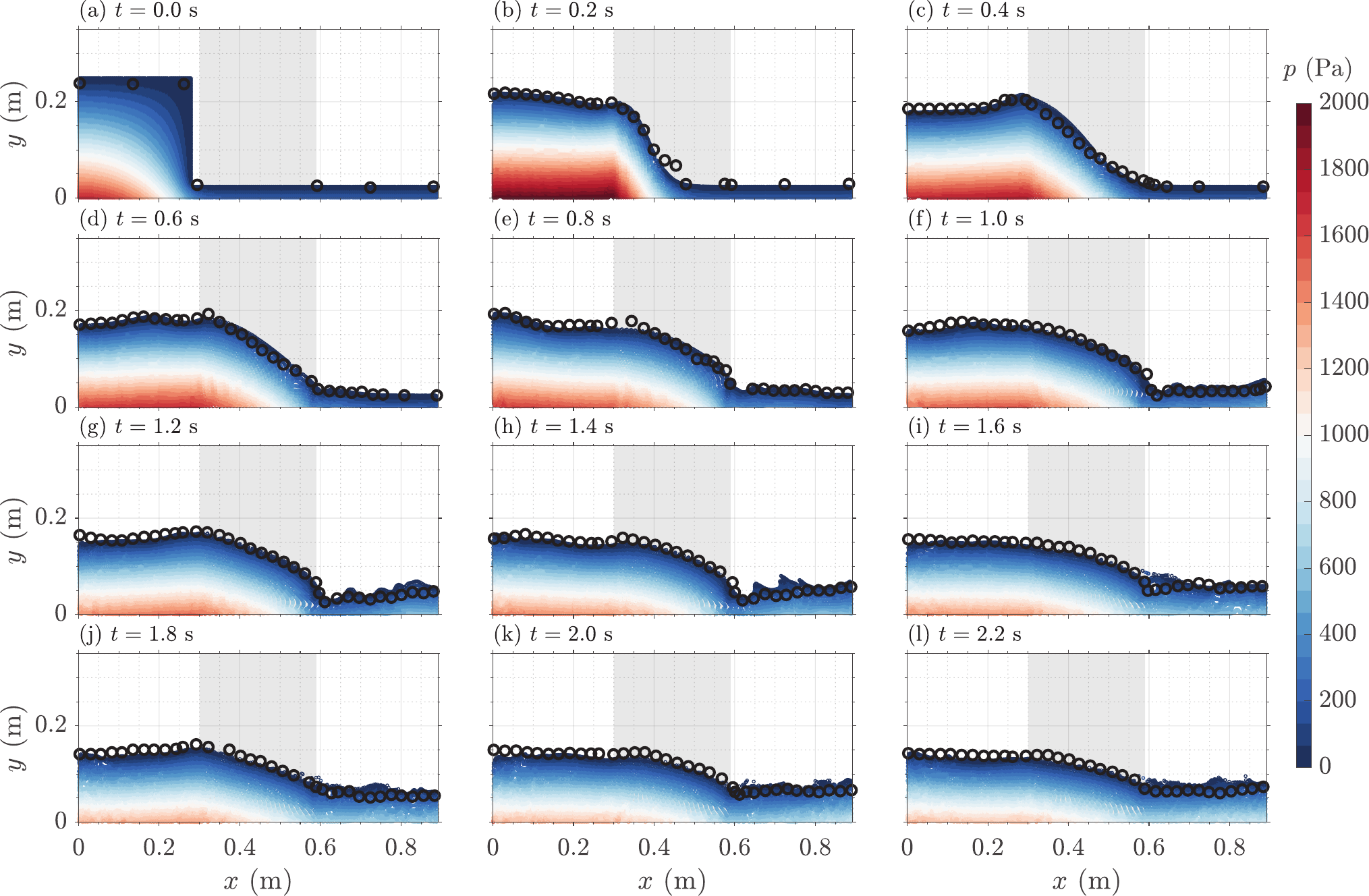}
    \caption{Dam break flows through rigid porous media - Case 2: evolution of pressure contour obtained by the FLIP scheme. The experimental data conducted by \citet{liu1999numerical} are denoted by the black circles.}
    \label{fig:3.3_case2_results}
\end{figure}

\begin{figure}[h!]
    \centering
     \begin{subfigure}[b]{0.323\textwidth}
         \centering
         \includegraphics[width=\textwidth]{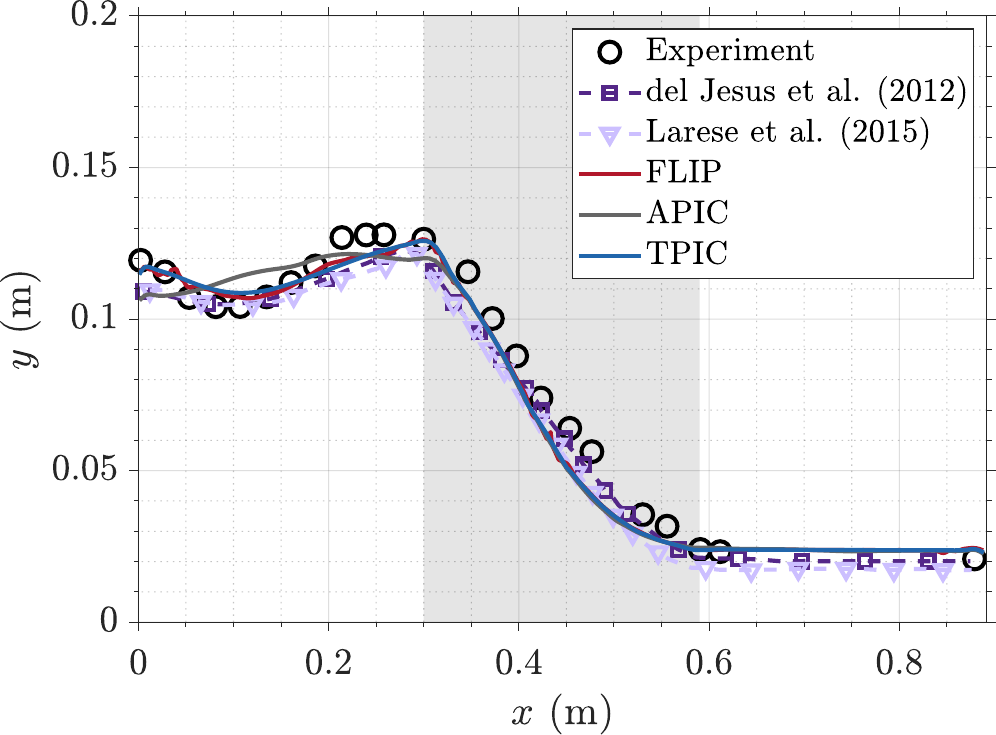}
         \caption{}
         \label{fig:3.3_case1_results_boundaries_compare}
     \end{subfigure}
     \hspace{0.05cm}
     \begin{subfigure}[b]{0.323\textwidth}
         \centering
         \includegraphics[width=\textwidth]{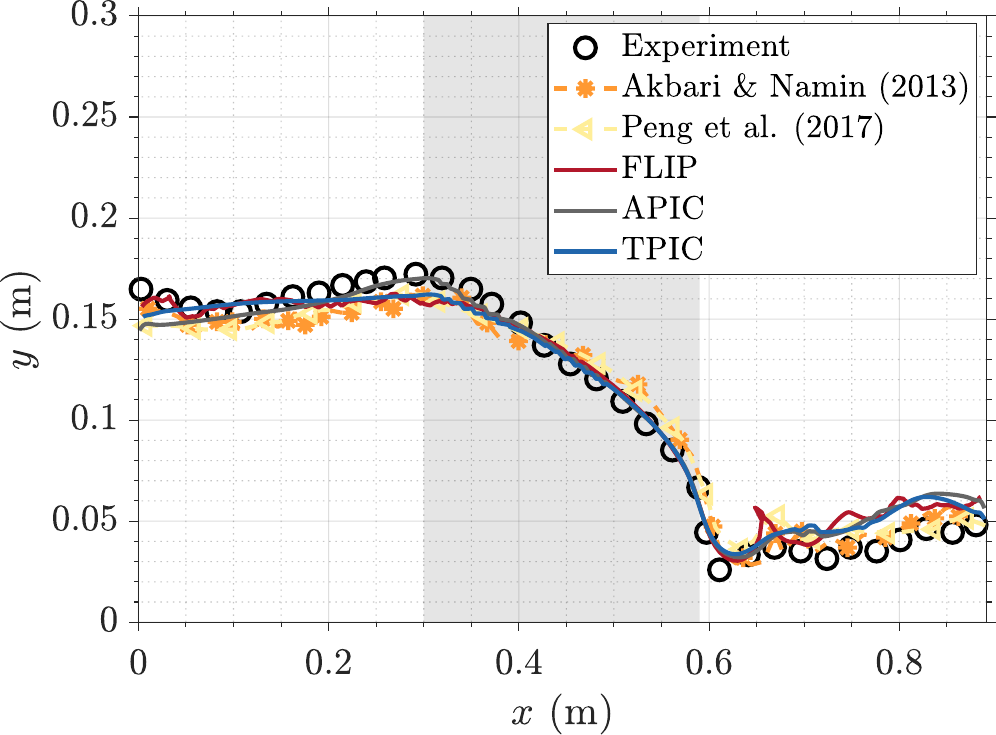}
         \caption{}
         \label{fig:3.3_case2_results_boundaries_compare}
     \end{subfigure}
     \hspace{0.05cm}
     \begin{subfigure}[b]{0.323\textwidth}
         \centering
         \includegraphics[width=\textwidth]{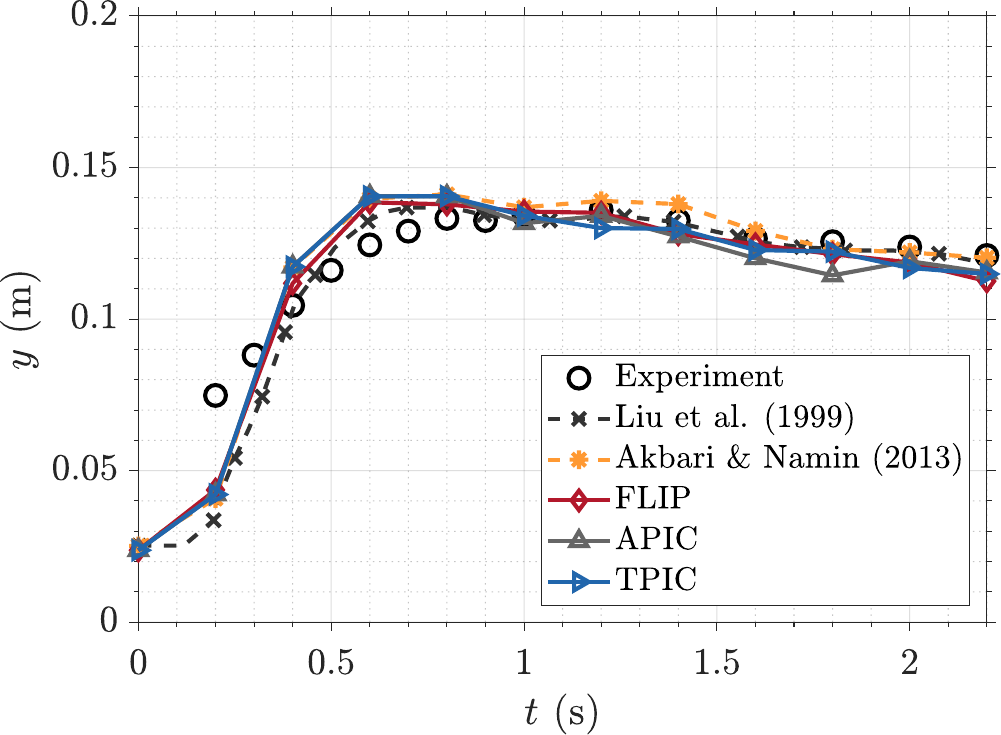}
         \caption{}
         \label{fig:3.3_case2_results_midporeheight_compare}
     \end{subfigure}
    \caption{Dam break flows through rigid porous media: free-surface profile at $t=1.2$ s for (a) Case 1 and (b) Case 2, and (c) time evolution of water height measured at the middle of the porous media, $x=0.445$ m, for Case 2. The obtained results of mixed MPM considering FLIP, APIC, and TPIC transfer schemes are compared with the experimental data conducted by \citet{liu1999numerical} and numerical results by \citet{liu1999numerical} (VOF-FDM), \citet{del2012three} (VOF-FVM), \citet{akbari2013moving} (SPH), \citet{larese2015finite} (edge-based FEM), and \citet{peng2017multiphase} (SPH). }
    \label{fig:3.3_spatial_temporal_plots}
\end{figure}

The snapshots of the pressure field and free-surface topology are depicted in Figs.~\ref{fig:3.3_case1_results} and \ref{fig:3.3_case2_results} for Cases 1 and 2, respectively. Here, we exclusively present results obtained through the FLIP scheme, as the overall behavior from the APIC and TPIC schemes shows only marginal differences. Both figures include a comparison with experimental data by \citet{liu1999numerical}, represented by black circles. Overall, there is good agreement between numerical and experimental measurements of the water height. The pressure profile is also smooth in all simulations. This was challenging to achieve without a smoothing algorithm in the weakly compressible MPM due to volumetric locking (e.g.~\cite{martinelli2016soil}). However, as evident in the figures, during the early stages of the simulations ($t\leq0.8$ s for Case 1 and $t\leq0.4$ s for Case 2), there is an apparent difference in water height between the two measurements. According to \citet{liu1999numerical}, these discrepancies are primarily caused by variations in the gate opening mechanisms considered. The numerical simulations assume instantaneous gate opening, whereas in the experiment, the gate was manually removed over about 0.1 s. The finite-duration removal in the experiment results in a relatively slower pressure release in the water column, accelerating the advancement of the water surface near the bottom inside of the porous dam compared to the upper surface. As the water rushes toward the porous medium, it generates an upward jet upon impacting the surface of the porous medium. The faster movement of water in the experiment leads to a higher water pile-up upon entering the porous medium compared to the simulation results. Moreover, since the simulations are performed in 2D, our simulation is unable to resolve the irregular 3D flow structures and vorticity induced by the gate geometry at the waterfront, e.g.~\cite{asai2021fluid}. In addition, uncertainties in empirical drag parameters ($A$ and $B$) in the numerical model may contribute to disagreements. This observation aligns with findings in many earlier studies (e.g., \cite{del2012three, akbari2013moving, larese2015finite, akbari2014modified, higuera2014three}), leading to calibration efforts for these parameters and the consideration of other drag force models \cite{losada2016modeling}. As the water movement stabilizes over time, the accuracy of the numerical model improves, resulting in a well-predicted free-surface configuration in the later stages of the simulation.


The obtained free-surface profiles at $t=1.2$ s are compared with several existing works \cite{del2012three, akbari2013moving, larese2015finite, peng2017multiphase} in Figs.~\ref{fig:3.3_case1_results_boundaries_compare} and \ref{fig:3.3_case2_results_boundaries_compare} for Cases 1 and 2, respectively. Additionally, the evolution of water height at the center of the porous domain ($x=0.445$ m) for Case 2 is compared with experimental recordings and prior numerical works \cite{liu1999numerical, akbari2013moving} in Fig.~\ref{fig:3.3_case2_results_midporeheight_compare}. In general, good agreements between the numerical results and the experiment are observed. However, due to the reasons explained earlier, the water height obtained by the proposed method in Fig.~\ref{fig:3.3_case2_results_midporeheight_compare} aligns better with past numerical works than the experimental data, especially at $t<0.4$ s, where they underpredict, and between $0.4$ s $< t< 1.0$ s, where they overpredict the water height. The measured height after $t>1.0$ s behaves very much in accordance with the experiments for all simulations.

\subsection{3D dam break of mixed free-surface--porous flow}
\label{subsec:3d_dam_break}

\begin{figure}[h!]
    \centering
     \includegraphics[width=0.7\textwidth]{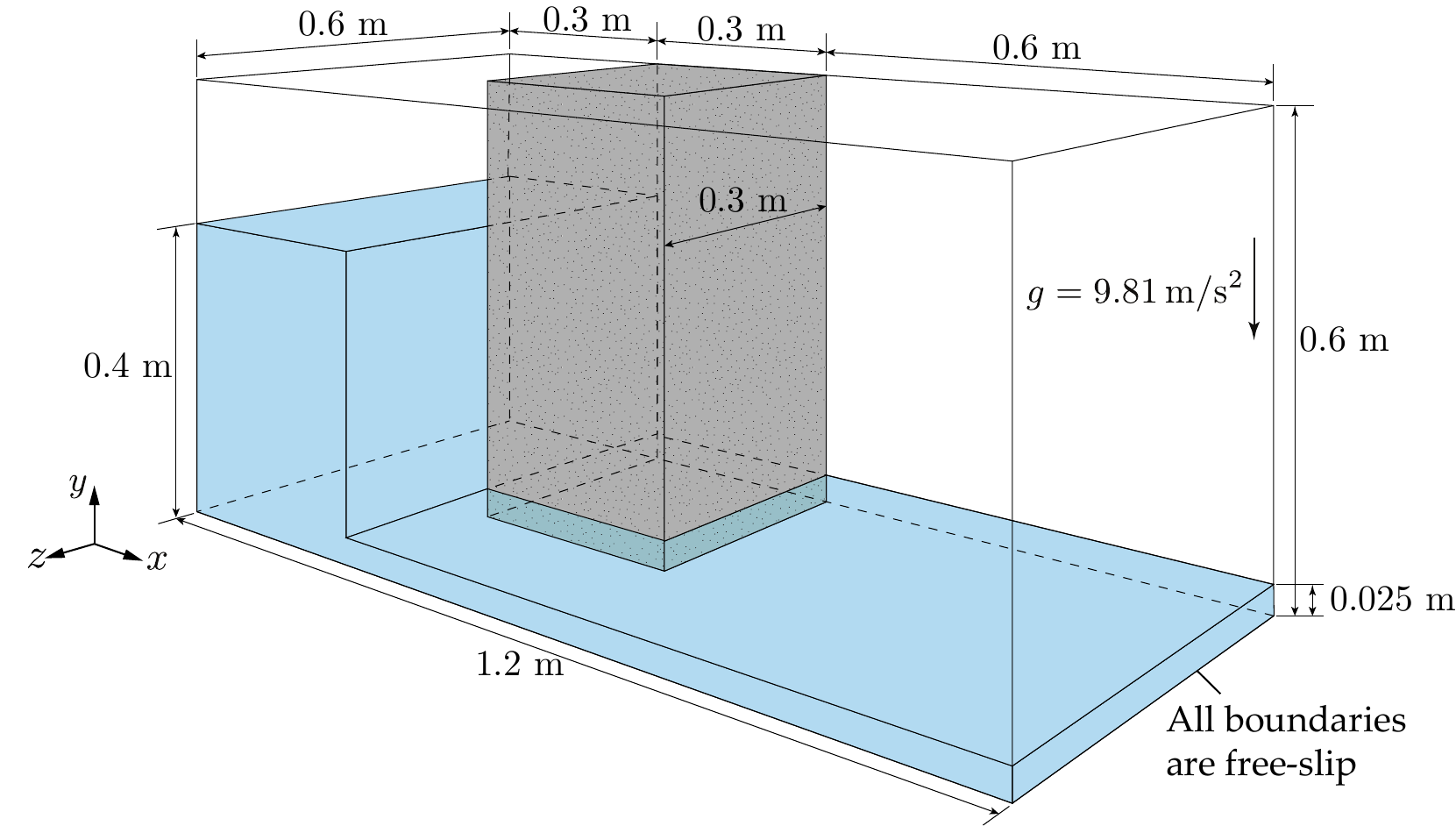}
    \caption{3D dam break of mixed free-surface--porous flow: geometrical settings and boundary conditions.}
    \label{fig:3.4_model}
\end{figure}

In the next validation test, we aim to demonstrate the capability of the proposed formulation in simulating complex three-dimensional fluid flows involving significant splashes, flow fragmentation, as well as seepage inside a rigid porous media. The geometrical settings of the problem are defined in Fig.~\ref{fig:3.4_model} following the numerical work of \citet{del2012three}. The porous medium is considered to be crushed rocks, which porosity, permeability, and other parameters defining the drag force model are set to be the same as the one utilized in the previous numerical example, see Table \ref{tab:3.3_cases} in Section \ref{subsec:2d_dam_break}. Similar to before, we also consider the effective viscosity $\mu_e$ inside the porous media to be the same as the free-fluid viscosity $\mu$. The gravitational acceleration $g=9.81$ $\mathrm{m/s^2}$ is considered and all side walls are set to be free-slip.

The MPM background grid is selected to be constructed of structured hexahedral elements with a cell size of $h=0.025$ m, where the material points are generated with $4\times4\times4$ PPC configuration. This results in a total of 27,648 background cells, 30,625 nodes, and 350,208 fluid material points. The initial volume of material points is assumed to be uniform. Moreover, the time step for this problem is set to be $\Delta t=0.002$ s, where all of the FLIP, APIC, and TPIC transfer schemes are considered for comparison purposes. 

\begin{figure}[htbp]
  \centering
   \includegraphics[width=0.8\textwidth]{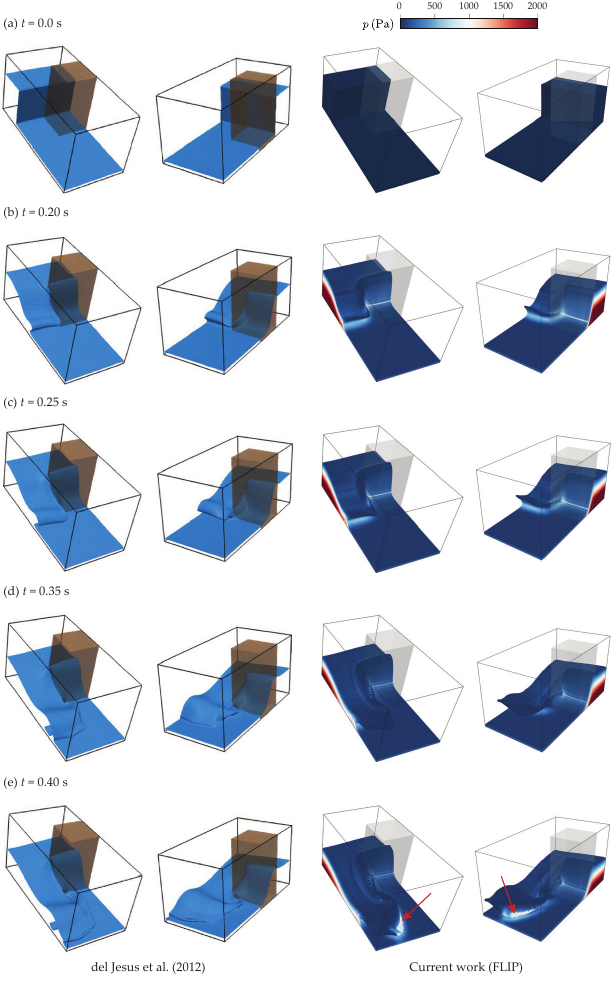}
    \caption{3D dam break of mixed free-surface--porous flow: snapshot comparison between the obtained results (right) and the results performed by \citet{del2012three} using VOF-FVM. Red arrows highlight the increase in pressure as the breaking wave crashes into the still water.}
    \label{fig:3.4_results_FLIP_a}
\end{figure}

\begin{figure}[htbp]
    \ContinuedFloat
    \centering
     \includegraphics[width=0.8\textwidth]{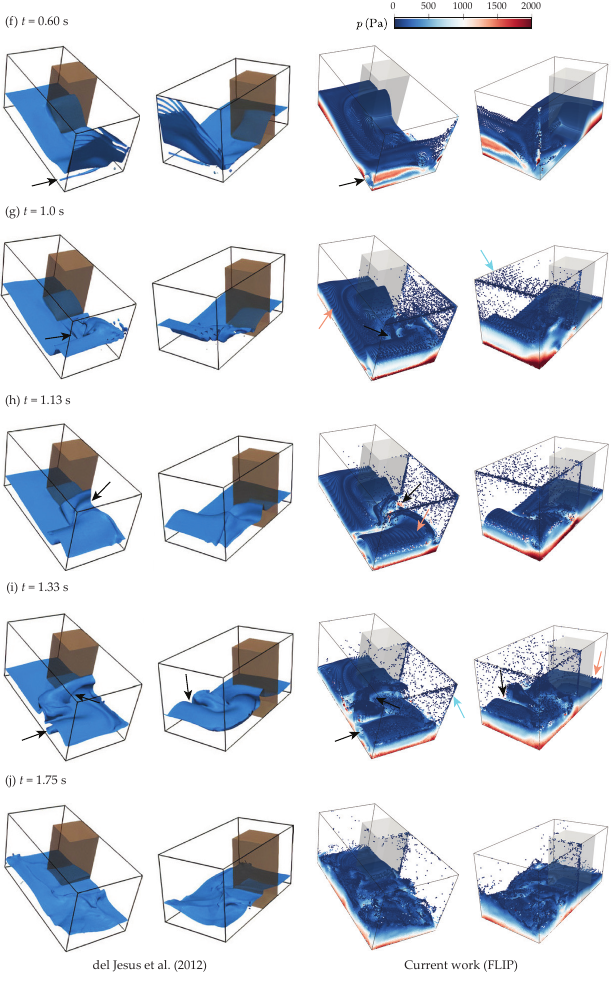}
    \caption{(continued). Light red arrows indicate the structure-like pattern resulting from the initial particle discretization. Light blue arrows emphasize water particles adhering to the upper and side walls due to normal wall constraints. Black arrows highlight similar free-surface flow features between the two models.}
    \label{fig:3.4_results_FLIP_b}
\end{figure}

The time evolution of the dam-break simulation, considering the FLIP scheme, is presented in Fig.~\ref{fig:3.4_results_FLIP_a} and qualitatively compared with the simulation results conducted by \citet{del2012three}. Overall, the mixed MPM results with the FLIP scheme exhibit similar behavior to the VOF-FVM results. Both models effectively capture similar free-surface flow features, such as the presence of vortical cavities after wave breaking, which are highlighted by the black arrows. From Fig.~\ref{fig:3.4_results_FLIP_a}, three notable observations arise. First, the obtained pressure profile demonstrates a consistently smooth and stable behavior throughout the simulation. However, since the MPM utilizes shared nodal kinematic fields during self-contact, early increases in the pressure field can be seen during the collapsing process of the wave crest into the still water (see Fig.~\ref{fig:3.4_results_FLIP_a}e as highlighted by the red arrows). Secondly, following the initial impact on the front wall and the subsequent upward movement, several fluid material points stuck to the upper and side walls (highlighted by the light blue arrows). Given that all walls are modeled as slip walls, with fixed normal displacement and velocity, this imposed boundary condition prevents fluid material points from immediately separating from the walls. However, over time, these adhering material points are expected to drop downward due to gravity. Lastly, owing to the consideration of structured particle spacing in the initialization process, structured-like striations are visible within the fluid body, particularly near the surface (indicated by the light red arrows). As has been previously discussed \cite{chandra2023stabilized}, this attribute is yet inevitable, as the employed $\delta$-correction method is designed to alleviate particle clumping rather than to correct overall numerical density and equalize particle spacing.

\begin{figure}[h!]
  \centering
   \includegraphics[width=0.45\textwidth]{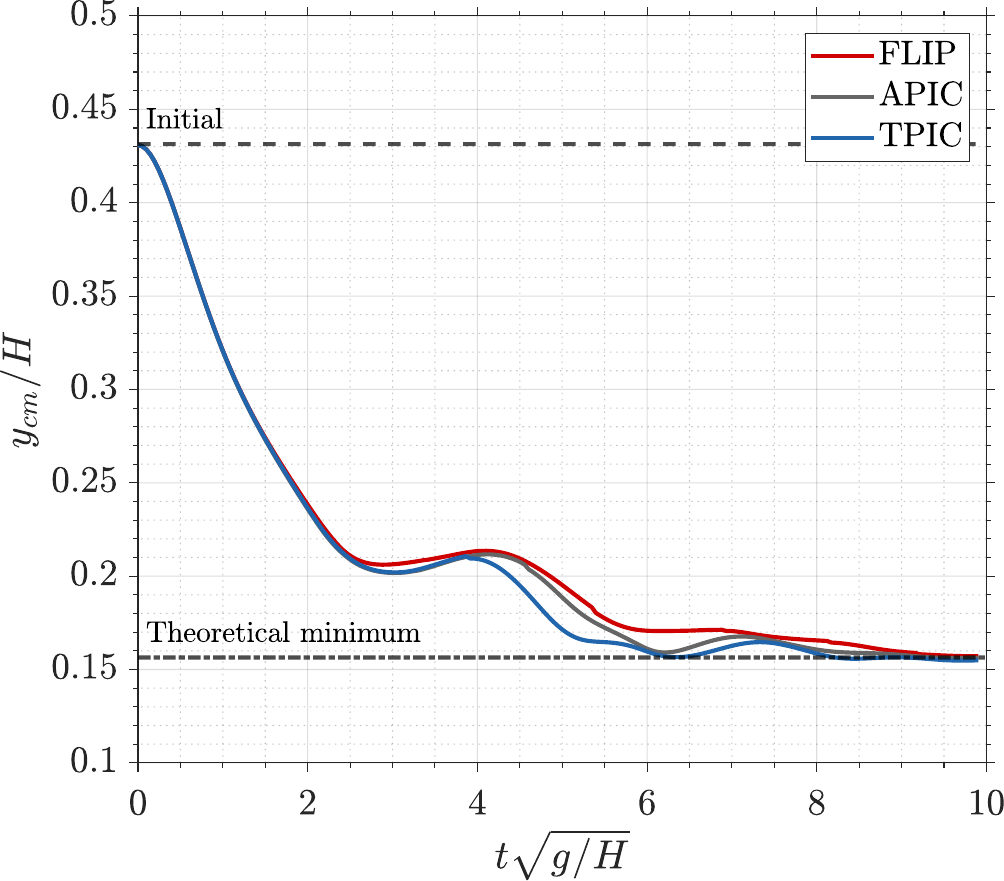}
    \caption{3D dam break of mixed free-surface--porous flow: $y$-coordinate of the computed center of mass over normalized time.}
    \label{fig:3.4_com}
\end{figure}

\begin{figure}[h!]
  \centering
   \includegraphics[width=0.85\textwidth]{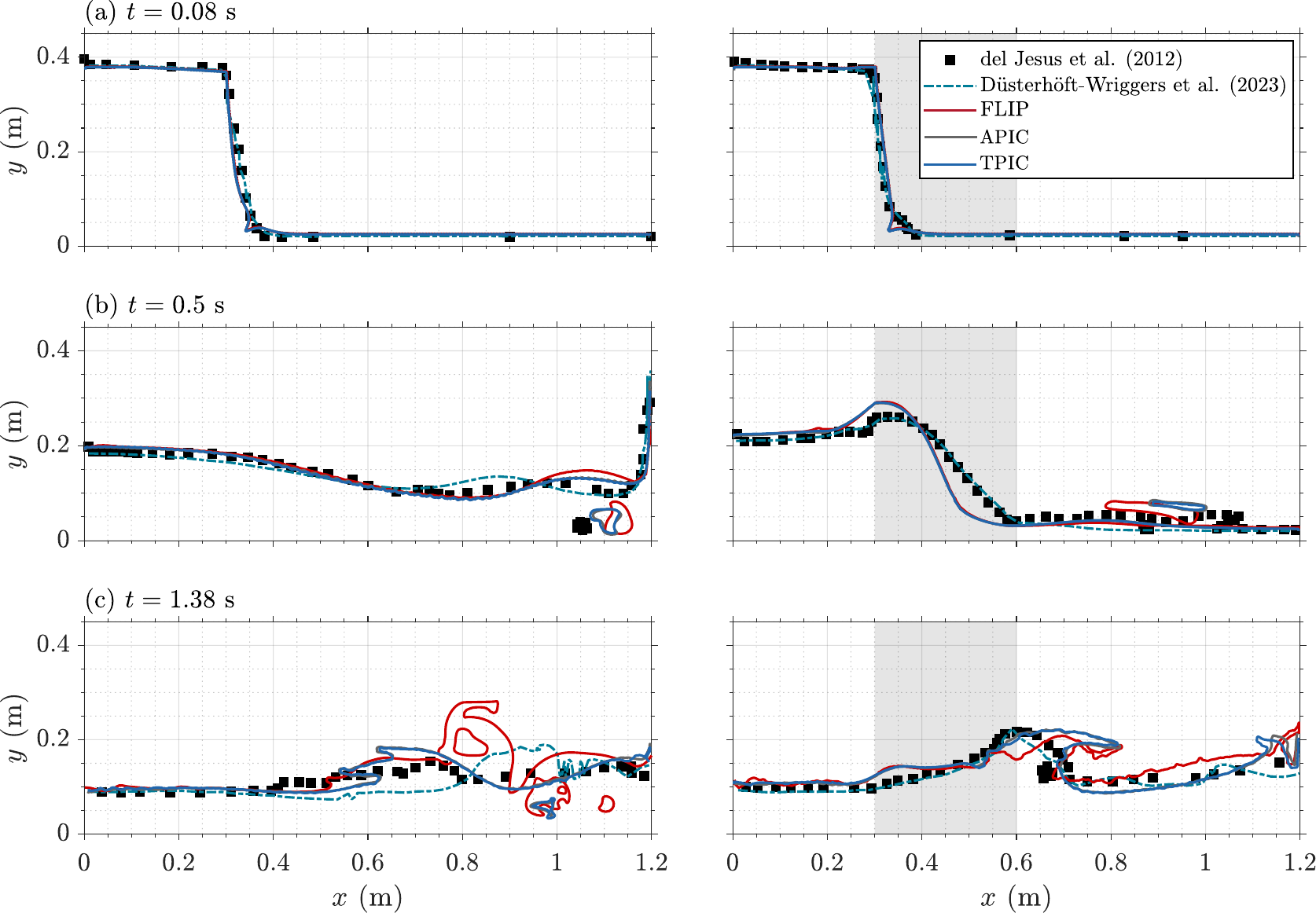}
    \caption{3D dam break of mixed free-surface--porous flow: comparison of free-surface profile at two different planes: (right) $z=0.1$ m and (left) $z=0.5$ m, between the obtained results of mixed MPM considering FLIP, APIC, and TPIC schemes and numerical results of other studies: \citet{del2012three} (VOF-FVM) and \citet{dusterhoft2023assessment} (VOF-FVM).}
    \label{fig:3.4_results_boundaries_compare}
\end{figure}

Next, the results obtained by the three different transfer schemes: FLIP, APIC, and TPIC, are compared. First, the overall volume conservation is checked by plotting the $y$-coordinate of the center of mass, which can be computed by Eq.~\eqref{eq:center_of_mass}. Considering the initial geometry and porosity distribution, we can calculate the theoretical limits of the mass center at the initial configuration and as the fluid approaches the resting state, which are approximately equal to $0.4314H$ and $0.1564H$, where $H=0.4$ m is the initial height of the water column. As shown by Fig.~\ref{fig:3.4_com}, the performed simulations show convergence towards the theoretical minimum value; this also verifies the performance of the $\delta$-correction method for 3D problems.

Furthermore, Fig.~\ref{fig:3.4_results_boundaries_compare} presents the free-surface profile in 2D slices at $z=0.1$ and $0.5$ m planes. Here, the obtained results are compared with the simulations performed using VOF-FVM done by \citet{del2012three} and \citet{dusterhoft2023assessment}. While there are some differences in the free-surface topology, particularly inside the porous domain as they utilized different drag-force models, the overall behaviors of all the simulations are in good agreement. Since the free-surface lines are drawn following the material point positions, the presence of vortical cavities as well as splashes during wave breaking can be tracked. Among the three MPM transfer schemes, it is obvious that the FLIP scheme generates a more splashy fluid motion. Meanwhile, the free-surface profiles obtained by the APIC and TPIC schemes are similar. As highlighted before in several previous works \cite{fu2017polynomial, fei2021revisiting, chandra2023stabilized}, the APIC and TPIC schemes produce a more apparent viscous behavior, which is shown by the free-surface profile that does not fragment and produce splashes as observed in the FLIP case.

\begin{figure}[h!]
  \centering
   \includegraphics[width=0.85\textwidth]{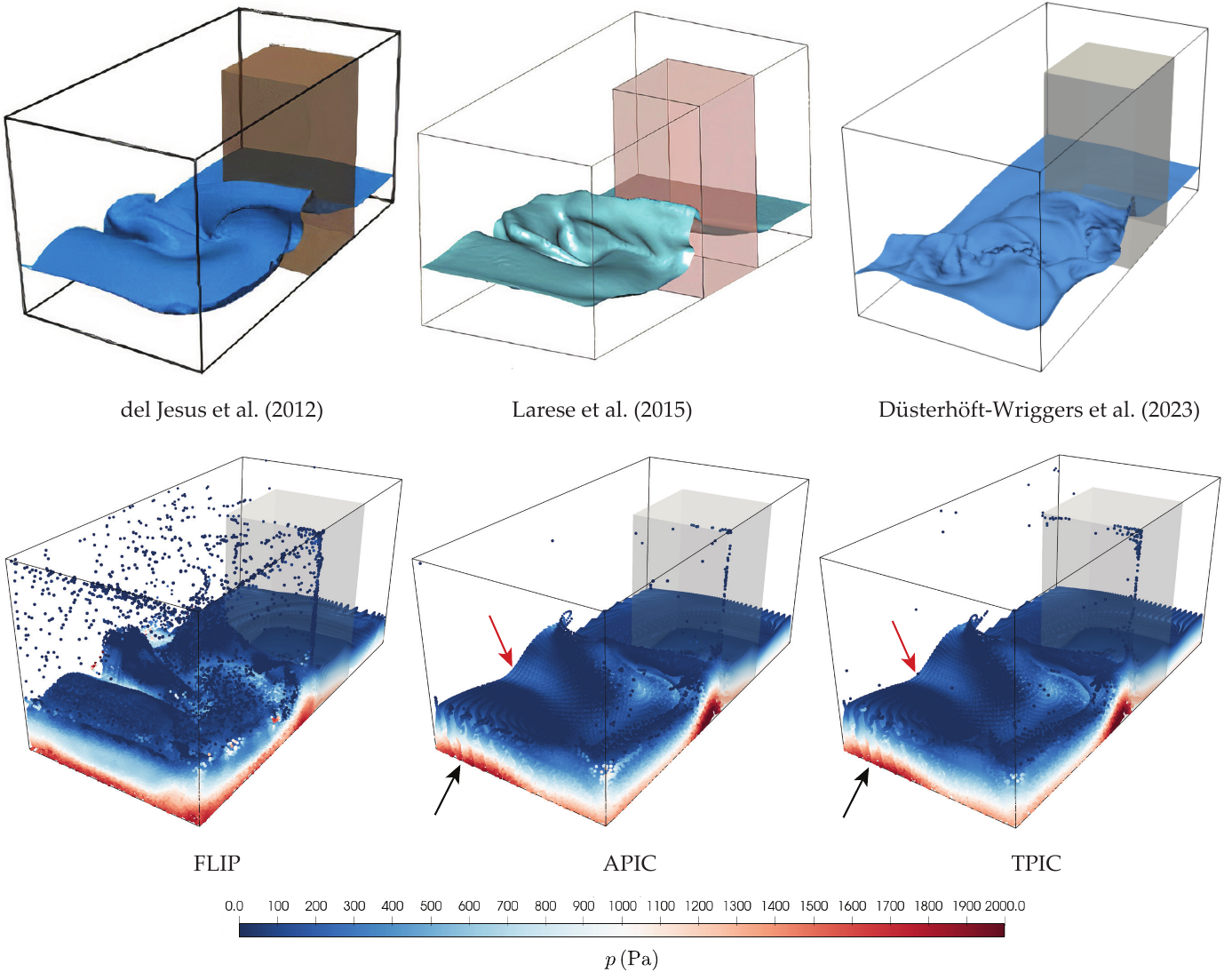}
    \caption{3D dam break of mixed free-surface--porous flow: comparison of free-surface profile at $t=1.33$ s between the obtained results of mixed MPM considering FLIP, APIC, and TPIC schemes (bottom row) and numerical results of other studies (top): \citet{del2012three} (VOF-FVM), \citet{larese2015finite} (edge-based FEM), and \citet{dusterhoft2023assessment} (VOF-FVM). Red arrows emphasize the relatively smooth free surface topology obtained by the APIC and TPIC schemes. Black arrows point out the structure-like pattern resulting from the initial particle discretization.}
    \label{fig:3.4_results_compare}
\end{figure}

While this attribute enhances numerical stability, the diffusive nature of the results often does not align visually with fluid flows observed in reality and experiments \cite{chandra2023stabilized}. The same features can be further observed in 3D through Fig.~\ref{fig:3.4_results_compare} which plots the free-surface profile at $t=1.33$ s. Here, in addition to the aforementioned numerical works, the result of edge-based FEM conducted by \citet{larese2015finite} is also presented. From Fig.~\ref{fig:3.4_results_compare}, the APIC and TPIC schemes are shown to produce a smoother free-surface profile compared to the other numerical results. The apparent viscosity of the two schemes also reduces the mixing of material points, and hence, a structure-like pattern can still be observed even after wave breaking. Enhancement methods, such as the integration of a combined APIC/FLIP or TPIC/FLIP approach, as well as the utilization of high-order mapping techniques like PolyPIC \cite{fu2017polynomial}, could potentially mitigate the observed viscosity issues. However, additional investigations into this aspect are necessary to achieve a more comprehensive understanding.

\section{Conclusions}
\label{sec:conclusions}

In this paper, we have developed a stabilized mixed MPM formulation to accurately and robustly model concurrent free-surface and seepage flow. The proposed framework is derived from mass and momentum balances that couple the Navier-Stokes equation with the Darcy-Brinkman-Forchheimer equation in a unified manner. To solve the balance equations, a mixed displacement-pressure coupled framework is derived by considering consistent linearization of the nonlinear drag force model. Here, a stabilization strategy employing the VMS method is incorporated into the proposed framework to satisfy the discrete inf-sup stability condition. To provide a seamless transition from free-surface flows to flows in porous media, $C^1$-continuous porosity and permeability fields are constructed using the blurred interface technique. This technique enhances numerical stability by reducing noise associated with sudden jumps in porosity. In addition, several other stabilization methods, such as the B-Spline basis function, the kernel correction method, and the $\delta$-correction method, have been incorporated into the formulation to reduce numerical artifacts associated with spatial discretization and particle quadrature. Specifically, we are able to improve numerical stability by iteratively performing the kernel-correction method to avoid negative nodal mass.

The formulation presented in this work has undergone extensive validation through a series of verification and validation tests. The 1D verification suite of gravity-driven flow confirms the accuracy, stability, and efficiency of the proposed method in comparison to the fractional-step MPM approach. Through the conducted tests, we can confirm the volume conservation property and the quadratic convergence of the Newton-Krylov solver. In certain conditions involving small permeability porous media, our method demonstrates a significant speedup ($>10,000$ times) compared to the fractional-step method, which is constrained by small critical time increments. Additionally, the obtained pressure profile is significantly smoother as the mixed MPM approach circumvents errors associated with free-surface detection and pressure imposition. From comparisons with (semi-)analytical solutions, the error of the numerical results is quadratically converged with mesh refinement in the case of pure seepage flow analysis. However, when the blurred interface is considered, the convergence rate decreases to a lower trend, though smoother and more stable results can be obtained. Furthermore, 2D Couette flow in a composite channel is simulated to demonstrate the capability of modeling accurate momentum transfer of shear flow in viscous fluid across the blurred porous interface. Again, our results have shown consistent performances in producing accurate results for different simulation settings. Lastly, the 2D and 3D dam break simulations are performed to validate the stability of the computed pressure field and the accuracy of free-surface topology compared with experimental data and numerical results from several previous studies.

For future work, it is necessary to develop and implement different types of conforming and nonconforming BCs, as well as a more detailed investigation of porous and free-surface interface conditions, e.g.~considering the tangential stress jump condition or the capillary suction condition. Additionally, ongoing research involves extending the proposed solver to incorporate history-dependent elasto-plastic solids, which may undergo large deformations and experience damage. For further validation and practical applications, emphasis should be placed on enhancing the efficiency and scalability of the mixed MPM solver for simulating large-scale 3D scenarios with millions of particles. This necessitates research efforts in the efficient implementation of dynamic load balancing and fast matrix assembly.

\section*{CRediT author contribution statement}
\textbf{B. Chandra:} Conceptualization, Methodology, Software, Investigation, Validation, Formal analysis, Visualization, Writing - original draft, review, and editing. \textbf{R. Hashimoto}: Conceptualization, Methodology, Supervision, Writing - review and editing. \textbf{K. Kamrin:} Formal analysis, Resources, Supervision, Writing - review and editing. \textbf{K. Soga:} Resources, Supervision, Writing - review and editing.

\section*{Acknowledgments}
This work is partially funded by the Donald H. McLaughlin Chair Fund at UC Berkeley. R. Hashimoto acknowledges the support of the Japan Society for the Promotion of Science, Fostering Joint International Research (A), Grant Number 21KK0254. This research used the Savio computational cluster resource provided by the Berkeley Research Computing program at the University of California, Berkeley (supported by the UC Berkeley Chancellor, Vice Chancellor for Research, and Chief Information Officer). B. Chandra thanks Hasitha Wijesuriya from UC Berkeley for their help in setting up the computational environment in Savio. The authors also acknowledge the assistance of students from Hiroshima University, Shinnosuke Matsumi and Kohei Mihara, supervised by R. Hashimoto in prototyping the numerical framework. 

\section*{Declaration of competing interest}
The authors declare that they have no conflict of interest.

\section*{Data availability statement}
The data that support the findings of this study are available from the corresponding authors upon reasonable request.

\appendix

\section{Equivalent drag force coefficients expressions for different input parameters}
\label{app:drag_force_coeff}

This appendix provides easy conversion and computation of the drag force coefficients given previously in Eq.~\eqref{eq:drag_force_simplifed}. Assuming the fluid is a viscous fluid with non-zero viscosity $\mu$, the coefficients $\widetilde{A}$ and $\widetilde{B}$ are given by default in this paper as:
\begin{eqnarray}
	{\widetilde{A}} = \frac{\theta \mu}{\rho k}\,, \qquad {\widetilde{B}} = {\frac{B \sqrt{\theta}}{\sqrt{Ak}}}\,,
	\label{eq:drag_force_coefficients_form_1}
\end{eqnarray}
where Kozeny-Carman\cite{carman1937fluid} and Ergun's\cite{ergun1952fluid} assumptions are considered (see \Cref{subsec:balance_laws}). Here, the input material parameters are the intrinsic permeability, $k$ ($L^2$), as well as the two dimensionless coefficients $A$ and $B$, whose value can be found by empirical relationship or via a calibration process \cite{del2012three, akbari2013moving, larese2015finite, losada2016modeling}. In the case of different initial input parameters, a conversion can be performed through simple algebraic operations. Several common cases considered in the literature, such as by using the hydraulic conductivity $K$ ($LT^{-1}$) or effective grain size $d$ ($L$), are given in Table \ref{tab:drag_force_coeff}.

\begin{table}[h!]
\small
\centering
\caption{Drag force coefficient conversion from different input parameters.}
\label{tab:drag_force_coeff}
\renewcommand{\arraystretch}{2.0}
\begin{tabular}{||c|c|c|c||}
\hline
\multirow{2}{*}{Term} & Input 1                                                                     & Input 2                                                                               &Input 3\\ \cline{2-4} 
                      & $k_{\mathrm{ref}}$($L^2$), $\theta_{\mathrm{ref}}$(-)                                                  & $K_{\mathrm{ref}}$($LT^{-1}$), $\theta_{\mathrm{ref}}$(-)                                               & $d_{\mathrm{ref}}$($L$)\\ \hline \hline
Kozeny-Carman         & $k = \dfrac{\theta^3}{\theta_{\mathrm{ref}}^3} \dfrac{(1-\theta_{\mathrm{ref}})^2}{(1-\theta)^2}k_{\mathrm{ref}}$ & $K = \dfrac{\theta^3}{\theta_{\mathrm{ref}}^3} \dfrac{(1-\theta_{\mathrm{ref}})^2}{(1-\theta)^2}K_{\mathrm{ref}}$  &$k=\dfrac{\theta^3}{(1-\theta)^2} \dfrac{d_{\mathrm{ref}}^2}{A}$\\ \hline
Linear, $\widetilde{A}$   & $\dfrac{\theta \mu}{\rho k}$                                                          & $\dfrac{\theta g}{K}$                                                                  &$ \dfrac{A \mu (1-\theta)^2}{ \theta^2 \rho d_{\mathrm{ref}}^2}$\\ \hline
Nonlinear, $\widetilde{B}$ & $\dfrac{B\sqrt{\theta}}{\sqrt{A k}}$                                                  & $\dfrac{B\sqrt{\theta \rho g}}{\sqrt{A K\mu}}$                                         &${\dfrac{B (1-\theta)}{ \theta d_{\mathrm{ref}}}}$\\ \hline
\end{tabular}
\end{table}

\section{Quadratic B-Spline basis functions and iterative kernel correction method}
\label{app:basis_funct_it_kc}

The internal zero-centered quadratic B-Spline basis function can be written in a closed-form expression as suggested by \citet{steffen2008analysis} as:
\begin{eqnarray}
    S_{Ip} &:=& S_I(\tb x^n_p) =  \prod_{d = 1}^{\rm dim} w\left(\frac{x^n_{p,d} - x_{I,d}}{h}\right) \,, \label{eq:bspline_func}\\
    w(\xi_d) &:=&
    \begin{dcases} 
      -|\xi_d|^2+\frac{3}{4}\,, & \qquad 0 \leq |\xi_d| < \frac{1}{2}\,,\\
      \frac{1}{2} |\xi_d|^2 - \frac{3}{2} |\xi_d| + \frac{9}{8}\,, & \qquad \frac{1}{2}\leq |\xi_d| < \frac{3}{2}\,,\\
      0\,, & \qquad |\xi_d| \geq \frac{3}{2}\,.
    \end{dcases}  
\end{eqnarray}
Here, $h$ and $\xi_d$ denote the cell spacing and the relative distance of direction $d$, respectively. Correspondingly, we can derive the direction-wise gradient as:
\begin{eqnarray}
    \frac{\partial S_{Ip}}{\pd x_{\alpha}} &:=&  \frac{\partial S_{I} (\tb x_p^n)}{\pd x_{\alpha}} = \left( \prod_{d \neq \alpha}^{\rm dim} w\left(\frac{x^n_{p,d} - x_{I,d}}{h}\right) \right) \frac{\pd w\left(\frac{x^n_{p,\alpha} - x_{I,\alpha}}{h}\right)}{\pd x_\alpha} \,, \label{eq:bspline_func_grad}\\
        \frac{\pd w\left(\xi_\alpha\right)}{\pd x_\alpha} &:=&
    \begin{dcases} 
      \frac{-2 \xi_\alpha}{h}\,, & \qquad 0 \leq |\xi_\alpha| < \frac{1}{2}\,,\\
      \frac{\mathrm{sgn}(\xi_\alpha)}{h}\left(|\xi_\alpha| - \frac{3}{2} \right)\,, & \qquad \frac{1}{2}\leq |\xi_\alpha| < \frac{3}{2}\,,\\
      0\,, & \qquad |\xi_\alpha| \geq \frac{3}{2}\,.
    \end{dcases}  
\end{eqnarray}


To satisfy both the partition of unity and linear field reproduction conditions, \citet{nakamura2023taylor} proposed the kernel correction method based on the weighted least squares (WLS) approach. The method is summarized as follows. First, the correction function $\widehat{C}$ is constructed using the WLS method as:
\begin{eqnarray}
    \widehat{C}(\tb x, \tb x_I):=\tb P(\tb x - \tb x_p) \tb M^{-1} \tb P^T(\tb x - \tb x_p)\,, 
    \label{eq:kernel_correction_function}
\end{eqnarray}
where the moment matrix $\tb M$ and linear polynomial basis $\tb P$ are expressed as:
\begin{eqnarray}
    \tb M:=\tb M(\tb x_p)=\sum_{I=1}^{n_n} S_{Ip} \tb P^T(\tb x_I - \tb x_p) \tb P(\tb x_I - \tb x_p)\,, \qquad \tb P(\tb x) = [1,\,x,\,y,\,z]\,\quad\mathrm{for}\,\mathrm{dim}=3\,. 
    \label{eq:moment_matrix}
\end{eqnarray}
The corrected basis function given can be obtained by multiplying the correction function with the standard basis function as:
\begin{eqnarray}
    \widehat{S}_{Ip} = \widehat{C}(\tb x_p, \tb x_I) {S}_{Ip}\,.
\end{eqnarray}
Knowing that the inverse of moment matrix $\tb M^{-1}$ can be expressed as a $(\mathrm{dim}+1)\times (\mathrm{dim}+1)$ block matrix:
\begin{eqnarray}
    \tb M^{-1}=
    \begin{bmatrix}
        C_1 & \tb C_2^T\\
        \tb C_2 & \tb C_3
    \end{bmatrix}\,,
    \label{eq:moment_matrix_inv}
\end{eqnarray}
the corrected basis functions $\widehat{S}_{Ip}$ and their gradient $\nabla \widehat{S}_{Ip}$ can be computed as \cite{nakamura2023taylor}:
\begin{eqnarray}
    \begin{Bmatrix}
        \widehat{S}_{Ip} \\
        \nabla \widehat{S}_{Ip}
    \end{Bmatrix} := \mathcal{K}_c (S_{Ip}) =  
    \begin{Bmatrix}
        C_1 + \tb C_2 \cdot \left(\tb x_I - \tb x_p \right) \\
        \tb C_2 + \tb C_3 \cdot \left(\tb x_I - \tb x_p \right)
    \end{Bmatrix} S_{Ip}\,,
   \label{eq:kernel_correction_function_2}
\end{eqnarray}
where $\mathcal{K}_c:\mathcal{R}\rightarrow\mathcal{R} \times \mathcal{R}^{\rm{dim}}$ is defined as the kernel-correction operator. It should be noted that the kernel correction procedure should be done only at the region where the partition of unity is not satisfied, e.g.~near the domain boundary. 

Upon further investigation, it has come to our attention that the WLS kernel correction method, as presented by \citet{nakamura2023taylor}, does not ensure non-negative values for basis functions even though the uncorrected kernels computed by Eq.~\eqref{eq:bspline_func} are non-negative. This issue arises due to the lack of a guarantee that the kernel-corrector operator is always positive. Specifically, the inequality $\left(C_1 + \tb{C}_2 \cdot (\tb{x}_I - \tb{x}_p)\right) \geq 0$ does not always hold. As highlighted by \citet{andersen2007material}, the utilization of negative basis functions is incompatible with the MPM formulation and may lead to significant issues, particularly when mass on a background node becomes negative. While such conditions are relatively rare in typical 2D problems, the likelihood of its occurrence significantly rises in scenarios involving splashes and fragmentation, such as those encountered in 3D free-surface fluid dynamics problems, e.g.~\Cref{subsec:3d_dam_break}. 

\begin{figure}[h!]
  \centering
   \includegraphics[width=1.0\textwidth]{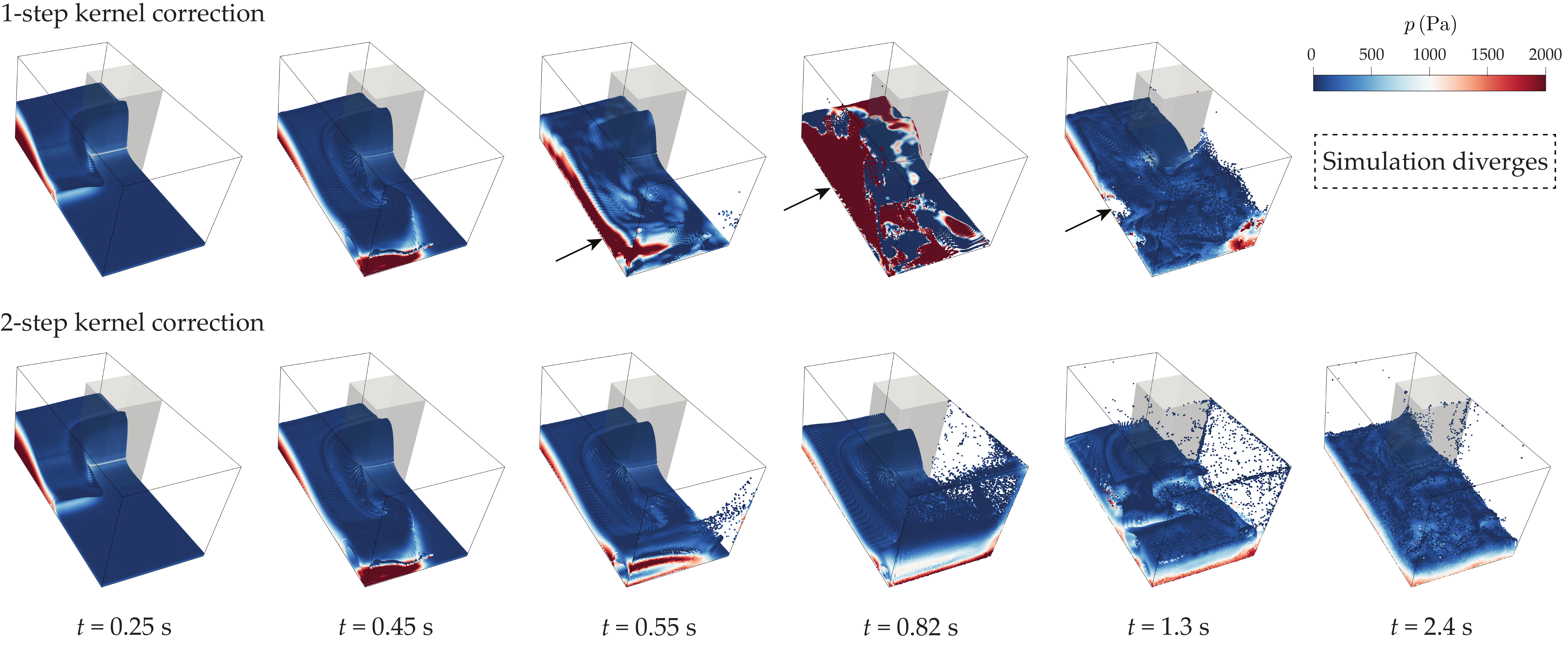}
    \caption{Snapshots of 3D dam-break simulation results performed with (top row) 1-step and (bottom row) 2-step kernel correction algorithm, all considering FLIP transfer scheme. Arrows highlight spurious pressure oscillations and numerical explosions.}
    \label{fig:B_it_kc}
\end{figure}

In the current work, a simple mitigation technique to avoid negative basis functions is considered. Here, we employ an iterative kernel correction algorithm, which performs the WLS kernel correction successively if there are any negative basis functions generated by the prior kernel correction routine. The algorithm is summarized in \Cref{al:iterative_kc} while Fig.~\ref{fig:B_it_kc} presents snapshots comparing simulation results without and with the iterative kernel correction algorithm for the 3D dam break problems discussed in Section \ref{subsec:3d_dam_break}. As evident from the figures, the iterative kernel correction method effectively maintains numerical stability, especially in scenarios where splashes are generated. 

\begin{algorithm}
\caption{Iterative WLS kernel correction algorithm.}
\textbf{Input:} \texttt{max\_iteration}. \Comment{Defines the maximum number of iterations to avoid an infinite loop, usually set as 2 or 3. Setting it as 1 recovers the kernel correction proposed by \citet{nakamura2023taylor}.}\\
\textbf{Output:} $\widehat{S}_{Ip}$ and $\nabla \widehat{S}_{Ip}$.
\begin{algorithmic}[1]

\State Compute quadratic B-Spline kernel $S_{Ip}$ and its gradient $\nabla S_{Ip}$: Eqs.~\eqref{eq:bspline_func} and \eqref{eq:bspline_func_grad}.

\If{$\sum_I S_{Ip} \neq 1$}
    \State Set \texttt{kernel\_correction} $\leftarrow$ \textit{true}.
\Else
    \State Return $S_{Ip}$ and $\nabla S_{Ip}$.
\EndIf

\State Set $k=0$, $\widehat{S}^0_{Ip} = S_{Ip}$, and $\nabla \widehat{S}^0_{Ip}=\nabla S_{Ip}$.
\While{\texttt{kernel\_correction}}
    \State Compute the corrected kernel $\widehat{S}^{k+1}_{Ip}$ and gradient $\nabla \widehat{S}^{k+1}_{Ip}$: Eqs.~\eqref{eq:moment_matrix} - \eqref{eq:kernel_correction_function_2}.
    \State Set \texttt{kernel\_correction} $\leftarrow$ \textit{false} and append $k=k+1$.
    \If{$k<\texttt{max\_iteration}$}
    \For{all connectivity $I=1,...,n_n$}
    \If{$\widehat{S}^{k}_{Ip} < 0.0$} 
    \State Set $\widehat{S}^{k}_{Ip}=0$ and \texttt{kernel\_correction} $\leftarrow$ \textit{true}.
    \EndIf
    \EndFor
    \EndIf
    
\EndWhile

\State Return $\widehat{S}^k_{Ip}$ and $\nabla \widehat{S}^k_{Ip}$.
    
\end{algorithmic}
\label{al:iterative_kc}
\end{algorithm}

\section{Analytical solutions of verification problems}
\label{app:analytical_solutions}

\subsection{1D gravity-driven flow - Case 3}
\label{app:1d_case3_derivation}
The analytical solution for the problem described in \Cref{subsubsec:1d_case3} can be derived directly from Eq.~\eqref{eq:equilibrium}. Since the fluid is assumed to free-flow inside a porous media, the stress gradient in all fluid parcels is zero, thus, the flow is governed entirely by gravity, as its driving force, and drag force, as the resisting force. One can rewrite Eq.~\eqref{eq:equilibrium} for one-dimensional case as:
\begin{eqnarray}
    a_y(t) = g - \widetilde{A} v_y(t) - \widetilde{B} |v_y(t)| v_y(t)\,.
    \label{eq:case3_acc_1d}
\end{eqnarray}
Integrating Eq.~\eqref{eq:case3_acc_1d} using analytical or numerical integration techniques, considering the given permeability conditions and with initial condition $v_y(0) = 0$ m/s, will yield the 1D velocity profile given in Fig.~\ref{fig:3.1.3_linear} and \ref{fig:3.1.3_nonlinear} for the Darcian and non-Darcian flow, respectively.

\subsection{1D gravity-driven flow - Case 4}
\label{app:1d_case4_derivation}

To derive the analytical solution for problems given in \Cref{subsubsec:1d_case4}, first, the governing equation of the free-fluid and pore-fluid, the mass and momentum balance, are written in the 1D-strong-form as:
\begin{alignat}{3}
    \frac{\pd v_{ff,y}}{\pd y}(t)&=0\,, \qquad& a_{ff,y}(t)&=g-\frac{p_i(t)}{h_{ff}(t) \rho} \,, &\qquad& \forall \tb x \in \Omega_{ff}\,,\label{eq:ff_governing_equations}\\
    \frac{\pd v_{fp,y}}{\pd y}(t)&=0\,, \qquad&  a_{fp,y}(t)&= g- \widetilde{A} v_{fp,y}(t) - \widetilde{B} |v_{fp,y}(t)| v_{fp,y}(t)+\frac{p_i(t)}{h_{fp}(t) \rho}\,, &\qquad& \forall \tb x \in \Omega_{fp}\,. \label{eq:pm_governing_equations}
\end{alignat}
Here, as previously given, $p_i$ is the interface pressure at $\Gamma_i$, whereas the subscripts $\square_{ff}$ and $\square_{fp}$ indicate the fluid variables in the non-porous and porous domain, respectively (see Fig.~\ref{fig:3.1_model} for notation).

To get into the solution, it is necessary to couple the two sets of equations given in Eqs.~\eqref{eq:ff_governing_equations} and \eqref{eq:pm_governing_equations}. Here, we consider a sharp and discontinuous porous interface as illustrated in Fig.~\ref{fig:2_porous_flow_MPM}a. The first coupling condition is to ensure mass conservation across the porous interface $\Gamma_i$, which is assumed to be permeable. This requires the continuity of normal velocity components as:
\begin{eqnarray}
    \tb v_{ff}\cdot \overline{\tb n} = \theta \tb v_{fp} \cdot \overline{\tb n} \qquad \mathrm{at}\, \Gamma_i\,,
    \label{eq:velocity_across_jump}
\end{eqnarray}
where $\overline{\tb n}$ is the normal vector pointing outward into the free-fluid region from the porous domain. Notice that, in virtue of Eqs.~\eqref{eq:ff_governing_equations}$_1$ and \eqref{eq:pm_governing_equations}$_1$, Eq.~\eqref{eq:velocity_across_jump} can be generalized for the whole fluid domain for this case as:
\begin{eqnarray}
    v_{ff,y}(t) = \theta v_{fp,y}(t)\,.
    \label{eq:velocity_across_jump_1d}
\end{eqnarray}
Furthermore, the second coupling condition ensures the equilibrium of normal forces across the interface, i.e.
\begin{eqnarray}
    p_{ff} = p_{fp} := p_i\qquad \mathrm{at}\, \Gamma_i\,.
    \label{eq:pressure_across_jump}
\end{eqnarray}

Given the initial height of the free-fluid as $h_{ff}(0)=h_0$, the height of the free-fluid and porous fluid domain can be expressed in relation to the respective velocity fields as:
\begin{eqnarray}
    h_{ff}(t) = h_0 - \int_0^t v_{ff, y}(t) \td t\,, \qquad
    h_{fp}(t) = \int_0^t v_{fp, y}(t) \td t\,,
\label{eq:height_change}%
\end{eqnarray}
which can be combined knowing Eq.~\eqref{eq:velocity_across_jump_1d} as:
\begin{eqnarray}
     h_{ff}(t) = h_0 - \theta h_{fp}(t)\,.
\label{eq:height_change_connected}%
\end{eqnarray}
Note that, this condition is valid only when $h_{ff}(t) \geq 0$. Through Eqs.~\eqref{eq:velocity_across_jump_1d}, \eqref{eq:height_change}, and \eqref{eq:height_change_connected}, all the velocity and acceleration terms written in Eqs.~\eqref{eq:ff_governing_equations} and \eqref{eq:pm_governing_equations} can be written in terms of the ODE of $h_{ff}$ as:
\begin{subequations}
    \begin{eqnarray}
        v_{ff,y}(t) = - \dot{h}_{ff}(t)\,, &\qquad& a_{ff,y} (t)= - \ddot{h}_{ff}(t)\,,\\
        v_{fp,y}(t) = - \dot{h}_{ff}(t)/{\theta}\,, &\qquad& a_{fp,y}(t) = - {\ddot{h}_{ff}}(t)/{\theta}\,.
    \end{eqnarray}
    \label{eq:kinematic_derivation_1d}%
\end{subequations}
Substituting the kinematic expressions given in Eq.~\eqref{eq:kinematic_derivation_1d} to Eqs.~\eqref{eq:ff_governing_equations}$_2$ and \eqref{eq:pm_governing_equations}$_2$, and equating the interface pressure $p_i$ from the two equations, results in the following nonlinear ODE:
\begin{eqnarray}
    h_{ff}\left(g+\ddot{h}_{ff}\right)+\left(\frac{h_0 - h_{ff}}{\theta}\right)\left(g+ \frac{\ddot{h}_{ff}}{\theta} + \left(\widetilde{A}+\widetilde{B}\frac{|\dot{h}_{ff}|}{\theta}\right)\frac{\dot{h}_{ff}}{\theta}\right)=0\,,
    \label{eq:case4_nonlinear_ode}
\end{eqnarray}
which can be solved by using a numerical ODE solver. Here, the interface pressure $p_i$ can be computed subsequently upon obtaining the height of free fluid $h_{ff}(t)$ as:
\begin{eqnarray}
    p_i(t) = \rho h_{ff}(t)\left( g + \ddot{h}_{ff}(t) \right)\,.
\end{eqnarray}
Considering the material and initial conditions given in \Cref{subsubsec:1d_case4}, the semi-analytical results plotted in Fig.~\ref{fig:3.1.4_results} can be obtained.

\subsection{1D gravity-driven flow - Case 5}
\label{app:1d_case5_derivation}

The analytical solution for the next problem, \Cref{subsubsec:1d_case5}, can be done similarly to \Cref{subsubsec:1d_case4}. Following the previous procedure, the mass and momentum balance of the pore-fluid and free-fluid are expressed in 1D-strong-form as (see Fig.~\ref{fig:3.1_model} for notation):
\begin{alignat}{3}
    \frac{\pd v_{fp,y}}{\pd y}(t)&=0\,, \qquad&  a_{fp,y}(t)&= g- \widetilde{A} v_{fp,y}(t) - \widetilde{B} |v_{fp,y}(t)| v_{fp,y}(t)-\frac{p_i(t)}{h_{fp}(t) \rho}\,, &\qquad& \forall \tb x \in \Omega_{fp}\,, \label{eq:pm_governing_equations_2} \\
    \frac{\pd v_{ff,y}}{\pd y}(t)&=0\,, \qquad& a_{ff,y}(t)&=g+\frac{p_i(t)}{h_{ff}(t) \rho} \,, &\qquad& \forall \tb x \in \Omega_{ff}\,.\label{eq:ff_governing_equations_2}
\end{alignat}
Here, the compatibility of velocity and pressure fields at the porous interface $\Gamma_i$, which are described by Eqs.~\eqref{eq:velocity_across_jump_1d} and \eqref{eq:pressure_across_jump}, are valid.

Given the initial height of the pore-fluid as $h_{fp}(0)=h_0$, we can express the heights of the porous and free fluid domain as follows:
\begin{eqnarray}
    h_{fp}(t) &=& h_0 - \int_0^t v_{fp, y}(t) \td t\,, \qquad
    h_{ff}(t) = \int_0^t v_{ff, y}(t) \td t\,,
\label{eq:height_change_2}\\
     h_{fp}(t) &=& h_0 - \frac{h_{ff}(t)}{\theta}\,,
\label{eq:height_change_connected_2}%
\end{eqnarray}
taking into account Eq.~\eqref{eq:velocity_across_jump_1d} and considering the constraint $h_{fp}(t) \geq 0$. Through Eqs.~\eqref{eq:velocity_across_jump_1d}, \eqref{eq:height_change_2}, and \eqref{eq:height_change_connected_2}, we can rewrite all the velocity and acceleration terms in Eqs.~\eqref{eq:pm_governing_equations_2} and \eqref{eq:ff_governing_equations_2} as:
\begin{subequations}
    \begin{eqnarray}
        v_{fp,y}(t) = - \dot{h}_{fp}(t)\,, &\qquad& a_{fp,y}(t) = - {\ddot{h}_{fp}}(t)\,,\\
        v_{ff,y}(t) = - \theta \dot{h}_{fp}(t)\,, &\qquad& a_{ff,y} (t)= - \theta \ddot{h}_{fp}(t)\,.
    \end{eqnarray}
    \label{eq:kinematic_derivation_1d_2}%
\end{subequations}
By having Eq.~\eqref{eq:kinematic_derivation_1d_2} in hand, we can derive and solve the following nonlinear ODE following the process done in \ref{app:1d_case4_derivation}:
\begin{eqnarray}
    \theta \left( h_0 - h_{fp}\right) \left(g+ \theta \ddot{h}_{fp} \right) + h_{fp} \left( g + \ddot{h}_{fp} + \left(\widetilde{A} + \widetilde{B}|\dot{h}_{fp}| \right)\dot{h}_{fp} \right)=0\,.
    \label{eq:case5_nonlinear_ode}
\end{eqnarray}
Subsequently, the interface pressure $p_i$ can be computed as:
\begin{eqnarray}
    p_i(t) = -\rho \theta \left(h_0 - h_{fp}(t) \right)\left( g + \theta \ddot{h}_{fp}(t) \right)\,.
\end{eqnarray}
Considering the given material and initial conditions, the semi-analytical results previously plotted in Fig.~\ref{fig:3.1.5_results} can be obtained.

\begin{remark}
\label{rem:stiff_equation_case5}
The solution to the differential equation presented in Eq.~\eqref{eq:case5_nonlinear_ode} exhibits significant stiffness, characterized by a rapid increase in acceleration within a remarkably short time span. To illustrate the extent of this stiffness, the semi-analytical solution of velocity, acceleration, and jerk (the rate of change of acceleration) is plotted in Fig.~\ref{fig:C.3_anals}, which considers the material parameters considered in \Cref{subsubsec:1d_case5}. The acceleration, as shown in the plots, experiences a rapid and substantial increase at the end of the simulation. Notably, the jerk reaches exceptionally high magnitudes, which are indicative of the highly dynamic nature of the system, and it emphasizes the challenge of numerically solving this ODE accurately.

\begin{figure}[h!]
    \centering
     \begin{subfigure}[b]{0.315\textwidth}
         \centering
         \includegraphics[width=\textwidth]{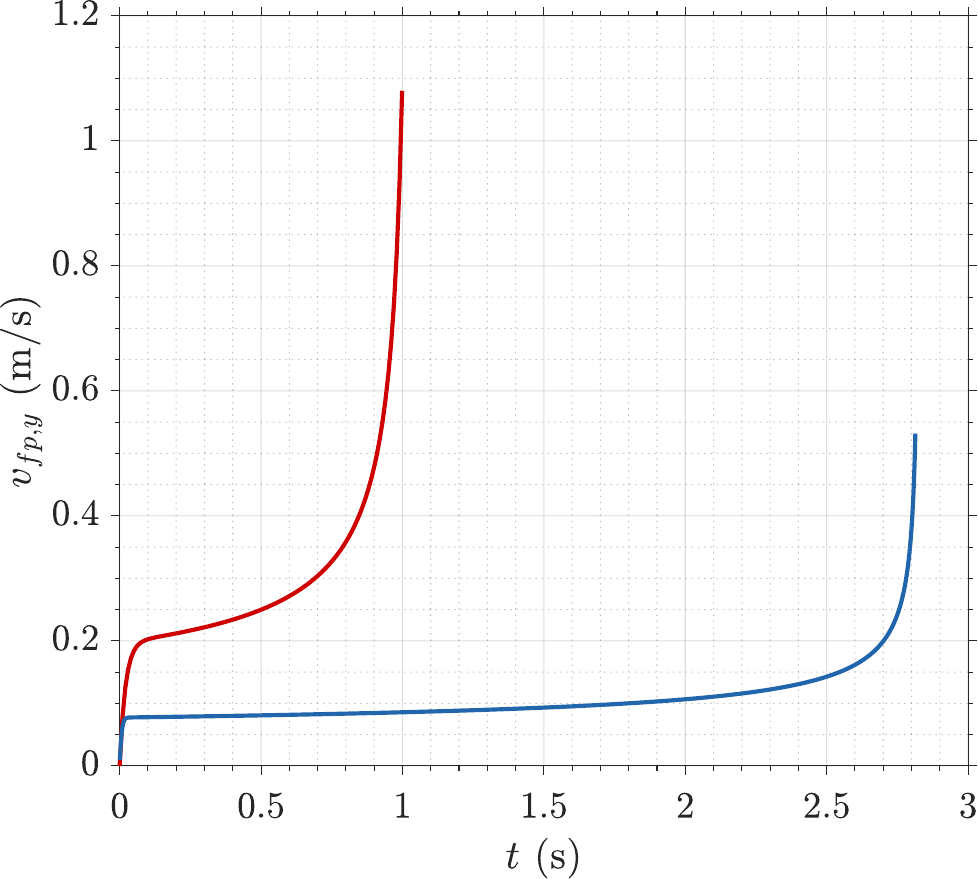}
         \caption{}
         \label{fig:C.3_anal_vy}
     \end{subfigure}
     \hspace{0.1cm}
     \begin{subfigure}[b]{0.312\textwidth}
         \centering
         \includegraphics[width=\textwidth]{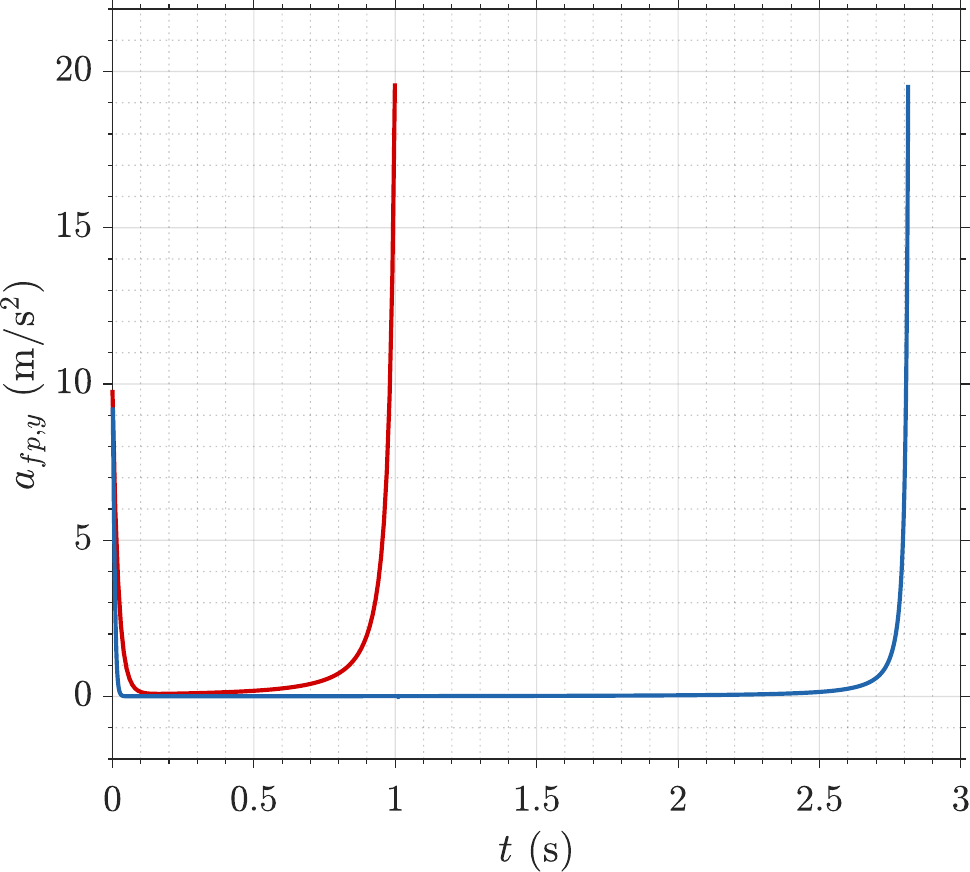}
         \caption{}
         \label{fig:C.3_anal_ay}
     \end{subfigure}
     \hspace{0.1cm}
     \begin{subfigure}[b]{0.32\textwidth}
         \centering
         \includegraphics[width=\textwidth]{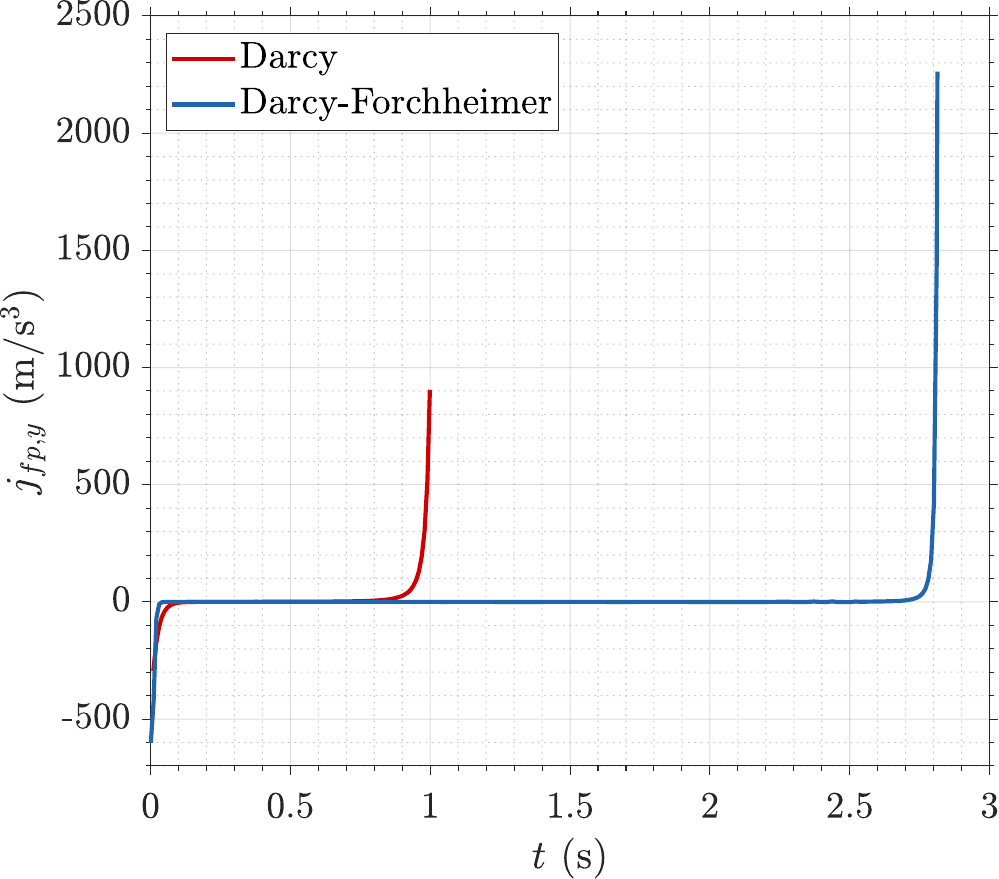}
         \caption{}
         \label{fig:C.3_anal_jy}
     \end{subfigure}
    \caption{1D gravity-driven flow - Case 5: semi-analytical solutions for vertical (a) velocity $v_{fp,y}$, (b) acceleration $a_{fp,y}$, and (c) jerk $j_{fp,y} \equiv \dot{a}_{fp,y}$ of fluid inside the porous domain for different drag force models.}
    \label{fig:C.3_anals}
\end{figure}
\end{remark}

\begin{remark}
In addressing a stiff ODE, it is important to ensure the reliability of the semi-analytical solution. We attempt to verify the behavior of higher-order derivatives obtained through automatic differentiation, particularly near the starting and endpoints, i.e.~at $t=0$ where $h_{fp}=h_0$ and $v_{fp}=0$, as well as when $h_{fp}\rightarrow0$ at $t\rightarrow t^*$, respectively. Here, $t^*$ is the end time where $h_{fp}=0$. To systematically assess correctness, we employ the series solution approach, expressing $h_{fp}(t)$ as $r$-th order time polynomial:
\begin{eqnarray}
    h_{fp}(t) := c_0 + c_1 t + c_2 t^2 + c_3 t^3 + c_4 t^4 + \dots + c_r t^r \,,
    \label{eq:hfp_polynomial}
\end{eqnarray}
where $c_r$ is the $r$-th order coefficient. Substituting the expression given by Eq.~\eqref{eq:hfp_polynomial} into Eq.~\eqref{eq:case5_nonlinear_ode} and solving for the leading coefficients of the polynomial at the two endpoints allows us to determine the limits of the higher-order derivatives as $t$ approaches $0$ and $t^*$. Here, we consider the material parameters previously described in \Cref{subsubsec:1d_case5}. The obtained theoretical limits are then compared to the limits obtained from the semi-analytical solutions for acceleration and jerk, considering both the Darcy and Darcy-Forchheimer drag force models, as presented in \Cref{tab:case5_high_order_limits}. As indicated in the table, the limits obtained from the semi-analytical solution align reasonably well with the theoretical limits. This alignment confirms that the higher-order derivatives are well-behaved functions.

\begin{table}[h!]
\small
\centering
\caption{The limits of acceleration and jerk as they approach the starting and endpoints, $t\rightarrow0$ and $t\rightarrow t^*$, obtained from the semi-analytical solutions and series solution approach (theory). The normalized error is computed as the difference between the two limits normalized by the theoretical limit.}
\label{tab:case5_high_order_limits}
\resizebox{\textwidth}{!}{%
\begin{tabular}{||c|c|c|c|c|c||}
\hline
Drag force model                  & Limiting conditions & Variables & Semi-analytical limit    & Theoretical limit & Normalized error (\%) \\ \hline \hline
\multirow{4}{*}{Darcy}            & \multirow{2}{*}{\begin{tabular}[c]{@{}c@{}}$t=0$ s, \\ $h_{fp}=h_0$,  $v_{fp,y}=0$ m/s\end{tabular}}  & $a_{fp,y}\,{\rm (m/s^2)}$  & 9.81                   & 9.81            & 0                 \\ \cline{3-6} 
                                  &                     & $j_{fp,y}\,{\rm (m/s^3)}$  & -527.4881                & -481.1805         & 9.6237            \\ \cline{2-6} 
                                  & \multirow{2}{*}{\begin{tabular}[c]{@{}c@{}}$t=t^*=0.9988$ s, \\ $h_{fp}=0$ m,  $v_{fp,y}=1.0803$ m/s\end{tabular}}  & $a_{fp,y}\,{\rm (m/s^2)}$  & 19.6197                  & 19.62           & -0.0015           \\ \cline{3-6} 
                                  &                     & $j_{fp,y}\,{\rm (m/s^3)}$  & 905.1797                 & 904.603          & 0.0638            \\ \hline
\multirow{4}{*}{Darcy-Forchheimer} & \multirow{2}{*}{\begin{tabular}[c]{@{}c@{}}$t=0$ s, \\ $h_{fp}=h_0$,  $v_{fp,y}=0$ m/s\end{tabular}}  & $a_{fp,y}\,{\rm (m/s^2)}$  & 9.81                   & 9.81            & 0                 \\ \cline{3-6} 
                                  &                     & $j_{fp,y}\,{\rm (m/s^3)}$  & -528.7478                & -481.1805         & 9.8855            \\ \cline{2-6} 
                                  & \multirow{2}{*}{\begin{tabular}[c]{@{}c@{}}$t=t^*=2.8136$ s, \\ $h_{fp}=0$ m,  $v_{fp,y}=0.5323$ m/s\end{tabular}}  & $a_{fp,y}\,{\rm (m/s^2)}$  & 19.619                  & 19.62           & -0.0053           \\ \cline{3-6} 
                                  &                     & $j_{fp,y}\,{\rm (m/s^3)}$  & 2262.8172                & 2267.358         & -0.2003           \\ \hline
\end{tabular}%
}
\end{table}

\end{remark}

\begin{remark}
    Even though the Newmark-$\beta$ method provides unconditional stability when $2\beta_N \geq \gamma_N\geq \frac{1}{2}$, the accuracy of the MPM solution of this problem remains sensitive to the time increments used. This sensitivity suggests the need and potential benefits of adaptive and higher-order time-stepping strategies, e.g.~higher-order implicit Runge–Kutta methods \cite{solin2012adaptive, butcher2016numerical}, to achieve reliable and efficient numerical solutions.
\end{remark}

\subsection{1D gravity-driven flow - Case 6}
\label{app:1d_case6_derivation}

Following the processes described in \ref{app:1d_case4_derivation} and \ref{app:1d_case5_derivation}, we can further derive an analytical expression for flow across different porous media, with possibly different permeability coefficients. First, let us introduce the assumed notations. Consider two porous media with domain described within volume $\Omega_{pm1}$ and $\Omega_{pm2}$, where the fluid flowing within these porous media are defined as $\Omega_{fp1}$ and $\Omega_{fp2}$, respectively. Here, the governing equations of the fluid flow in the problem described in \Cref{subsubsec:1d_case6} are given as:
\begin{alignat}{3}
    \frac{\pd v_{1,y}}{\pd y}(t)&=0\,, \qquad&  a_{1,y}(t)&= g- \widetilde{A}_1 v_{1,y}(t) - \widetilde{B}_1 |v_{1,y}(t)| v_{1,y}(t)-\frac{p_i(t)}{h_{1}(t) \rho}\,, &\qquad& \forall \tb x \in \Omega_{fp1}\,, \label{eq:pm1_governing_equations} \\
    \frac{\pd v_{2,y}}{\pd y}(t)&=0\,, \qquad&  a_{2,y}(t)&= g- \widetilde{A}_2 v_{2,y}(t) - \widetilde{B}_2 |v_{2,y}(t)| v_{2,y}(t)+\frac{p_i(t)}{h_{2}(t) \rho}\,, &\qquad& \forall \tb x \in \Omega_{fp2}\,, \label{eq:pm2_governing_equations}
\end{alignat}
  where the variables $\square_\alpha$, $\alpha=\{1,2\}$, denote the respective variables associated with the pore fluid $\Omega_{fp\alpha}$ or porous domain $\Omega_{pm\alpha}$  (c.f.~Fig.~\ref{fig:3.1_model}). 

The kinematic and pressure compatibility conditions within the pore fluid interface, also denoted as $\Gamma_i$, are expressed as follows:
\begin{eqnarray}
    \theta_1 \tb v_{1}\cdot \overline{\tb n} = \theta_2 \tb v_{2} \cdot \overline{\tb n} \qquad \mathrm{at}\, \Gamma_i\,,
    \label{eq:velocity_across_jump_3}
\end{eqnarray}
\begin{eqnarray}
    p_{1} = p_{2} := p_i\qquad \mathrm{at}\, \Gamma_i\,,
    \label{eq:pressure_across_jump_3}
\end{eqnarray}
where $0<\theta_1<1$ and $0<\theta_2<1$ are the porosity fields for the two porous domains. Considering the continuity conditions for each pore fluid in 1D, Eq.~\eqref{eq:velocity_across_jump_3} is equivalent to:
\begin{eqnarray}
    \theta_1 v_{1,y}(t) = \theta_2 v_{2,y}(t)\,.
    \label{eq:velocity_across_jump_1d_3}
\end{eqnarray}

The pore fluid is assumed to flow from $\Omega_{pm1}$ to $\Omega_{pm2}$, with initial height $h_0$. The evolution of the heights of the pore fluids is written as:
\begin{eqnarray}
    h_{1}(t) = h_0 - \int_0^t v_{1, y}(t) \td t\,, \qquad
    h_{2}(t) = \int_0^t v_{2, y}(t) \td t\,,
\label{eq:height_change_3}%
\end{eqnarray}
or, in virtue of Eq.~\eqref{eq:velocity_across_jump_1d_3}, can be written as:
\begin{eqnarray}
     h_{1}(t) = h_0 - \frac{\theta_2}{\theta_1} h_{2}(t)\,.
\label{eq:height_change_connected_3}%
\end{eqnarray}
Following the same process as before, the kinematic fields of the pore fluids over time can be related to the height $h_1$ and its time derivatives, i.e.
\begin{subequations}
    \begin{eqnarray}
        v_{1,y}(t) = - \dot{h}_{1}(t)\,, &\qquad& a_{1,y} (t)= - \ddot{h}_{1}(t)\,,\\
        v_{2,y}(t) = - \frac{\theta_1}{\theta_2}\dot{h}_{1}(t)\,, &\qquad& a_{2,y}(t) = - \frac{\theta_1}{\theta_2}\ddot{h}_{1}(t)\,.
    \end{eqnarray}
    \label{eq:kinematic_derivation_1d_3}%
\end{subequations}

Substituting these equations to the governing equations, Eqs.~\eqref{eq:pm1_governing_equations} and \eqref{eq:pm2_governing_equations}, and performing certain mathematical manipulations will yield the following nonlinear ODE to be solved numerically:
\begin{eqnarray}
    h_1 \left(g + \ddot{h}_1 + \widetilde{A}_1\dot{h}_1 + \widetilde{B}_1 |\dot{h}_1|\dot{h}_1 \right) +
    \frac{\theta_1}{\theta_2}\left(h_0 - h_1\right) \left(g + \frac{\theta_1}{\theta_2}\ddot{h}_1 + \frac{\theta_1}{\theta_2}\widetilde{A}_2\dot{h}_1 + \frac{\theta_1^2}{\theta_2^2}\widetilde{B}_2 |\dot{h}_1|\dot{h}_1 \right)
    =0\,.
    \label{eq:case6_nonlinear_ode}
\end{eqnarray}
Upon solving for $h_1$ and its time derivatives, the pressure at the interface $p_i$ can be subsequently computed,
\begin{eqnarray}
    p_i(t) = \rho h_1(t) \left(g + \ddot{h}_1(t) + \widetilde{A}_1\dot{h}_1(t) + \widetilde{B}_1 |\dot{h}_1(t)|\dot{h}_1(t) \right)\,.
\end{eqnarray}
Finally, by considering the material and initial conditions, the semi-analytical results can be obtained as plotted in Fig.~\ref{fig:3.1.6_results} for two different sets of porosity and permeability values.

\subsection{2D steady Couette flow in a composite channel}
\label{app:couette}

The steady-state analytical solution for Couette flow in a composite channel was derived and proposed by \citet{kuznetsov1998analytical}, which considered the non-Darcian effects with Brinkman's and Fochheimer's corrections. In this appendix, the analytical solution of \citet{kuznetsov1998analytical} is rewritten considering the used notation. Considering the geometry defined in Fig.~\ref{fig:3.2_model}, the governing equations of the fully-developed and steady-state flow for both the non-porous and porous regions are given as:
\begin{alignat}{3}
    \frac{\pd v_{ff,x}}{\pd x}&=0\,, \qquad&  \mu \frac{\pd^2 v_{ff,x}}{\pd y^2} &=0 \,, \qquad& -\delta \leq y \leq 0&\,,\label{eq:ff_governing_equations_couette}\\
    \frac{\pd v_{fp,x}}{\pd x}&=0\,, \qquad&  \mu_e \frac{\pd^2 v_{fp,x}}{\pd y^2} - \rho \theta \left( \widetilde{A} + \widetilde{B} |v_{fp,x}|\right) v_{fp,x}&= 0\,, \qquad& 0 \leq y \leq L&\,. \label{eq:pm_governing_equations_couette}
\end{alignat}
Here, the porous interface $\Gamma_i$ is located at $y=0$ and the subscripts $\square_{ff}$ and $\square_{fp}$ denote the fluid variables in the non-porous and porous domain, respectively. 

These governing equations are subject to the following boundary conditions:
\begin{subequations}
    \begin{eqnarray}
        v_{ff,x} = v_\text{bottom} &\quad& \mathrm{at}\,\, y=-\delta\,, \label{eq:bc_couette_analytical_a}\\
        v_{fp,x} = 0 &\quad& \mathrm{at}\,\, y=L\,,\label{eq:bc_couette_analytical_b}\\
        v_{ff,x} = v_{fp,x} = v_{i,x} &\quad& \mathrm{at}\,\, y=0\,,\label{eq:bc_couette_analytical_c}\\
        \mu_e \left.\frac{\pd^2 v_{fp,x}}{\pd y^2}\right|_{y=0^+} - \mu \left.\frac{\pd^2 v_{ff,x}}{\pd y^2}\right|_{y=0^-} = \beta' \frac{\mu}{\sqrt{\Tilde{k}}} v_{i,x} &\quad& \mathrm{at}\,\, y=0\label{eq:bc_couette_analytical_d}\,.
    \end{eqnarray}
    \label{eq:bc_couette_analytical}%
\end{subequations}
Here, Eqs.~\eqref{eq:bc_couette_analytical_a} and \eqref{eq:bc_couette_analytical_b} are the imposed BCs, whereas Eq.~\eqref{eq:bc_couette_analytical_c} dictates the kinematic compatibility of velocity tangent to the interface $\Gamma_i$, where $v_{i,x}$ is the horizontal velocity at the interface. Furthermore, the last condition, Eq.~\eqref{eq:bc_couette_analytical_d}, denotes the tangential stress jump condition, where $\beta'$ is an adjustable coefficient of order one \cite{ochoa1995momentum_i} and $\Tilde{k}=k/\theta^2$ is the intrinsic permeability scaled by porosity-squared. As shown by \citet{ochoa1995momentum_ii}, this tangential stress jump BC is different from the one proposed by \citet{beavers1967boundary}. In the current work, the tangential traction across the interface is considered to be in equilibrium (non-slip), and thus, $\beta'=0$ is considered. 

While deriving the solution to the governing equations, Eqs.~\eqref{eq:ff_governing_equations_couette} and \eqref{eq:pm_governing_equations_couette}, \citet{kuznetsov1998analytical} considered that the momentum boundary layer due to the moving plate at $y=-\delta$ does not reach the fixed plate at $y=L$, and hence, the non-slip velocity condition defined by Eq.~\eqref{eq:bc_couette_analytical_b} can be replaced with the far-field boundary conditions:
\begin{eqnarray}
    v_{fp,x} \rightarrow 0\,,  \quad \frac{\pd v_{fp,x}}{\pd y} \rightarrow 0 \quad \mathrm{as}\,\, y\rightarrow\infty\,.\label{eq:bc_couette_analytical_e}
\end{eqnarray}
This condition is an acceptable assumption as long as the Darcy number is small, i.e.~$Da_H=\Tilde{k}/H^2<1$\cite{kuznetsov1998analytical}. 

Denoting the normalized vertical coordinate as $\Tilde{y}=y/H$ and the normalized velocity as $\Tilde{v}=v/v_\text{bottom}$, the steady-state analytical solution for this problem can be expressed as:
\begin{eqnarray}
\begin{dcases}
    \Tilde{v}_{ff,x} = \Tilde{v}_{i,x} - \left(1-\Tilde{v}_{i,x}\right)\ddfrac{H}{\delta} \Tilde{y}\,, & -\ddfrac{\delta}{H} \leq \Tilde{y} \leq 0\,,\\
    \Tilde{v}_{fp,x} = \ddfrac{3}{2} \ddfrac{A'}{B'} \left\lbrace \left[ \ddfrac{D' \exp\left( \ddfrac{\sqrt{B'}}{\gamma_e}\Tilde{y} \right)-1}{D' \exp\left( \ddfrac{\sqrt{B'}}{\gamma_e}\Tilde{y} \right)+1} \right]^2 - 1  \right\rbrace \,, & 0 \leq \Tilde{y} \leq \ddfrac{L}{H}\,,
\end{dcases}
\end{eqnarray}
where
\begin{eqnarray}
    A' = \frac{Re_H F}{\sqrt{Da_H}}\,,\qquad B' = \frac{1}{Da_H}\,, \qquad D' = \frac{1+\left( 1 + \ddfrac{2}{3} \ddfrac{A'}{B'}\Tilde{v}_{i,x} \right)^{1/2}}{1-\left( 1 + \ddfrac{2}{3} \ddfrac{A'}{B'}\Tilde{v}_{i,x} \right)^{1/2}}\,.
    \label{eq:couette_dimensional}
\end{eqnarray}
Here, the dimensionless velocity at the porous interface $\Gamma_i$,  i.e.~$\Tilde{v}_{i,x}$, can be obtained by solving the following transcendental equation:
\begin{eqnarray}
    -\gamma_e \Tilde{v}_{i,x} \left( \frac{2}{3} A' \Tilde{v}_{i,x} + B' \right)^{1/2} + \left(1-\Tilde{v}_{i,x}\right) \frac{H}{\delta} = \beta' \sqrt{B'} \Tilde{v}_{i,x}\,.
\end{eqnarray}
In the expressions above, additional nondimensional parameters are introduced, which are defined as follows:
\begin{eqnarray}
    \gamma_e=\left(\frac{\mu_e}{\mu}\right)^{1/2}\,, \qquad Re_H = \frac{\rho v_\text{bottom} H}{\mu} \,, \qquad F = \frac{B\sqrt{\theta}}{\sqrt{A}}\,.
    \label{eq:couette_nondimensional}
\end{eqnarray}
Notice that the parameters $A$ and $B$ in Eq.~\eqref{eq:couette_nondimensional}$_3$ are the nonlinear drag force coefficients (see \ref{app:drag_force_coeff}), which should be distinguished from $A'$ and $B'$ defined in Eqs.~\eqref{eq:couette_dimensional}$_1$ and \eqref{eq:couette_dimensional}$_2$, respectively.

\section{Fractional-step approach}
\label{app:fractional_step}

In the fractional-step approach, the velocity and pressure fields are weakly coupled through an intermediate state of velocity denoted as $\tb v^*\equiv \dot{\tb u}^*$. The strong-form unified momentum balance (Eq.~\eqref{eq:equilibrium}) can be rewritten following the pressure projection proposed by \citet{chorin1968numerical} as:
\begin{eqnarray}
    \theta^n \rho \frac{\dot{\tb u}^* - \dot{\tb u}^{n}}{\Delta t} &=&  \nabla\cdot \tb s^n + \theta^n \rho \tb b^n - \theta^n \rho \left( \widetilde{A} + \widetilde{B} ||\dot{\tb u}^{n}|| \right) \dot{\tb u}^{n} \,,\label{eq:predictor} \\
    \theta^n \rho \frac{\dot{\tb u}^{n+1} - \dot{\tb u}^{*}}{\Delta t} &=& - \theta^n \nabla p^{n+1} \,, \label{eq:corrector} \\
     \nabla \cdot \left( \theta^n \nabla p^{n+1} \right) &=& \frac{\rho}{\Delta t}\nabla \cdot \left( \theta^n \dot{\tb u}^* \right)\,.
 \label{eq:PPE}
\end{eqnarray}
Here, we follow the explicit drag force assumption\cite{yamaguchi2020solid}, and hence, the critical time step is acknowledged to be dependent on the value of permeability coefficients (see \citet{kularathna2021semi}). 

The corresponding weak form of Eqs.~\eqref{eq:predictor}-\eqref{eq:PPE} can be derived by including the two test functions, $\delta \tb u$ and $\delta p$, for the two primary variables. After applying Gauss' theorem and performing several manipulations, the weak-form predictor-corrector scheme can be obtained:
\begin{linenomath}
\begin{align}
\begin{split}
    \int_\Omega \theta^n \rho \frac{\dot{\tb u}^* - \dot{\tb u}^{n}}{\Delta t} \cdot \delta \tb u \,\td\Omega &=  -\int_\Omega \tb s^n : \nabla\delta \tb u \,\td\Omega + \int_\Omega \theta^n {\rho} \tb b^n \cdot \delta \tb u \,\td\Omega  \\
    & \qquad - \int_{\Omega} \theta^n \rho \left( \widetilde{A} + \widetilde{B}||\dot{\tb u}^n||
    \right) \dot{\tb u}^n \cdot \delta \tb u \,\td \Omega + \int_{\Gamma_{N}}\overline{\tb t}_s \cdot \delta  \tb u\,\td\Gamma\,,\label{eq:inc_predictora_weak}
\end{split} \\
      \int_\Omega \theta^n \rho \frac{\dot{\tb u}^{n+1} - \dot{\tb u}^{*}}{\Delta t} \cdot \delta \tb u \,\td\Omega &= - \int_\Omega \theta^n \nabla p^{n+1} \cdot \delta \tb u \,\td \Omega \,,\label{eq:inc_correctora_weak}\\
   \int_\Omega \theta^n \nabla p^{n+1} \cdot \nabla \delta p \,\td \Omega &=  -\int_\Omega \frac{\rho}{\Delta t} \nabla \cdot \left(\theta^n \dot{\tb u}^* \right) \delta p \,\td \Omega + \int_{\Gamma_{N}}\overline{t}_p \delta p\,\td\Gamma\,. \label{eq:PPE_weak}
\end{align}
\end{linenomath}
Here, $\bar{\tb t}_s = \tb s^n \cdot \tb n$ and $ \bar t_p = \theta^n \nabla p^{n+1} \cdot \tb n$ denote the surface deviatoric traction and the pressure gradient traction, respectively. Similar to the previous assumptions, these terms will be omitted.

Taking into account the spatial discretization assumption described in \Cref{subsec:spatial_disc}, the blurred porosity interface in \Cref{subsec:porosity_field}, and by invoking the arbitrariness of the nodal test functions, $\delta \tb u_{I}$ and  $\delta p_{I}$, Eqs.~\eqref{eq:inc_predictora_weak}-\eqref{eq:PPE_weak} can be expressed in their fully discrete form by the following series of equations\cite{yamaguchi2020solid, kularathna2021semi}:
\begin{enumerate}
    \item \textit{Predictor}: to compute intermediate nodal acceleration and velocity, $\ddot{\tb u}^*_I$ and $\dot{\tb u}^*_I$, i.e.
    \begin{eqnarray}
        \ddot{\tb u}^*_I &=& \frac{1}{m^n_I}\left(\tb f_I^{\rm int} + \tb f_I^{\rm ext} + \tb f_I^{\rm drag}\right)\,, \label{eq:predictor_mpm}\\
        \dot{\tb u}^*_I &=& \dot{\tb u}^n_I + \Delta t \ddot{\tb u}^*_I \label{eq:predictor_mpm_vel}\,,
    \end{eqnarray}
    where the nodal internal, external, and drag forces can be expressed as:
    \begin{subequations}
    \begin{eqnarray}
        \tb f_I^{\rm int} &=& - \sum_{p=1}^{n_p} \tb B_{Ip}^T \tb{s}^n_p V^n_p\,, \\
        \tb f_I^{\rm ext} &=& \sum_{p=1}^{n_p} m_p S_{Ip} \tb b^n_p \,, \\ 
        \tb f_I^{\rm drag} &=& - \Biggl( \sum_{p=1}^{n_{p}} \theta^n_p \rho_p \left( \widetilde{A}_p + \widetilde{B}_p \left\Vert \dot{\tb u}^n_p \right \Vert
\right) S_{Ip} V_p^n\Biggl) \dot{\tb u}^n_I\,.
    \end{eqnarray}
    \label{eq:ext_int_drag_forces}%
    \end{subequations}
        
    \item \textit{Pressure Poisson equation}: to solve for the updated nodal pressure $p^{n+1}_I$, i.e.
    \begin{eqnarray}
        \tb L_{IJ} {p}^{n+1}_J = \tb D_{IJ} \dot{\tb u}^{*}_J \,.
        \label{eq:ppe_matrix}
    \end{eqnarray}
    Here, the Laplacian matrix $\tb L$ and the divergence matrix $\tb D$ are expressed as,
    \begin{eqnarray}
     \tb L_{IJ} = \sum_{p=1}^{n_{p}} \theta_p^n \nabla S_{Ip} \nabla S_{Jp}  V^n_p \,, 
     \qquad
     \tb D_{IJ} = - \frac{{\rho_p}}{{\Delta t}} \sum_{p=1}^{n_{p}} \left( \theta^n_p {S_{Ip} \, \nabla S_{Jp}} + \nabla \theta^n_p {S_{Ip} \, S_{Jp}} \right) \, V^n_p \,. \label{eq:divergence_matrix_mpm}
    \end{eqnarray}
    In order to obtain a unique solution of Eq.~\eqref{eq:ppe_matrix}, Dirichlet boundary conditions should be properly imposed, particularly at the free surface nodes. The current work utilizes the volume fraction detection approach proposed by \citet{kularathna2018splitting}, which has been further extended to account for higher order basis function \cite{chandra2023stabilized} and varying porosity field. Further details on the detection algorithm, boundary imposition, and associated numerical errors in MPM can be found in \cite{kularathna2018splitting, chandra2023stabilized}.
    
    \item \textit{Corrector}: to compute the updated acceleration and velocity, $\ddot{\tb u}_I^{n+1}$ and $\dot{\tb u}_I^{n+1}$, i.e.
    \begin{eqnarray}
        \dot{\tb u}^{n+1}_I = \dot{\tb u}^{*}_I + \frac{\Delta t}{{m}^n_I} \tb G_{IJ} {p}^{n+1}_J \,,\label{eq:corrector_mpm} 
    \end{eqnarray}
    where, the gradient matrix $\tb G_{IJ}$ can be expressed as,
    \begin{eqnarray}
    \tb G_{IJ} = - \sum_{p=1}^{n_{p}} \theta^n_p S_{Ip} \nabla S_{Jp} \, V^n_p
     \,.\label{eq:gradient_matrix_mpm}
    \end{eqnarray}
    A new nodal acceleration based on the corrected velocity can be obtained as:
    \begin{eqnarray}
        \ddot{\tb u}^{n+1}_I = \frac{\dot{\tb u}^{n+1}_I - \dot{\tb u}^{n}_I}{\Delta t} =  \frac{\dot{\tb u}^{*}_I - \dot{\tb u}^{n}_I}{\Delta t} + \frac{\dot{\tb u}^{n+1}_I - \dot{\tb u}^{*}_I}{\Delta t} = \ddot{\tb u}^{*}_I + \frac{1}{m_I^n}\tb G_{IJ} {p}^{n+1}_J\,.\label{eq:corrector_mpm_acc} 
    \end{eqnarray}
\end{enumerate}

\bibliography{mybibfile}

\end{document}